\documentclass[11pt, oneside]{article}
\usepackage{geometry}

\usepackage{amssymb}
\usepackage{amsthm}
\usepackage{amsmath}
\usepackage{mathtools}
\usepackage{bm}
\usepackage{mathrsfs} 
\usepackage{enumerate}
\usepackage{float} 
\usepackage{tikz-cd} 
\usepackage{todonotes}

\usepackage{xcolor}
\definecolor{ForestGreen}{RGB}{34, 139, 34}
\definecolor{frenchrose}{rgb}{0.96, 0.29, 0.54}
\definecolor{ballblue}{rgb}{0.13, 0.67, 0.8}
\definecolor{Crimson}{rgb}{0.86, 0.08, 0.24}
\definecolor{Plum}{rgb}{0.56, 0.27, 0.6}

\usepackage{tikz}

\usepackage{hyperref}
 \hypersetup{
   colorlinks   = true, 
   linkcolor    =  Crimson, 
   citecolor   = ForestGreen 
 }

\usepackage{relsize}
\usetikzlibrary{shapes.misc}
\usetikzlibrary{calc}

\newtheorem{theorem}{Theorem}[section]
\newtheorem{lemma}[theorem]{Lemma}
\newtheorem{proposition}[theorem]{Proposition}
\newtheorem{corollary}[theorem]{Corollary}
\newtheorem{conjecture}[theorem]{Conjecture}

\theoremstyle{definition}
\newtheorem{definition}[theorem]{Definition}
\newtheorem{example}[theorem]{Example}
\newtheorem{setup}[theorem]{Setup}
\newtheorem{dfnprop}[theorem]{Definition/Proposition}
\newtheorem{remark}[theorem]{Remark}

\newtheorem*{notation*}{Notation}

\DeclareMathOperator{\C}{{\mathbb{C}}} 
\newcommand{\CR}{\mathbb{R}}
\newcommand{\Z}{\mathbb{Z}}
\newcommand{\kk}{\mathbf{k}}
\newcommand{\D}{\mathbb{D}}
\newcommand{\F}{\mathcal{F}}
\newcommand{\W}{\mathcal{W}}
\newcommand{\Fuk}{\textnormal{Fuk}}
\newcommand{\Hom}{\textnormal{Hom}}
\newcommand{\Ext}{\textnormal{Ext}}
\newcommand{\End}{\textnormal{End}}
\DeclareMathOperator{\Mod}{{\mathbf{mod}}} 
\newcommand{\perf}{\textnormal{\textbf{perf}}}
\DeclareMathOperator{\Proj}{{\textnormal{Proj}}} 
\newcommand{\Sym}{\textnormal{Sym}}
\newcommand{\bpi}{\bm{\pi}}
\DeclareMathOperator{\Real}{{\textnormal{Re}}} 
\DeclareMathOperator{\Imaginary}{{\textnormal{Im}}} 

\DeclareMathOperator{\gldim}{{\mathrm{gl.\!dim}}}
\DeclareMathOperator{\domdim}{{\mathrm{dom.\!dim}}}

\renewcommand\S{\mathfrak{S}} 
\renewcommand\L{\mathcal{L}} 

\newcommand{\Ihat}{\hat{\mathcal{I}}} 
\newcommand{\I}{\mathcal{I}} 
\newcommand{\N}{\mathcal{N}} 

\renewcommand{\P}{P} 
\newcommand{\CTmod}{M} 
\newcommand{\op}{\text{op}}

\newcommand{\fhat}{\hat{f}} 
\newcommand{\f}{f} 
\newcommand{\p}{\mathbf{p}} 
\newcommand{\what}{\hat{\mathbf{w}}} 
\newcommand{\w}{\mathbf{w}} 
\newcommand{\Ghat}{\hat{B}} 
            \newcommand{\Bhat}{\Ghat}
\newcommand{\A}{A} 
\newcommand{\G}{B} 
            \newcommand{\B}{\G}
\newcommand{\ghat}{\hat{\lambda}} 

\newcommand{\whatcirc}{\hat{\mathbf{w}}^{\circ}} 
\newcommand{\SODhat}{\hat{\mathcal{S}}} 
\newcommand{\SODhati}{\hat{S}} 
\newcommand{\SOD}{\mathcal{S}} 
\newcommand{\SODi}{S} 

\newcommand{\BFShat}{\hat{\mathcal{B}}} 
\newcommand{\BFS}{\mathcal{B}} 
\newcommand{\AW}{\mathscr{A}} 
\newcommand{\BW}{\mathscr{B}} 

\newcommand{\Dhat}{\hat{D}} 
\newcommand{\Dt}{D} 
\newcommand{\M}{F} 
\newcommand{\Cone}{\textnormal{Cone}}
\newcommand{\Core}{\textnormal{Core}}
\newcommand{\lambdabold}{\bm{\lambda}}
\newcommand{\lambdaboldhat}{\bm{\hat{\lambda}}}
\newcommand{\gammabold}{\bm{\gamma}}
\newcommand{\Dbold}{\bm{D}}
\newcommand{\Vbold}{\bm{V}}
\newcommand{\Deltabold}{\bm{\Delta}}
\newcommand{\Upsilonbold}{\bm{\Upsilon}}

\newcommand{\X}{X} 
\newcommand{\HH}{\mathcal{H}} 

\newcommand{\K}{K} 

\DeclarePairedDelimiter\floor{\lfloor}{\rfloor}

\newcommand{\mycomment}[1]{}

\title{\bf{{Symplectic higher Auslander correspondence for type A}}}
\date{}		
\author{Ilaria Di Dedda}

\begin{document}


\maketitle

\textbf{Abstract}. We prove that the Fukaya--Seidel categories of a certain family of singularities on $\C^d$ are equivalent to the perfect derived categories of higher Auslander algebras of Dynkin type A. We relate these to the Fukaya--Seidel categories of Brieskorn--Pham singularities and to the partially wrapped Fukaya categories considered in \cite{DJL}. We provide a symplectic interpretation to higher Auslander correspondence of type A in terms of Fukaya\textendash Seidel categories of Lefschetz fibrations.

\setcounter{tocdepth}{2}
{\hypersetup{
    colorlinks   = true, 
    linkcolor    =  black, 
  }
\tableofcontents
}

\section{Introduction}
\subsection{Main results}\label{sec:mainresults}
A polynomial function ${\p:\C^d \to \C}$ is said to be a \emph{Brieskorn\textendash Pham} polynomial if it is of the form
\begin{equation*}
    \p(x_1,\dots,x_d)=x_1^{n_1+1}+\dots+x_d^{n_d+1}
\end{equation*}
for a tuple $(n_1,\dots,n_d)$ of positive integers. A complex-valued holomorphic function on $\C^d$ with an isolated singularity at the origin is said to be a \emph{Brieskorn--Pham singularity} if its germ is equivalent to some Brieskorn--Pham polynomial.

A natural invariant of any isolated hypersurface singularity $f$ is the derived Fukaya--Seidel category $\F(f)$ constructed in \cite[Chapter~3]{SeidelBk}. The Fukaya\textendash Seidel category of a Brieskorn\textendash Pham singularity is well-understood \cite[Theorem~1.1]{FU11}; we adapt the known construction due to Futaki\textendash Ueda, and use a Sebastiani\textendash Thom-type technique to re-prove the following.

\begin{proposition}[Theorem~\ref{thm:mainthm1}]\label{thm:introthm1}
    The Fukaya\textendash Seidel category of the Brieskorn\textendash Pham singularity $\fhat_{n,d}$ defined by
    \begin{equation*}
        \fhat_{n,d}(x_1,\dots,x_d):= x_1^{n+1} +\dots + x_d^{n+1}
    \end{equation*}
    is quasi-equivalent to the perfect derived category $\perf(\Bhat_{n,d})$ of the algebra $\Bhat_{n,d}:=(\Bhat_{n})^{\otimes d}$, where $\Bhat_{n}$ is the path algebra of the $A_{n}$ quiver over a field $\kk$ with alternating orientation of the arrows, as in Figure~\ref{fig:Bhat} (left).
\end{proposition}

The $\kk$-algebra $\Bhat_{n,d}$ can equivalently be defined as the path algebra over $\kk$ of a quiver, whose set of vertices is indexed by a tuple ${(i_1,\dots,i_d)}$ of integers between $1$ and $n$, with arrows ${(i_1,\dots,i_h,\dots,i_d) \to (i_1,\dots,i_h\pm 1,\dots,i_d)}$ whenever $i_h$ is odd, and which is equipped with all possible commutativity relations. See Figure~\ref{fig:Bhat} (centre), and Definition~\ref{def:Bhatnd}.

\begin{figure}[t]
    \centering
    \begin{minipage}{.3 \textwidth}
        \centering
        \begin{tikzpicture}
            [scale=.8,align=center, v/.style={draw,shape=circle, fill=black, minimum size=1.2mm, inner sep=0pt, outer sep=0pt},
            font=\small,
            ]

            \node (0) at (0,0) {$1$};
            \node (1) at (1,0) {$2$};
            \node (2) at (2,0) {$3$};
            \node (3) at (3,0) {$4$};
 
            \path
            (0) edge[color=black,->] (1)
            (2) edge[color=black,->] (1)
            (2) edge[color=black,->] (3);
        \end{tikzpicture}
    \end{minipage}
    \begin{minipage}{.3 \textwidth}
        \centering
        \begin{tikzpicture}
            [scale=1,align=center, v/.style={draw,shape=circle, fill=black, minimum size=1.2mm, inner sep=0pt, outer sep=0pt},
            font=\small,
            ]
       
            \node (00) at (0,0) {$11$};
            \node (10) at (1,0) {$21$};
            \node (20) at (2,0) {$31$};
            \node (30) at (3,0) {$41$};

            \node (01) at (0,1) {$12$};
            \node (11) at (1,1) {$22$};
            \node (21) at (2,1) {$32$};
            \node (31) at (3,1) {$42$};

            \node (02) at (0,2) {$13$};
            \node (12) at (1,2) {$23$};
            \node (22) at (2,2) {$33$};
            \node (32) at (3,2) {$43$};

            \node (03) at (0,3) {$14$};
            \node (13) at (1,3) {$24$};
            \node (23) at (2,3) {$34$};
            \node (33) at (3,3) {$44$};

            \path
            (00) edge[color=black,->] (10)
            (20) edge[color=black,->] (10)
            (20) edge[color=black,->] (30)

            (01) edge[color=black,->] (11)
            (21) edge[color=black,->] (11)
            (21) edge[color=black,->] (31)

            (02) edge[color=black,->] (12)
            (22) edge[color=black,->] (12)
            (22) edge[color=black,->] (32)

            (03) edge[color=black,->] (13)
            (23) edge[color=black,->] (13)
            (23) edge[color=black,->] (33)

            (00) edge[color=black,->] (01)
            (02) edge[color=black,->] (01)
            (02) edge[color=black,->] (03)

            (10) edge[color=black,->] (11)
            (12) edge[color=black,->] (11)
            (12) edge[color=black,->] (13)

            (20) edge[color=black,->] (21)
            (22) edge[color=black,->] (21)
            (22) edge[color=black,->] (23)

            (30) edge[color=black,->] (31)
            (32) edge[color=black,->] (31)
            (32) edge[color=black,->] (33)

            ;
        \end{tikzpicture}
        \end{minipage}
        \begin{minipage}{.3 \textwidth}
            \centering
            \begin{tikzpicture}
                [scale=1,align=center, v/.style={draw,shape=circle, fill=black, minimum size=1.2mm, inner sep=0pt, outer sep=0pt},
                font=\small,
                ]

                \node (01) at (0,1) {$12$};           
    
                \node (02) at (0,2) {$13$};
                \node (12) at (1,2) {$23$};
             
                \node (03) at (0,3) {$14$};
                \node (13) at (1,3) {$24$};
                \node (23) at (2,3) {$34$};

                \node at (2,0) {};

                \path
           
                (02) edge[color=black,->] (12)
                (03) edge[color=black,->] (13)
                (02) edge[color=black,->] (01)
                (02) edge[color=black,->] (03)
                (12) edge[color=black,->] (13)
                (23) edge[color=black,->] (13)
                ;
            \end{tikzpicture}
            \end{minipage}
\caption{The quivers underlying $\Bhat_{n}$ (left), $\Bhat_{n,d}$ (centre) and $\B_{n,d}$ (right), with all possible commutativity relations, for $n=4$ and $d=2$.}
\label{fig:Bhat}
\end{figure}

When ${n:=n_1=\dots=n_d}$, a Brieskorn\textendash Pham polynomial is invariant with respect to the action of the symmetric group $\S_d$ on $\C^d$ given by permuting its coordinates. Let ${\Sym^d(\C)=\C^d / \S_d}$ be the $d$-fold symmetric product of $\C$, and 
\begin{equation*}
    \bpi: \C^d \to \Sym^d(\C), \qquad (x_1,\dots, x_d) \mapsto \{x_1,\dots, x_d\}
\end{equation*}
its branched cover. Let $\f_{n,d}$ be the map
\begin{equation*}
    \f_{n,d}: \Sym^d(\C)\to \C, \qquad\quad \{x_1,\dots,x_d\} \mapsto x_1^{n+1} +\dots + x_d^{n+1},
\end{equation*}
which satisfies
\begin{equation}\label{eq:liftfnd}
    \fhat_{n,d}=\f_{n,d} \circ \bpi.
\end{equation}
Under the identification ${\varphi:\Sym^d(\C)\xrightarrow{\cong}\C^d}$, which we discuss in Section~\ref{sec:FSbackground}, and for $n\geq d$, $\f_{n,d}\circ \varphi^{-1}$ defines an isolated hypersurface singularity. We define $\B_{n,d}$ to be the path algebra over $\kk$ of the quiver whose set of vertices is indexed by a tuple ${(i_1,\dots,i_d)}$ of strictly increasing integers between $1$ and $n$, with arrows ${(i_1,\dots,i_h,\dots,i_d) \to (i_1,\dots,i_h\pm 1,\dots,i_d)}$ whenever $i_h$ is odd, and which is equipped with all possible commutativity relations. See Figure~\ref{fig:Bhat} (right), and Definition~\ref{def:Bnd}. Our first main result concerns the computation of ${\F(\f_{n,d}):=\F(\f_{n,d}\circ \varphi^{-1})}$.

\begin{theorem}[Theorem~\ref{thm:mainthm2}]\label{thm:introthm2}
    The Fukaya\textendash Seidel category of $\f_{n,d}$ is quasi-equivalent to the perfect derived category $\perf(\B_{n,d})$.
\end{theorem}

The relation between Proposition~\ref{thm:introthm1} and Theorem~\ref{thm:introthm2} is the following. We prove in Section~\ref{section:2} (specifically, a combination of Theorem~\ref{thm:mainthm1}, Proposition~\ref{prop:actiononthimbles}, and Theorem~\ref{thm:mainthm2}) that there is a strong full exceptional collection $\SODhat$ in $\F(\fhat_{n,d})$ such that:
\begin{enumerate}[i.]
    \item The endomorphism algebra of $\SODhat$ is isomorphic to $\Bhat_{n,d}$ as a graded $\kk$-algebra;
    \item The action of $\S_d$ on $\C^d$ restricts to the collection of geometric objects underlying $\SODhat$;
    \item On the subcollection  $\SODhat^{\circ} \subset \SODhat$ on which the $\S_d$-action is free, there is a choice of geometric object $\Dt_i$ ($i=1,\dots,{n \choose d}$), representative of each equivalence class in $\SODhat^{\circ}/\S_d$, for which $\SOD:=\{\Dt_i\}_i$ is a strong full exceptional collection in $\F(\f_{n,d})$;
    \item The derived endomorphism algebra of $\SOD$ is concentrated in degree 0 and is isomorphic to $\B_{n,d}$.
\end{enumerate}
This is discussed in Sections~\ref{section:BP} and \ref{section:quotientFC}.

Following \cite{Auroux2}, we equip $\Sym^d(\D)$ (a ball of real dimension $2d$) with the structure of a Liouville domain, and with a collection of stops in the symplectic hypersurfaces:
\begin{equation}\label{eq:DJLstops}
    \Lambda_n^{(d)}:= \bigcup_{p \in \Lambda_{n}} \{p\} \mathrel{\boldsymbol{\times}} \Sym^{d-1} (\D),
\end{equation}
where ${\Lambda_n \subset \partial \D}$ is a set of $n+1$ marked points, and each point ${\left(p, \{x_1,\dots,x_{d-1}\}\right)}$ on the right-hand side of \eqref{eq:DJLstops} is identified with a point in $\Sym^{d} (\D)$ by
\begin{equation*}
    \left(p, \{x_1,\dots,x_{d-1}\}\right) \mapsto  \{p, x_1,\dots,x_{d-1}\}.
\end{equation*}
We consider the \emph{partially wrapped Fukaya category} defined by Auroux in \cite{Auroux1} of the pair $(\Sym^d(\D),\Lambda_n^{(d)})$, which we denote by $\W_n^d:=\W^{\textnormal{Aur}}(\Sym^d(\D),\Lambda_n^{(d)})$. Our second result is the following.

\begin{theorem}[Theorem~\ref{thm:mainthm3}]\label{thm:introthm3}
    There is a quasi-equivalence of triangulated $A_{\infty}$-categories
    \begin{equation*}
        \F(\f_{n,d}) \xrightarrow{\text{ $\simeq$ }} \W_n^{d}.
    \end{equation*}
\end{theorem}

Theorem~\ref{thm:introthm3} amounts to an explicit computation of the derived endomorphism algebra of a preferred set of generators of $\W_n^{d}$. Our computations provide an isomorphism between this algebra and the $\kk$-algebra $\B_{n,d}$ given by Theorem~\ref{thm:introthm2}; this is discussed in Section~\ref{section:F-W}.

The main result in \cite{DJL} provides a quasi-equivalence of triangulated $A_{\infty}$-categories between $\W_n^d$ and the perfect derived $A_{\infty}$-category of the \emph{higher Auslander algebra} $\A_{n,d}$; we can define the latter as follows. Given $n$ and $d$ positive integers satisfying $n\geq d$, consider the quiver whose set of vertices is indexed by a tuple $(i_1,\dots,i_d)$ of strictly increasing integers between $1$ and $n$, and arrows given by
\begin{equation*}
    (i_1,\dots,i_h,\dots,i_d) \to (i_1,\dots,i_h+1,\dots,i_d)
\end{equation*}
for all $h=1,\dots,d$. We define $\A_{n,d}$ to be the path algebra of this quiver over $\kk$, modulo certain relations, which we will detail in Definition~\ref{def:And}. Equivalently (\cite[Theorem~6.12]{Iyama2}), each $\kk$-algebra $\A_{n,d}$ can be inductively constructed from $\A_{n-1,d-1}$ in the following way. $\A_{n,1}$ denotes the path algebra of the linearly oriented $A_n$ quiver. By work of Iyama (\cite[Section~6]{Iyama2}), the iteratively constructed algebra $\A_{n-1,d-1}$ admits a unique \emph{$(d-1)$-cluster tilting module} $M$, that is a (left) $\A_{n-1,d-1}$-module such that:
\begin{align*}
    \textnormal{add}M=&\{X \in \Mod \A_{n-1,d-1} \mid  \Ext_{\A_{n-1,d-1}}^i(X,M)=0 \quad \text{for any } 0 < i <d-1\}= \\
    =&\{X \in \Mod \A_{n-1,d-1} \mid  \Ext_{\A_{n-1,d-1}}^i(M,X)=0 \quad \text{for any } 0 < i <d-1\},
\end{align*}
where $\textnormal{add}\mathcal{X}$ denotes the smallest full subcategory of $\Mod \mathcal{X}$ which is idempotent complete, closed under finite directed sums, and contains $\mathcal{X}$. $\A_{n,d}$ is defined as:
\begin{equation}\label{eq:higherAuslA}
    \A_{n,d}:=\End_{\A_{n-1,d-1}}(M).
\end{equation}

Immediate consequence of \cite[Theorem~1]{DJL} and Theorem~\ref{thm:introthm3} is the following result.

\begin{corollary}\label{cor:maincor}
    There is a quasi-equivalence of triangulated $A_{\infty}$-categories
\begin{equation*}
    \pushQED{\qed} \F(\f_{n,d}) \xrightarrow{\text{ $\simeq$ }}  \perf(\A_{n,d}).\qedhere \popQED
\end{equation*}
\end{corollary}

As an application of Corollary~\ref{cor:maincor}, we provide in Section~\ref{sec:tiltingobject} a \emph{tilting object} for each $\A_{n,d}$ (i.e.\ an object $T$ generating $\perf(\A_{n,d})$ as a triangulated category, and satisfying $\Hom(T,T[i])=0$ for $i\neq 0$), which is not isomorphic to $\A_{n,d}$ itself. Tilting objects are of special interest in \emph{tilting theory}, as they constitute a rich source of derived equivalences, due to a theorem of Happel \cite[Chapter~III]{Happel}, later generalised by Rickard \cite[Theorem~6.4]{Rickard} (see also \cite{Keller}). We interpret the derived equivalence provided by this tilting object as the higher-dimensional analogue of the well-known  derived equivalence between the linearly oriented $A_n$ quiver and the alternating one (see for example a combination of \cite[Theorem~3.2]{KY11} and \cite[Theorem~8.6]{FZ03}).

Finally, we formulate what we call \emph{symplectic higher Auslander correspondence of type A}, which provides a geometric motivation for the algebraic definition \eqref{eq:higherAuslA}. We present this as an inductive construction that geometrically mimics a fundamental result in representation theory, which we will contextualise in Section~\ref{sec:Auslbackground}.

\begin{theorem}[Theorem~\ref{prop:uniquetopfibr}, Corollary~\ref{mainthmsympAus}, Definition~\ref{def:Sigma}]\label{thm:introthm4}
    For a fixed pair of integers $n\geq d\geq 1$, suppose there is a symplectic Lefschetz fibration $\Phi_{n,d}: \X_{n,d} \to \D$ on a symplectic manifold with boundary $\X_{n,d}$, and a quasi-equivalence of triangulated $A_{\infty}$-categories
    \begin{equation*}
        \F(\Phi_{n,d}) \xrightarrow{\text{ $\simeq$ }}  \perf(\A_{n,d}).
    \end{equation*}
    Then, there is an algorithmically prescribed Lefschetz fibration $\Phi_{n+1,d+1}:\X_{n+1,d+1} \to \D$ on a symplectic manifold $\X_{n+1,d+1}$, and a quasi-equivalence of triangulated $A_{\infty}$-categories:
    \begin{equation*}
        \F(\Phi_{n+1,d+1}) \xrightarrow{\text{ $\simeq$ }} \perf(\A_{n+1,d+1}).
    \end{equation*}
    Moreover, the generic fibre of $\Phi_{n+1,d+1}$ is algorithmically prescribed via a sequence of symplectic handle attachments to $\X_{n,d}$.
\end{theorem}

The relation between Theorem~\ref{thm:introthm4} and the higher Auslander correspondence is one of the main results of this article, and it will be explained in Section~\ref{section:3}.

The identification of $\perf(\A_{n,d})$ with $\F(\Phi_{n,d})$ is a potential source of derived equivalences for the $\kk$-algebras $\A_{n,d}$, arising from the monodromy of $\Phi_{n,d}$. We also expect, but do not prove at this time, the construction we describe in Section~\ref{section:3} to generalise to a wider class of $\kk$-algebras (i.e.\ beyond the type A case), for which currently there is no geometric counterpart. This falls into the more general framework of finding connections between representation theory and geometry, with the goal of using geometric techniques to prove algebraic results, which for example \cite{DJL} do in the type A context. In this sense, the generalisability of Theorem~\ref{thm:introthm4} makes this construction more versatile than the one in \cite[Theorem~1]{DJL}.

\begin{notation*}
    Throughout the article, we use the following notation:
    \begin{itemize}
        \item[--] $\A_{n,d}$ and $\AW_{n,d}$ denote two quasi-isomorphic $A_{\infty}$-algebras. They are endomorphism algebras computed in $\perf(\A_{n,d})$ and $\W_n^d$ respectively, defined in Definition~\ref{def:And} and Section~\ref{sec:tiltingobject};
        \item[--] $\B_{n,d}$, $\BFS_{n,d}$ $\BW_{n,d}$ denote three pairwise isomorphic $A_{\infty}$-algebras. They and are endomorphism algebras computed in $\perf(\B_{n,d})$, $\F(\f_{n,d})$ and $\W_n^d$ respectively, defined in Definitions~\ref{def:Bnd}, \ref{def:BSFnd} and \ref{def:BWnd};
        \item[--] $\Bhat_{n,d}$ and $\BFShat_{n,d}$ denote two isomorphic $A_{\infty}$-algebras. They and are endomorphism algebras computed in $\perf(\Bhat_{n,d})$ and $\F(\fhat_{n,d})$ respectively, defined in Definitions~\ref{def:Bhatnd} and \ref{def:BFShatnd};
        \item[--] $\N_{n,d}$, $\I_{n,d}$ and $\Ihat_{n,d}$ denote three different posets, defined in Definitions~\ref{def:Nnd}, \ref{def:Ind} and \ref{def:Ihatnd}.
    \end{itemize} 
\end{notation*}

\section{Background}
Throughout this article, we fix $\kk$ to be any field. Unless otherwise specified, our algebras are finite dimensional, all modules are left modules, and all $\kk$-linear triangulated categories are $\Hom$-finite. Given $\Lambda$ a $\kk$-algebra, we denote by $\Mod\Lambda$ and $\Proj \Lambda$ respectively the categories of finitely generated and of projective $\Lambda$-modules. We follow the conventions of \cite{SeidelBk} for $A_{\infty}$-categories. For two objects $X,Y$ in an $A_{\infty}$-category $\mathcal{C}$, we denote by $\hom_{\mathcal{C}}(X,Y)$ the graded vector space of morphisms from $X$ to $Y$ at the chain level, equipped with $A_{\infty}$-operations $\{\mu_n\}_{n\geq 1}$, and by $\Hom_{\mathcal{C}}(X,Y)$ the cohomology groups of ${(\hom_{\mathcal{C}}(X,Y),\mu_1)}$ in the graded cohomological category $H^*(\mathcal{C})$.

\subsection{Higher Auslander\textendash Reiten theory}\label{sec:Auslbackground}
The classical Auslander correspondence (\cite[Section~III.4]{Auslander71}, \cite[Section~VI.5]{ARS95}) establishes a bijection between Morita equivalence classes of certain $\kk$-algebras:
\begin{equation}\label{eq:Auslcorr}
    \left\{ \begin{array}{c} 
        \text{Finite-dimensional}\\ 
        \text{$\kk$-algebras of finite}\\ 
        \text{representation type} 
    \end{array}\right\}
    \quad\longleftrightarrow\quad
    \left\{ \begin{array}{c}
        \text{Finite-dimensional} \\
         \text{$\kk$-algebras  $\Gamma$  such that}\\  
         \text{$\gldim \Gamma\leq 2 \leq \domdim \Gamma$}
    \end{array}\right\}
\end{equation}
where $\gldim \Gamma$ and $\domdim\Gamma$ respectively denote the \emph{global} and \emph{dominant} dimensions of a $\kk$-algebra $\Gamma$ in the sense of \cite{Auslander55} and \cite{Tachikawa64}. We recall that the latter is defined as the largest number $d$ such that, in a minimal injective coresolution of any $\Gamma$-module $M$
\begin{equation*}
    0 \to M \to I_0 \to I_1 \to \dots \to I_{d-1} \to I_d \to \dots ,
\end{equation*}
the injective modules $I_0, I_1,\dots, I_{d-1}$ are also projective. Algebras satisfying the right-most inequalities of \eqref{eq:Auslcorr} are called \emph{Auslander algebras}.  Given a $\kk$-algebra $\Lambda$ of finite representation type, the Auslander correspondence \eqref{eq:Auslcorr} is explicitly realised by
\begin{equation}\label{eq:Auslcorrmap}
    \Lambda \mapsto \Gamma_{\Lambda}:=\End_{\Lambda}(M),
\end{equation}
for any additive generator $M$ of $\Mod\Lambda$. The \emph{Auslander\textendash Reiten quiver} associated to $\Lambda$, detailed in  \cite[Section~VII]{ARS95}, gives a combinatorial way to describe how the algebraic structure of Auslander algebras encodes part of the information of $\Mod\Lambda$. It is known \cite[Section~VII.1]{ARS95} that the Auslander\textendash Reiten quiver of $\Lambda$ coincides with the quiver underlying $\Gamma_{\Lambda}^{\op}$.

Building on the foundations of this classical theory, Iyama introduced in \cite{Iyama1,Iyama3} a new class of finite-dimensional algebras $\Gamma$, called \emph{$d$-dimensional Auslander algebras} or \emph{higher Auslander algebras}, satisfying the inequalities
\begin{equation*}
    \gldim \Gamma\leq d \leq \domdim \Gamma.
\end{equation*}
Higher Auslander algebras allowed for the formulation of the $d$-dimensional analogue of \eqref{eq:Auslcorr}, known as \emph{higher Auslander correspondence} \cite[Theorem~0.2]{Iyama1}. This powerful result gives a bijection between the set of equivalence classes of finite \emph{$d$-cluster-tilting subcategories} $\mathcal{U}^{(d)}$ of $\Mod \Lambda$ for finite-dimensional algebras $\Lambda$ admitting an additive generator, and Morita-equivalence classes of $(d+1)$-dimensional Auslander algebras. The correspondence is explicitly realised by
\begin{equation*}
    \mathcal{U}^{(d)} \mapsto \Gamma_{\Lambda}^{(d+1)}:=\End_{\Lambda}(M)
\end{equation*}
for any additive generator $M$ of $\mathcal{U}^{(d)}$, which recovers \eqref{eq:Auslcorrmap} for $d=1$. We recall (\cite[Definition~1.1]{Iyama2}) that a functorially finite subcategory $\mathcal{U}^{(d)}$ of $\Mod \Lambda$ is called \emph{$d$-cluster-tilting} if
\begin{align*}
    \mathcal{U}^{(d)}=&\{X \in \Mod \Lambda \mid  \Ext_{\Lambda}^i(X,\mathcal{U}^{(d)})=0 \quad \text{for any } 0 < i <d\} \nonumber \\
   = &\{X \in \Mod \Lambda \mid  \Ext_{\Lambda}^i(\mathcal{U}^{(d)},X)=0 \quad \text{for any } 0 < i <d\}.
\end{align*}

The construction of higher Auslander algebras is inductive in the following sense. Given a $d$-dimensional Auslander algebra $\Lambda$ (which is \emph{$d$-complete}, see \cite[Definition~1.11]{Iyama2}), Iyama proved the existence of a $d$-cluster tilting object $M$ of $\Lambda$, so that $\End_{\Lambda}(M)$ is a $(d+1)$-dimensional Auslander algebra \cite[Theorem~1.14]{Iyama2}. The latter can be explicitly described in terms of a quiver with relations, see \cite[Theorem~6.7]{Iyama2}.

The focus of this article is on \emph{higher Auslander algebras of type A}, arising from ${\Lambda=A_{n,1}}$, the path algebra over $\kk$ of the linearly oriented $A_{n}$ quiver. We denote by $\A_{n,d}$ the \linebreak $d$-dimensional Auslander algebra associated to the linearly oriented $A_{n}$ quiver, which we define below. We first define the following ordered set.

\begin{definition}\label{def:Nnd}
    For $n$ and $d$ natural numbers, with $n\geq d$, we consider the index set:
    \begin{equation*}
        \N_{n,d}:=\left\{I=(i_1,\dots,i_d)\mid 1\leq i_1<\dots <i_d\leq n \right\}.
    \end{equation*}
    We equip $\N_{n,d}$ with a lexicographic order defined by $I < J$ if there exists an ${a=1,\dots,d}$ such that ${i_h=j_h}$ for all $1 \leq h < a$ and $i_{a}<j_a$. 
\end{definition}

We recall Iyama's definition of higher Auslander algebras of type A (\cite[Section~6]{Iyama2}), though we follow the equivalent definition of Oppermann--Thomas \cite[Section~3]{OT12}.
\begin{definition}[{\cite[Theorem/Construction 3.4]{OT12}}]\label{def:And}
    The \emph{$d$-dimensional Auslander algebra of type $A_n$}, denoted by $\A_{n,d}$, is the path algebra over $\kk$ of the quiver whose set of vertices is indexed by $\N_{n,d}$, with arrows 
    \begin{equation*}
        \alpha_{I;h}: (i_1,\dots,i_{h-1},i_h,i_{h+1},\dots,i_d) \to (i_1,\dots,i_{h-1},i_h+1,i_{h+1},\dots,i_d)
    \end{equation*}
    for all $h=1,\dots,d$, and which is equipped with the following relations for all $h,k=1,\dots,d$ (note that we compose from right to left):
    \begin{itemize}
        \item If $ (i_1,\dots,i_h+1,\dots,i_d) \in \N_{n,d}$ and $ (i_1,\dots,i_k+1,\dots,i_d) \notin \N_{n,d}$,
        \begin{equation*}
            \alpha_{i_1,\dots,i_{h}+1,\dots,i_d; k} \cdot \alpha_{I; h}=0
        \end{equation*}
    \item If $ (i_1,\dots,i_h+1,\dots,i_d) \notin \N_{n,d}$ and $ (i_1,\dots,i_k+1,\dots,i_d) \in \N_{n,d}$,
    \begin{equation*}
        \alpha_{i_1,\dots,i_{k}+1,\dots,i_d; h} \cdot \alpha_{I; k}=0
    \end{equation*}
    \item If $ (i_1,\dots,i_h+1,\dots,i_d), (i_1,\dots,i_k+1,\dots,i_d) \in \N_{n,d}$,
    \begin{equation*}
        \alpha_{i_1,\dots,i_{k}+1,\dots,i_d; h} \cdot \alpha_{I; k}= \alpha_{i_1,\dots,i_{h}+1,\dots,i_d; k} \cdot \alpha_{I; h}.
    \end{equation*}
    \end{itemize}
\end{definition}
Higher Auslander algebras of type A have been the object of study in a great number of instances due to their rich combinatorial nature, see for example \cite{HI11, OT12, DJW19, DJL}. It is also known that they are $d$-complete \cite[Theorem~5.6]{IO11}. The quiver underlying $\A_{n,d}$ looks like a $(d+1)$-simplex of side length $n-d$, with commutativity relations for every small square and zero relations for each small half square \cite[Sections~1.1 and 6.2]{Iyama2}. Each $\A_{n,d}$ can be inductively constructed using the higher Auslander correspondence, and is realised by
\begin{equation*}
    \A_{n,d} \mapsto \End_{\A_{n,d}}(M)=:\A_{n+1,d+1}
\end{equation*}
for a \emph{basic} (as in \cite[Chapter~I.6]{ASS06}) $d$-cluster-tilting $\A_{n,d}$-module $M$ \cite[Theorem~5.7]{IO11} (see also \cite{Iyama2}). The quiver underlying $\A_{n+1,d+1}$ is also known as the \emph{Auslander\textendash Reiten quiver} of $\A_{n,d}$ \cite[Definition~6.1.1]{Iyama1}.

We study derived invariants of higher Auslander algebras of type A. The perfect derived category $\perf(\A_{n,d})$ has a distinguished collection of generators given by indecomposable projective $\A_{n,d}$-modules, which we label as follows. For $I\in \N_{n,d}$, denote by $e_I$ the idempotent element of $\A_{n,d}$ corresponding to the path of length $0$ starting at the vertex of the quiver underlying $\A_{n,d}$ at the lattice point $I$. Further denote by $I^{\op}$ the tuple:
\begin{equation*}
    I^{\op}:=(n+1-i_d,\dots,n+1-i_1)
\end{equation*}
and by $\P_I:=\A_{n,d} e_{I^{\op}}$ the (left) indecomposable projective module associated to the vertex $I^{\op}$. The endomorphism algebra of this collection of projectives is isomorphic to $\A_{n,d}$ as graded $\kk$-algebras \cite[Section~II.2]{ARS95}. Using the combinatorial presentation of higher Auslander algebras of type A, we note there is a non-zero morphism from $\P_I$ to $\P_J$, for $I,J \in \N_{n,d}$, exactly when 
\begin{equation}\label{eq:intertwining}
    i_1 \leq j_1 < i_2 \leq j_2 < \dots < i_d \leq j_d
\end{equation}
(\cite[Theorem~3.6]{OT12}); $J$ is said to be \emph{intertwining} $I$ whenever \eqref{eq:intertwining} holds. We define a total order of the collection of indecomposable projective $\A_{n,d}$-modules by $\P_I<\P_J$ whenever $I<J$ in $\N_{n,d}$, as in Definition~\ref{def:Nnd}.

Following \cite{Iyama1}, we identify the category of projective $\A_{n+1,d+1}$-modules \linebreak $\Proj(\A_{n+1,d+1})$ with the $d$-cluster tilting subcategory $\mathcal{U}^{(d)} \subset \Mod \A_{n,d}$. We also identify $\Proj(\A_{n,d})$ with the subcategory of $\Mod \A_{n,d}$ of projective $d$-cluster tilting modules. Under these identifications and for $I \in \N_{n,d}$, $\P_{1,I+1}$ denotes the indecomposable $d$-cluster tilting $\A_{n,d}$-module identified with the indecomposable projective $\A_{n,d}$-module $\P_I$, where $I+1:=(i_1+1,\dots,i_d+1)$. For $J\in \N_{n+1,d+1}$ such that $j_1 \neq 1$ and for $h=1,\dots,d+1$, denote by $\CTmod_J$ the indecomposable $d$-cluster tilting $\A_{n,d}$-module identified with the projective $\A_{n+1,d+1}$-module $\P_J$. Denote by $\P_{\widehat{j_h}}$ the indecomposable projective $\A_{n,d}$-module, where
\begin{equation*}
    \widehat{j_h}:=J_{\widehat{j_h}}=(j_1-1,\dots,j_{h-1}-1,j_{h+1}-1,\dots,j_{d+1}-1).
\end{equation*}
The following proposition is due to \cite[Proposition~3.17]{OT12}, see also \cite[Proposition~2.7]{JK19}.

\begin{proposition}[{\cite[Proposition~3.17]{OT12}}]\label{prop:projresCT}
    For $J \in \N_{n+1,d+1}$ such that $j_1\neq 1$, let $\CTmod_J$ be the indecomposable direct summand of the $d$-cluster tilting $\A_{n,d}$-module $M$ as above. It admits a minimal resolution by indecomposable projective $\A_{n,d}$-modules
    \begin{equation*}
        \P_{\widehat{j_{d+1}}}\to  \P_{\widehat{j_{d}}}\to \dots \to  \P_{\widehat{j_{1}}} \to  \CTmod_J \to 0
    \end{equation*}
    where the morphisms are the unique ones determined by the intertwining conditions \eqref{eq:intertwining}.\qed
\end{proposition}

\subsection{Fukaya categories of symmetric products of surfaces}\label{sec:Fuksym}
The setting we describe in this section is the same as that in \cite{DJL}, due to Auroux \cite{Auroux1} and based on the construction in bordered Heegaard Floer homology by Lipshitz\textendash Ozsv\'{a}th\textendash Thurston \cite{LOT18}. We review the relevant background, and refer to these for a more detailed account.

Let $\Sigma$ be a compact Riemann surface with non-empty boundary, and $d$ a positive integer. We will denote by $\Sym^d(\Sigma)$ the $d$-fold symmetric product of $\Sigma$. $\Sym^d(\Sigma)$ can be equipped with a natural symplectic structure $\omega$ coming from a choice of positive area form $\alpha$ on $\Sigma$, as prescribed by Perutz in \cite[Corollary~7.2]{Perutz}. Away from the diagonal, $\omega$ is the smooth pushforward $\bpi_*(\alpha^{\times d})$, where $\bpi: \Sigma^{\times d} \to \Sym^d(\Sigma)$ is the branched covering map which ramifies along the big diagonal of $\Sigma^{\times d}$. This choice of symplectic structure allows for a combinatorially nice description of certain exact Lagrangians in $\Sym^d(\Sigma)$, given by products of pairwise disjoint exact Lagrangian submanifolds of $\Sigma$. In particular, $\Sym^d(\Sigma)$ is a symplectic manifold with corners, and a Liouville domain (as in \cite{EG91, Eliashberg}). For $\Lambda \subset \partial \Sigma$ a non-empty finite collection of points on the boundary of $\Sigma$, the divisor
\begin{equation}\label{eq:DJLstopsSigma}
    \Lambda^{(d)}:= \bigcup_{p \in \Lambda} \{p\} \mathrel{\boldsymbol{\times}} \Sym^{d-1} (\Sigma)
\end{equation}
defines the union of symplectic hypersurfaces in $\Sym^d(\Sigma)$ under the identification 
\begin{equation*}
    \left(p, \{x_1,\dots,x_{d-1}\}\right) \mapsto  \{p, x_1,\dots,x_{d-1}\}
\end{equation*}
of each point on the right-hand side of \eqref{eq:DJLstopsSigma} with a point in $\Sym^d(\Sigma)$. We have the following generation result for Auroux' partially wrapped Fukaya category $\W^{\textnormal{Aur}}(\Sym^d(\Sigma), \Lambda^{(d)})$, due to Auroux. 

\begin{theorem}[{\cite[Theorem~1]{Auroux1}}]\label{thm:Auroux}
    Let $\Lambda$ be a finite set of points on $\partial \Sigma$ as above, $L_1,...,L_m$ a collection of disjoint properly embedded arcs in $\Sigma$ with endpoints in $\partial \Sigma\setminus \Lambda$. Assume that $\Sigma \setminus (L_1 \cup ...\cup L_m)$ is a disjoint union of discs, each of which contains at most one point of $\Lambda$. Then, the partially wrapped Fukaya category $\W^{\textnormal{Aur}}(\Sym^d(\Sigma), \Lambda^{(d)})$ is generated by the $m \choose d$ Lagrangian submanifolds $L_{i_1,\dots,i_d} := L_{i_1} \times \dots \times L_{i_d}$, products of distinct pairs of arcs.\qed
\end{theorem}

Morphisms between products of Lagrangian arcs satisfying the hypotheses of Theorem~\ref{thm:Auroux} arise from products of Reeb chords induced by the Reeb flows along $\partial \Sigma$, which are the rotational flows in the counter-clockwise direction of the boundary \cite[Sections~4.1, 4.2]{Auroux2} (see also \cite[Definition~8]{Auroux1}). It will sometimes be useful to identify the endomorphism algebra $\mathcal{A}$ of products of a given collection of arcs with a certain \emph{strands algebra} arising from work of \cite{LOT18} on bordered Heegaard Floer homology.

We briefly recall, and refer to \cite[Section~3.1.2]{LOT18} for a more detailed construction, that a strands algebra is defined in terms of \emph{strands diagrams}. These are combinatorially given by a pair of subsets $S$ and $T$ of $\{1,\dots,\ell\}$, each of order $k$, and a bijection $\phi: S \to T$.
\begin{dfnprop}[{\cite[Section~3.1]{LOT18}}]
    Let $\kk$ be a field of characteristic 2. The \emph{strands algebra with $\ell$ strands in $k$ places} is the vector space over $\kk$ freely generated by strands diagrams. It has the structure of a differential graded algebra. Composition is given by horizontal concatenation (vanishing if two diagrams cannot be concatenated), with the condition that the product is zero if two strands cross more than once. The \linebreak $\Z$-grading is given by the number of crossings, and the differential corresponds to resolving crossings.
\end{dfnprop}

See Figures~\ref{fig:strands} and \ref{fig:strandsdifferential} for a depiction of composition and differential of strands diagrams (see also \cite[Equation~3.6]{LOT18} for the latter).

\begin{figure}[t]
    \centering
    \begin{minipage}{1 \textwidth}
        \centering
        \begin{tikzpicture}[scale=.7, v/.style={draw,shape=circle, fill=black, minimum size=1.2mm, inner sep=0pt, outer sep=0pt}]
            \node[v,label=180:$4$] at (0-1,4){};
            \node[v,label=180:$3$] at (0-1,3){};
            \node[v,label=180:$2$] at (0-1,2){};
            \node[v,label=180:$1$] at (0-1,1){};
            
            \node[v,label=0:$4$] at (1-1,4){};
            \node[v,label=0:$3$] at (1-1,3){};
            \node[v,label=0:$2$] at (1-1,2){};
            \node[v,label=0:$1$] at (1-1,1){};
    
            \draw  (0-1,1) to (1-1,3);
            \draw  (0-1,2) to (1-1,2);
    
            \node at (2-.5,2.5) {$\boldsymbol{\cdot}$};
    
            \node[v,label=180:$4$] at (3,4){};
            \node[v,label=180:$3$] at (3,3){};
            \node[v,label=180:$2$] at (3,2){};
            \node[v,label=180:$1$] at (3,1){};
            
            \node[v,label=0:$4$] at (4,4){};
            \node[v,label=0:$3$] at (4,3){};
            \node[v,label=0:$2$] at (4,2){};
            \node[v,label=0:$1$] at (4,1){};
    
            \draw  (3,2) to (4,4);
            \draw  (3,3) to (4,3);
    
            \node at (5+.5,2.5) {$=$};
    
            \node[v,label=180:$4$] at (6+1,4){};
            \node[v,label=180:$3$] at (6+1,3){};
            \node[v,label=180:$2$] at (6+1,2){};
            \node[v,label=180:$1$] at (6+1,1){};
            
            \node[v,label=0:$4$] at (7+1,4){};
            \node[v,label=0:$3$] at (7+1,3){};
            \node[v,label=0:$2$] at (7+1,2){};
            \node[v,label=0:$1$] at (7+1,1){};
    
            \draw  (6+1,1) to[out=60,in=180] (7+1,3);
            \draw  (6+1,2) to[out=0,in=240] (7+1,4);
    
            \node at (8+1.5,2.5) {$=$};
            \node at (8+2.5,2.5) {$0$};
        \end{tikzpicture}
\end{minipage}
    \caption{Left-to-right composition of strands diagrams, in the strands algebra with $2$ strands in $4$ places.}
    \label{fig:strands}
\end{figure}

\begin{figure}[t]
    \centering
    \begin{minipage}{1 \textwidth}
        \centering
        \begin{tikzpicture}[scale=.55, v/.style={draw,shape=circle, fill=black, minimum size=1.2mm, inner sep=0pt, outer sep=0pt}]
            \node at (-2,3) {$\partial$};
            \draw[line width=1.3] (-.75,5.3) to[out=240,in=120,looseness=.8] (-.75,.7);
            \draw[line width=1.3] (2.75,5.3) to[out=-60,in=60,looseness=.8] (2.75,.7);

            \node[v,label=180:$5$] at (0,5){};
            \node[v,label=180:$4$] at (0,4){};
            \node[v,label=180:$3$] at (0,3){};
            \node[v,label=180:$2$] at (0,2){};
            \node[v,label=180:$1$] at (0,1){};
            
            \node[v,label=0:$5$] at (2,5){};
            \node[v,label=0:$4$] at (2,4){};
            \node[v,label=0:$3$] at (2,3){};
            \node[v,label=0:$2$] at (2,2){};
            \node[v,label=0:$1$] at (2,1){};
    
            \draw  (0,1) to[out=60,in=180] (2,5);
            \draw  (0,2) to[out=0,in=230] (2,4);
            \draw  (0,3) to (2,3);
    
            \node at (4,3) {$=$};

            \node[v,label=180:$5$] at (5.5,5){};
            \node[v,label=180:$4$] at (5.5,4){};
            \node[v,label=180:$3$] at (5.5,3){};
            \node[v,label=180:$2$] at (5.5,2){};
            \node[v,label=180:$1$] at (5.5,1){};
            
            \node[v,label=0:$5$] at (7.5,5){};
            \node[v,label=0:$4$] at (7.5,4){};
            \node[v,label=0:$3$] at (7.5,3){};
            \node[v,label=0:$2$] at (7.5,2){};
            \node[v,label=0:$1$] at (7.5,1){};
    
            \draw  (5.5,1) to[out=80,in=180,looseness=1.5] (7.5,3);
            \draw  (5.5,2) to[out=0,in=230] (7.5,4);
            \draw  (5.5,3) to[out=0,in=180] (7.5,5);

            \node at (9,3) {$+$};

            \node[v,label=180:$5$] at (10.5,5){};
            \node[v,label=180:$4$] at (10.5,4){};
            \node[v,label=180:$3$] at (10.5,3){};
            \node[v,label=180:$2$] at (10.5,2){};
            \node[v,label=180:$1$] at (10.5,1){};
            
            \node[v,label=0:$5$] at (12.5,5){};
            \node[v,label=0:$4$] at (12.5,4){};
            \node[v,label=0:$3$] at (12.5,3){};
            \node[v,label=0:$2$] at (12.5,2){};
            \node[v,label=0:$1$] at (12.5,1){};
    
            \draw  (10.5,1) to[out=60,in=180] (12.5,4);
            \draw  (10.5,2) to[out=60,in=180] (12.5,5);
            \draw  (10.5,3) to (12.5,3);

            \draw (4.5,1) to (4.5,.5) to (8.5,.5) to (8.5,1);
            \node at (6.5,0) {$=0$};

            \node at (14,3) {$+$};

            \node[v,label=180:$5$] at (15.5,5){};
            \node[v,label=180:$4$] at (15.5,4){};
            \node[v,label=180:$3$] at (15.5,3){};
            \node[v,label=180:$2$] at (15.5,2){};
            \node[v,label=180:$1$] at (15.5,1){};
            
            \node[v,label=0:$5$] at (17.5,5){};
            \node[v,label=0:$4$] at (17.5,4){};
            \node[v,label=0:$3$] at (17.5,3){};
            \node[v,label=0:$2$] at (17.5,2){};
            \node[v,label=0:$1$] at (17.5,1){};

            \draw  (15.5,1) to[out=60,in=180] (17.5,5);
            \draw  (15.5,2) to[out=0,in=180] (17.5,3);
            \draw  (15.5,3) to[out=0,in=180] (17.5,4);

        \end{tikzpicture}
\end{minipage}
    \caption{Differential in the strands algebra  with $3$ strands in $5$ places.}
    \label{fig:strandsdifferential}
\end{figure}

When each component of $\partial \Sigma$ has exactly one marked point, the identification between $\mathcal{A}$ and the appropriate strands algebra is done by Auroux in \cite{Auroux2} (Theorem~1.2 therein) over a ring of characteristic 2. In his description, generators of the morphism space between two objects in $\W^{\textnormal{Aur}}(\Sym^d(\Sigma), \Lambda^{(d)})$ of the form given in Theorem~\ref{thm:Auroux} are in bijection with diagrams with $m$ strands in $d$ places. If the arcs in the Theorem~\ref{thm:Auroux} satisfy the additional assumption that each disc in $\Sigma \setminus (\bigcup L_i)$ contains \emph{exactly} one stop, Auroux further proves that the higher products in the partially wrapped Fukaya category involving these generators identically vanish \cite[Proposition~3.6]{Auroux2}. This is done by choosing particularly nice perturbation data for this set of generators, and by computing the Maslov index of holomorphic $n$-gons ($n\geq 3$) with boundary on them and verifying that this is non-negative. As the $\mu_n$-higher product involves a count of \emph{rigid} holomorphic discs, i.e.\ discs of Maslov index $2-n$, Auroux concludes that there are no rigid holomorphic $n$-gons with boundary on these generators.

We can replicate the above construction also if we allow each component of $\partial \Sigma$ to contain more than one element of $\Lambda$, which Auroux does in \cite[Proposition~11]{Auroux1}, still over a ring of characteristic 2. In particular, he provides the following.
\begin{proposition}[{\cite[Proposition~11]{Auroux1}}]
    Let $\kk$ be a field of characteristic 2, and fix a collection of arcs in $\Sigma$ satisfying the hypotheses of Theorem~\ref{thm:Auroux}. Let $S=\{s_i\}_i$ and $T=\{t_j\}_j$ be two subsets of $\{1,\dots,m\}$ of order $d$, and $\Phi$ the set of bijective maps from $S$ to $T$. If $L=\prod_i L_{s_i}$ and $L'=\prod_j L_{t_j}$, then the $\kk$-vector space $\hom(L,L')$ in $\W^{\textnormal{Aur}}(\Sym^d(\Sigma), \Lambda^{(d)})$ admits a basis indexed by the elements of
    \begin{equation}\label{eq:strandshoms}
        \bigsqcup_{\phi \in \Phi} \left(\hom(L_{s_1},L'_{\phi(s_1)}) \otimes \dots \otimes \hom(L_{s_d},L'_{\phi(s_d)})\right),
    \end{equation}
    where each $\hom(L_{i},L'_{j})$ is computed in $\W^{\textnormal{Aur}}(\Sigma,\Lambda)$. Moreover, the differential corresponds to resolving crossings.
\end{proposition}

We provide an explicit computation of the strands algebra associated to the endomorphism algebra of a full collection of a partially wrapped Fukaya category in Example~\ref{ex:strands} below. Auroux' construction can be lifted to arbitrary characteristic for particularly nice choices of Lagrangian arcs, as is carefully done in \cite[Section~3.1]{DJL}. As in the particular case of one marked point per boundary component of $\Sigma$, if the arcs satisfy the additional assumption that each disc in $\Sigma \setminus (\bigcup L_i)$ contains exactly one stop, the higher products identically vanish. The existence of rigid holomorphic $n$-gons (for $n \geq 3$) with boundary on the above generators is independent of the characteristic of the ground ring, so the higher products vanish in arbitrary characteristic \cite[Section~2.1.6]{DJL}.

For $n \in \mathbb{N}_{\geq 1}$ and following \cite{DJL}, we consider the disc $\Sigma = \D$ (with standard orientation of the boundary) and a set of $n+1$ points $\Lambda_n$ on $\partial \D$. Let $\Lambda^{(d)}_n$ define the collection of symplectic hypersurfaces in $\Sym^d(\D)$ (which we call \emph{stops}) as in \eqref{eq:DJLstops}; following \cite{Auroux2}, we denote by
\begin{equation*}
    \W_n^d=\W^{\textnormal{Aur}}(\Sym^d(\D), \Lambda_n^{(d)})
\end{equation*}
the partially wrapped Fukaya category of $\Sym^d(\D)$ with stops in $\Lambda^{(d)}_n$. The obstructions $2c_1(\Sym^d(\D))$ and $H^1(\Sym^d(\D))$ to the existence and uniqueness of a canonical grading on $\Sym^d(\D)$ both vanish, and we equip $\W_n^d$ with this $\Z$-grading \cite[Section~4]{Seidel99}. Products of arcs in $\D \setminus \Lambda_n$ are contractible, hence admit a unique choice of $\Z$-grading, up to a global shift; a choice of grading on such Lagrangian subspaces determines an object of the partially wrapped Fukaya category.  Note that if $n<d$, it can be shown that $\W_n^d$ is trivial (see \cite[Section~2.2]{DJL}).

\begin{notation*}
    Following the notation in \cite{DJL}, we introduce a labelling on the points $p_i \in \Lambda_n$ induced by the opposite orientation of the boundary $\partial \D$, so that $p_0 < \dots < p_n$ are ordered counter-clockwise. Denote by $i$ the boundary component of $\D \setminus \Lambda_n$ which lies between $p_i$ and $p_{i+1}$. Let $L_{ij}$ be the properly embedded arc in $\D \setminus \Lambda_n$ with endpoints lying on the components $i$ and $j$ of $\D \setminus \Lambda_n$. We identify these arcs (equipped with unique gradings) with objects of the partially wrapped Fukaya category $\W_n^1$.
\end{notation*}

\begin{example}\label{ex:strands}
    Let $\Sigma=\D$ be the standard disc, and $\Lambda=\Lambda_4$ be the set of $5$ points on its boundary. We consider the following set of arcs: $\{L_{1},L_{2},L_{3},L_{4}\}$, where ${L_i:=L_{0i}}$ following the above notation. We remark that this is a particularly nice collection, which Dyckerhoff\textendash Jasso\textendash Lekili consider in \cite{DJL}, see Figure~1 therein. Let us compute $\hom(L,L')$ in $W_4^3$, for $L=L_1\times L_2\times L_3$, $L'=L_2 \times L_3\times L_4$ (for simplicity, over a ring of characteristic $2$). The collection of arcs satisfies the hypothesis that $\D\setminus (\cup L_i)$ contains at least one element of $\Lambda$, so its endomorphism algebra in $W_4^3$ is a differential graded algebra. By \cite[Proposition~11]{Auroux1}, $\hom(L,L')$ admits a basis indexed by strands diagrams $\phi$ of $4$ strands in $3$ places, from $S=\{1,2,3\}$ to $T=\{2,3,4\}$, as subsets of $\{1,2,3,4\}$. There are four such diagrams:
    \begin{itemize}
        \item $\phi_1$, given by $\phi_1(1)=4$, $\phi_1(2)=2$, and $\phi_1(3)=3$;
        \item $\phi_2$, given by $\phi_2(1)=3$, $\phi_2(2)=2$, and $\phi_2(3)=4$;
        \item $\phi_3$, given by $\phi_3(1)=2$, $\phi_3(2)=4$, and $\phi_3(3)=3$;
        \item $\phi_4$, given by $\phi_4(1)=2$, $\phi_4(2)=3$, and $\phi_4(3)=4$.
    \end{itemize}
    In particular, each of these corresponds to the generator
    \begin{equation*}
        \hom_{\W_4^1}(L_1,L_{\phi_h(1)})\otimes\hom_{\W_4^1}(L_2,L_{\phi_h(2)})\otimes\hom_{\W_4^1}(L_3,L_{\phi_h(3)})
    \end{equation*}
    for $h=1,2,3,4$, and this is a complete set of generators for $\hom(L,L')$ in $\W_4^3$. The differential in the strands algebra is given by $\partial (\phi_1)=\phi_2 + \phi_3$, $\partial (\phi_2)=\partial (\phi_3)=\phi_4$, and this gives the corresponding differential for $\hom(L,L')$.
\end{example}

With respect to suitable grading structures, Auroux provides the following exact triangles in $\W^d_n$.

\begin{lemma}[{\cite[Lemma~5.2]{Auroux2}}]\label{lem:Aurouxtriangles}
    Let $0\leq i < j < k \leq n$ and let $d\geq 2$. The following statements hold.
    \begin{enumerate}[i.]
        \item Let $L$ be the product of $d-2$ distinct arcs satisfying the hypotheses of Theorem~\ref{thm:Auroux}, which are not homotopic to any of the arcs $L_{ij}$, $L_{jk}$ and $L_{ik}$. There is a quasi-isomorphism
        \begin{equation*}
            L \times L_{ij} \times L_{ik} \xrightarrow{\text{ $\simeq$ }} L \times L_{ij} \times L_{jk}.
        \end{equation*}
        \item Let $L$ be the product of $d-1$ distinct arcs satisfying the hypotheses of Theorem~\ref{thm:Auroux}, which are not homotopic to any of the arcs $L_{ij}$, $L_{jk}$ and $L_{ik}$. There is a non-split exact triangle
        \begin{equation*}
            L \times L_{ij} \rightarrow L \times L_{ik} \rightarrow L \times L_{jk} \rightarrow (L \times L_{ij})[1]
        \end{equation*}
        where $[1]$ denotes a shift in grading.\qed
    \end{enumerate}
\end{lemma}

The main result in \cite{DJL} provides the quasi-equivalence of triangulated $A_{\infty}$-categories between the perfect derived $A_{\infty}$-category of higher Auslander algebras of type A and Auroux' partially wrapped Fukaya category of $\Sym^d(\D)$ with stops in \eqref{eq:DJLstops}. Underlying this equivalence, they provide a quasi-isomorphism of differential $\Z$-graded algebras between $\A_{n,d}$ and the endomorphism algebra of the $d$-fold product of a preferred collection of arcs satisfying the assumptions of Theorem~\ref{thm:Auroux}, as given in \cite[Figure~1]{DJL}. We will refer to this result in Section~\ref{section:F-W}.

\subsection{Fukaya\textendash Seidel categories}\label{sec:FSbackground}
We consider an exact symplectic, convex at infinity manifold with corners $(X,\omega)$. Following \cite[Section~15]{SeidelBk}, we study exact symplectic Lefschetz fibrations $f:X\to \D$, whose base is the standard disc $\D\subset \C$. Fixing a regular value $*\in \partial \D$, the corresponding regular fibre $F_*:=f^{-1}(*)$ is a smooth manifold with corners carrying an exact symplectic structure inherited from $\omega$. As per \cite[Section~16b]{SeidelBk}, we consider a distinguished collection of embedded paths $\gamma_i: [0,1] \to (\D\setminus \partial \D)\cup \{*\}$ (called \emph{vanishing paths}) indexed by the critical values of $f$. We require $\gamma_i(0)=*$ for all $i$, $\{\gamma_i(1)\}_i$ to be the set of critical values of $f$, and the paths to be pairwise disjoint away from $\gamma_i(0)$ (Figure~\ref{fig:diagramvpaths}). For $t$ close to critical values of $f$, we construct \emph{vanishing cycles} as Lagrangian submanifolds of $F_t=f^{-1}(t)$ which collapse to a point as we approach each critical fibre $F_{\gamma_i(1)}=f^{-1}(\gamma_i(1))$. Using symplectic parallel transport in $(X,\omega)$, we can transport all vanishing cycles along their respective vanishing paths so that they all lie in the fibre $F_*$. To each vanishing path we can associate a \emph{Lefschetz thimble} $D_i:=D_{\gamma_i}$ as the union of all corresponding vanishing cycles above that path: these are embedded Lagrangian discs in $X$, whose boundaries $\partial D_i = D_i \cap F_*$ are the vanishing cycles $V_i \subset F_*$. Vanishing cycles and thimbles are naturally ordered by the (clockwise) ordering of the vanishing paths, as given by the clockwise ordering of the angles at the common intersection point $\gamma_i(0)$. We fix an indexing of the cycles such that $V_i < V_j$ whenever $i < j$.

\begin{figure}[t]
    \centering
    \begin{tikzpicture}
         \draw (0,0) ellipse (1.5 and 1);
 
         \draw[] (.1,.6) -- (-.1,.4);
         \draw[] (-.1,.6) -- (.1,.4);
  
         \draw[] (-.4,.1) -- (-.6,-.1);
         \draw[] (-.6,.1) -- (-.4,-.1);
 
         \draw[] (.1,-.6) -- (-.1,-.4);
         \draw[] (-.1,-.6) -- (.1,-.4);
 
         \node[] (star) at (1.5,0) {$*$};
 
         \draw[] (1.5,0) -- (0,.5);
         \draw[] (1.5,0) -- (-.5,0);
         \draw[] (1.5,0) -- (0,-.5);
     \end{tikzpicture}
     \caption{Base of a Lefschetz fibration, with regular value on $\partial \D$.}
     \label{fig:diagramvpaths}
 \end{figure}

As detailed in \cite[Chapter~3]{SeidelBk}, one can define the Fukaya\textendash Seidel category $\F(f)$ associated to a symplectic Lefschetz fibration $f$. Objects of this category are twisted complexes of \emph{Lagrangian branes} $D^{\#}_{i}$, each consisting of a Lefschetz thimble $D_i$ equipped with a brane structure (a choice of a spin structure and grading in the appropriate group). Morphisms between Lagrangian branes are given by Floer complexes. The ordered collection of Lagrangian thimbles associated to a choice of a distinguished collection of vanishing paths is a full exceptional collection in the Fukaya--Seidel category $\F(f)$, and the latter is independent (up to quasi-equivalence) of the choices of paths after taking twisted complexes \cite[Section~18j]{SeidelBk}. We can also equip each vanishing cycle $V_{i}$ with a brane structure $V^{\#}_i$ induced by that of its Lefschetz thimble; this turns each vanishing cycle into an object of the (compact) Fukaya category $\Fuk(F_*)$ of the regular fibre of $f$. This allows us to define a restriction functor
\begin{equation}\label{eq:restrictionfunctor}
    \F(f)\to \Fuk(F_*)
\end{equation}
given by restricting Lefschetz thimbles to their associated vanishing cycles. Essentially by construction (due to Seidel \cite[Section~18e]{SeidelBk}), we have an isomorphism of Floer complexes:
\begin{equation}\label{eq:Seidelisom}
    CF^{*}_{\F(f)}(D^{\#}_i,D^{\#}_j) \cong  CF^{*}_{\Fuk(F_*)}(V^{\#}_i,V^{\#}_j)
\end{equation}
whenever $i<j$, while each morphism space $CF^{*}_{\F(f)}(D^{\#}_i,D^{\#}_j)$ vanishes whenever $i>j$ and is isomorphic to the base ring for $i=j$. Concretely, this allows us to carry out Floer cohomology computations in the directed $A_{\infty}$-subcategory of $\Fuk(F_*)$ associated to an ordered collection of vanishing cycles. Unless otherwise specified, our Fukaya categories are all derived, so that they are invariants of the Lefschetz fibration.

Fix an ordered collection of vanishing paths $\gamma_1,\dots,\gamma_{\mu}$, where $\mu$ is the number of critical values of $f$. The braid group $Br_{\mu}$ acts freely on the set of all distinguished collections, and the action of the standard $(i-1)^{\textnormal{th}}$ generator of $Br_{\mu}$ gives rise to a Hurwitz-type move:
\begin{equation*}
    (\gamma_1, \dots, \gamma_{i-2}, \gamma_{i-1}, \gamma_{i}, \gamma_{i+1},\dots \gamma_{\mu})\mapsto (\gamma_1, \dots, \gamma_{i-2}, \tau_{\gamma_{i-1}}(\gamma_i), \gamma_{i-1}, \gamma_{i+1},\dots \gamma_{\mu})
\end{equation*}
where $\tau_{\gamma_{i-1}}(\gamma_i)$ is the vanishing path obtained by precomposing $\gamma_i$ with a clockwise loop around $\gamma_{i-1}$. This lifts to a Hurwitz-type move on the vanishing cycles. Furthermore, if $V'_{i}$ is the vanishing cycle associated to $\tau_{\gamma_{i-1}}(\gamma_i)$, this is realised as $V'_{i}=\tau_{V_{i-1}}(V_i)$ where $\tau$ denotes the symplectic Dehn twist as in \cite[Section~16c]{SeidelBk}. A Hurwitz move on two consecutive and disjoint vanishing cycles leaves the geometric objects unchanged, but it inverts their order. It is a non-trivial result by Seidel \cite[Corollary~17.17]{SeidelBk} that, in the derived Fukaya category of the regular fibre,
\begin{equation*}
    T_{V_{i-1}^{\#}}(V_{i}^{\#})\cong \tau_{V_{i-1}^{\#}}(V_{i}^{\#}),
\end{equation*}
where $T$ denotes the twist functor around a spherical object of a triangulated $A_{\infty}$-category \cite[Section~5h]{SeidelBk}.

At least when the Lagrangians are orientable, the morphism spaces in $\F(f)$ are \linebreak $\Z_2$-gradable. When $2c_1(X)$ vanishes, this can be extended to a $\Z$-grading \cite[Section~12a]{SeidelBk}. If $H^1(X)$ also vanishes, $X$ carries a canonical grading \cite[Section~4]{Seidel99}, which induces a unique (up to shifts) $\Z$-grading on the Floer complexes.

The setting of Lefschetz fibrations on $(X,\omega)$ described above is the one used in \cite[Chapter~3]{SeidelBk}, is that of Liouville domains. By considering a compact symplectic manifold $X$ with boundary, we can use the Liouville flow to construct the non-compact Liouville completion $X^{\circ}$ of $X$. Using the boundary regularity conditions of an exact symplectic Lefschetz fibration $f$ \cite[Section~15a]{SeidelBk}, we can extend this to a fibration $X^{\circ} \to \C$, and we can use symplectic parallel transport to ``push'' the generic fibre to lie at a point far away from the critical locus (``at infinity''). We will use the two settings (of symplectic Lefschetz fibrations on compact manifolds with boundary and on their Liouville completions) interchangeably.

We think of Lefschetz fibrations as the complex analogue of Morse functions \cite{Kas80}. Specifically, given $f:X \to \D$ with generic fibre $F_*$, $X$ is the symplectic manifold obtained by attaching top-dimensional handles to $F_*$, in the following way. Consider a small disc $\D_{\epsilon}\subset \D$ on the base of $f$, centred at a regular value, which does not contain any critical value. Above $\D_{\epsilon}$, the fibration is the trivial $\D_{\epsilon}\times F_*$. By increasing $\epsilon$, the disc will cross the critical values of $f$; whenever it does, the topology of the total space above it will change. This is prescribed by the attachment of a top dimensional handle to the vanishing cycle associated to the critical value we crossed \cite{Kas80}. In particular, this gives a handlebody decomposition of $X$, provided we have one for $F_*$.

In this article, we consider the symplectic manifold $X^{\circ}=\C^d$, equipped with standard symplectic form $\omega$ and complex structure $J$. We consider an exact symplectic Lefschetz fibration $f:\C^d \to \C$ with critical locus contained in the preimage of a compact set in the base (which we can think to be the standard disc $\D$, up to an isotopy of the base). The map
\begin{equation*}
    H(x_1,\dots,x_d)= \lvert x_1\rvert ^2+\dots+\lvert x_d\rvert^2
\end{equation*}
defines a plurisubharmonic function (\cite[Section~19a]{SeidelBk}), whose negative gradient flow points strictly inwards along $f^{-1}(\partial \D)$. The preimage of $\D$ inherits the structure of exact symplectic manifold with corners, which is also convex at infinity with respect to the negative gradient flow of $H$, and we define the Fukaya--Seidel category of $f$ as that of $f|_{f^{-1}(\D)}$. Note that $2c_1(\C^d)$ and $H^1(\C^d)$ both vanish, and $\C^d$ carries a canonical grading (\cite[Section~12]{SeidelBk} and \cite[Section~4]{Seidel99}), which induces a $\Z$-grading on the Fukaya\textendash Seidel category of $f$.

Seidel's definition of Fukaya\textendash Seidel categories can be extended to the context of isolated singularities, in the following sense. Let $f_0: \C^d \to \C$ be a holomorphic function with an isolated singularity at the origin. An \emph{isolated hypersurface singularity} is the holomorphic complex function germ of such $f_0$. By a \emph{Morsification} $f$ of $f_0$ we mean a representative $\C^d \times \C, (z,\epsilon) \mapsto f_{\epsilon}(z)$ of a holomorphic function germ satisfying $f(z,0)=f_0(z)$ for $z\in \C^d$ and such that $f_{\epsilon}:\C^d \to \C$ is a Morse function for almost all ${\epsilon}$ in a neighbourhood of zero. We often call $f$ a Morsification of $f_0$. For an isolated singularity $f_0$, the \emph{Milnor number} is the dimension of the Jacobi algebra, defined as
\begin{equation*}
    \C[x_1,\dots,x_n]/\text{Jac}f_0:=\C[x_1,\dots,x_d]/(\partial_1 f_0,\dots,\partial_d f_0).
\end{equation*}
The \emph{Milnor fibre} of $f_0$ is the fibre above a (sufficiently small) regular value of $f_{\epsilon}$.

If $f$ is a Morsification of $f_0$, its Milnor number $\mu(f)=\mu(f_0)$ coincides with the number of the (finitely many) isolated singularities near the origin. If the perturbation is generic and sufficiently small, the Milnor fibre can be identified with (the Liouville completion of) the standard notion of the Milnor fibre of $f_0$, defined as the intersection between a suitably small ball around the origin, and the regular fibre of $f_0$ above a point close to 0 (see \cite[Lemma~2.18]{Keating15}, \cite[Lemma~3.3]{Dimca} for the smooth version). With slight abuse of notation, we use the same name interchangeably.

Given an isolated hypersurface singularity $f_0$, we define the Fukaya\textendash Seidel category $\F(f_0):=\F(f)$ to be that of any Morsification $f$ of $f_0$. We note that the space of deformations of $f_0$ which are Morse is path connected \cite[Section~2.3.1]{Keating15}, so any two Morsifications are equivalent. The $A_{\infty}$-category $\F(f_0)$ is independent of a choice of such representative. Moreover, holomorphicity of the perturbations ensures that Morsifications of $f_0$ are symplectic Lefschetz fibrations. As with Lefschetz fibrations, when computing the Fukaya\textendash Seidel category of a hypersurface singularity we will interchangeably study the singularity on a compact manifold with boundary (a small ball around the origin), and its extension on the Liouville completion, which is compatible with the outwards pointing Liouville flow.

\begin{remark}{\label{rk:furtherperturbation}}
    A Lefschetz fibration $f$ is usually required to have distinguished singularities, and this is generically obtainable. However, it will sometimes be useful to relax this condition, and allow distinct nodal singularities to lie above the same critical value (see Definition~\ref{def:nongenericness}). However, when necessary, one can further perform a small perturbation of $f$ to separate the critical values \cite[Remark~2.32]{Keating15}.
\end{remark}

In Section~\ref{section:quotientFC}, we will study isolated singularities on $\Sym^d(\C)$. This inherits a symplectic structure under the natural identification:
\begin{equation}\label{eq:iso}
    \varphi: \Sym^d(\C) \xrightarrow{\text{ $\cong$ }} \C^d, \quad \mathbf{x}\mapsto (e_1(\mathbf{x}),\dots,e_d(\mathbf{x})),
\end{equation}
where $\mathbf{x}=\{x_1,\dots,x_d\}$ is the unordered collection of $d$ points in $\C$, and $e_i(\mathbf{x})$ is the $i^{\textnormal{th}}$ elementary symmetric polynomial. The Lagrangian submanifolds we construct in Section~\ref{section:quotientFC} are with respect to this structure. Note that $\Sym^d(\C)$ also inherits the natural symplectic structure from an area form on $\C$, as prescribed in Section~\ref{sec:Fuksym}. For $\Sigma=\C$, this turns out to be equivalent to the standard form on $\C^d$, as they both tame the complex structure $\Sym^d(J)$ induced by the standard complex structure $J$ on $\C$ \cite[Remark~2.1.1]{DJL}. The same fact holds for $\Sigma=\D$, under the identification of $\Sym^d(\D)$ with a $2d$-dimensional ball.

\subsection{Partially wrapped Fukaya categories}\label{sec:pwFc}
Throughout the article we make reference to three different versions of ``partially wrapped Fukaya categories''. Firstly, $\W_n^d$ and, more generally, partially wrapped Fukaya categories $\W^{\textnormal{Aur}}(\Sym^d(\Sigma), \Lambda^{(d)})$ of symmetric products of surfaces follow the Floer-theoretic setup of Auroux \cite{Auroux2}. Secondly, given a \emph{stopped Liouville manifold} $(X,\mathfrak{f})$ (as in \cite{GPS2}) we denote by $\W(X,\mathfrak{f})$ the partially wrapped Fukaya category defined in \cite{Sylvan} and \cite{GPS1, GPS2}. Auroux' setup of partially wrapped Fukaya categories can be seen a precursor of the more general framework developed by Sylvan and Ganatra--Pardon--Shende, who laid the technical foundations for partially wrapped Fukaya categories for \emph{Liouville sectors}, as in \cite[Definition~1.1]{GPS1}. However, unless otherwise specified, our definition of partially wrapped Fukaya categories of symmetric products of surfaces follows that of Auroux \cite{Auroux1}. To emphasize this, we have chosen to decorate Auroux' partially wrapped Fukaya categories with the superscript ``Aur'', so to avoid any confusion in notation with the  partially wrapped Fukaya categories of Sylvan and Ganatra--Pardon--Shende. Thirdly, we write $\textnormal{FS}(f)$ for the partially wrapped Fukaya category associated to a Lefschetz fibration $f$, defined in \cite[Section~8.6]{GPS2}. In the remainder of this section we give a brief overview of the definition of $\textnormal{FS}(f)$ and of its relation with the Fukaya--Seidel category $\F(f)$.

Following the setup of \cite{GPS2}, we define a \emph{Liouville pair} $(X^{\circ},F^{\circ})$ to consist of a Liouville manifold $X^{\circ}$ and of the completion $F^{\circ}$ of a Liouville hypersurface $F$ of the boundary at infinity of $X^{\circ}$. In this setting, we define the \emph{core} of $F^{\circ}$ to be
\begin{equation*}
    \mathfrak{f}:=\Core(F):=\bigcap_{t>0} Z^{-t}(F),
\end{equation*}
the set of points which do not escape to the boundary under the flow of the Liouville vector field $Z$ \cite[Section~1.1]{GPS2}. We also define the \emph{relative core} to be
\begin{equation*}
    \Core(X^{\circ}, F^{\circ}):=\Core(X^{\circ}) \cup (\mathfrak{f} \times \CR),
\end{equation*}
where $\mathfrak{f} \times \CR \subset X^{\circ}$ is the locus of points which do not escape to the complement of $F^{\circ}$ under positive Liouville flow. In particular, $\Core(X^{\circ}, F^{\circ})$ can be identified with the cone, under positive Liouville flow, of $F^{\circ}$.

We will only use the setup of Ganatra--Pardon--Shende in Section~\ref{sec:appendix}. Our definition of Fukaya\textendash Seidel categories of Lefschetz fibrations and of isolated singularities follows that of Seidel \cite[Chapter~3]{SeidelBk}. There are many existing definitions in the literature of ``Fukaya--Seidel categories'' (\cite{GPS2, MaydanskiySeidel, GP17}), which are technically different, although often conjectured to be equivalent. In particular, under the theory developed by Ganatra--Pardon--Shende (see \cite[Section~8.6]{GPS2}), it is expected that Fukaya\textendash Seidel categories of Lefschetz fibrations can be described as partially wrapped Fukaya categories of Liouville sectors. In this setting, the regular fibre $F=f^{-1}(*)$ of a Lefschetz fibration $f:X \to \D$ is a Liouville hypersurface of the boundary at infinity of $X$. The authors of \cite{GPS2} consider a Liouville embedding $\sigma: F \times \C_{\Real \geq 0}\hookrightarrow X$, and the core ${\mathfrak{f}=\Core(F)\times \{i\infty\}\subset \partial_{\infty}( F \times \C_{\Real \geq 0})}$. Proposition~8.20 of \cite{GPS2} states that there is an equivalence of pre-triangulated $A_{\infty}$-categories
\begin{equation*}
    \textnormal{FS}(f):=\W(X,\sigma) \xrightarrow{\simeq} \W(X,\mathfrak{f}),
\end{equation*}
where the left-hand side is Ganatra--Pardon--Shende's definition of Fukaya--Seidel categories.

\section{The derived equivalence \texorpdfstring{$\F(\f_{n,d}) \simeq \perf(A_{n,d})$}{FSf perfA}}\label{section:2}
We prove Proposition~\ref{thm:introthm1} in Section~\ref{section:BP}. We use the constructive description of $\F(\fhat_{n,d})$ to prove Theorem~\ref{thm:introthm2} in Section~\ref{section:quotientFC}. In Section~\ref{section:F-W}, we prove Theorem~\ref{thm:introthm3}.

\subsection{An invariant description of \texorpdfstring{$\F(\fhat_{n,d})$}{Ffhat}}\label{section:BP}
Brieskorn\textendash Pham singularities have been extensively studied in the context of singularity theory and symplectic topology \cite{Pham, Milnor68, SebastianiThom, Gabrielov73, Okada, FU09}. It is known (for example, a combination of \cite[Theorem~1.1]{FU11}, \cite[Theorem~3.2]{KY11}, \cite[Theorem~8.6]{FZ03}, and \cite[Th\'eor\`eme]{SebastianiThom}) that, if $A^i$ ($i=1,\dots,d$) is the path algebra of the $A_n$ quiver over $\kk$ with any orientation of its arrows, then there is a quasi-equivalence of triangulated $A_{\infty}$-categories
\begin{equation*}
    \F(\fhat_{n,d}) \xrightarrow{\text{ $\simeq$ }}\perf(A^1\otimes \dots \otimes A^d).
\end{equation*}

In this section, we re-prove the above when, for each $i$, $A^i:=\Bhat_{n}$ is the path algebra of the $A_n$ quiver with alternating orientation of the arrows. In particular, we construct a Morsification of $\fhat_{n,d}$ and a collection of vanishing paths such that the endomorphism algebra of the corresponding Lefschetz thimbles is isomorphic (as an $A_{\infty}$-algebra concentrated in degree 0) to $\Ghat_{n,d}$. 

\begin{definition}\label{def:Bhatnd}
    We define $\Bhat_{n,d}$ to be the path algebra over $\kk$ of the quiver whose set of vertices is indexed by a tuple $(i_1,\dots,i_d)$ of integers between $1$ and $n$, with arrows ${(i_1,\dots,i_h,\dots,i_d) \to (i_1,\dots,i_h\pm 1,\dots,i_d)}$ whenever $i_h$ is odd, and which is equipped with all possible commutativity relations.
\end{definition}

Throughout this section, we slightly relax the usual notions of ``Morsification'' and ``distinguished collection of paths'', in the following sense (compare with \cite[Remark~6.6]{Seidel15}).
\begin{definition}\label{def:nongenericness}
    We say a Morsification of an isolated hypersurface singularity is \emph{non-generic} if it has distinct non-degenerate singularities with coinciding critical values. A \emph{non-generic} distinguished collection of vanishing paths for a non-generic Morsification is one which has coinciding vanishing paths whenever two or more singularities have coinciding critical values, and satisfies the usual conditions in \cite[Section~16d]{SeidelBk} otherwise.
\end{definition}
We note that, by construction, two Lefschetz thimbles above coinciding vanishing paths will be disjoint, so the constructed collection of thimbles will be distinguished in the usual sense.

We construct a full exceptional collection of thimbles for $\fhat_{n,d}$ inductively on $d$, using a Thom\textendash Sebastiani-type technique, based on work of Gabri\'elov \cite{Gabrielov73} (see \cite{Keating15} for the symplectic version, which directly inspired our constructions). This will be Propositions~\ref{lemma:atomicMors} and \ref{prop:Gabrielovthimbles}. We compute the endomorphism algebra of such collections using techniques developed by Futaki\textendash Ueda \cite{FU11}, see Proposition~\ref{prop:Gabungraded}. We first review the more general constructions of Gabri\'elov and Futaki\textendash Ueda.

\subsubsection{A distinguished basis for decoupled singularities}\label{sec:Gabrielov}
Let $P_0(x)$ and $Q_0(y)$ be two singularities at the origin with $x \in \C^n$ and $y \in \C^m$, and suppose we want to study the singularity $P_0(x)+Q_0(y)$. Fix $P(x)$ and $Q(y)$ respectively two generic Morsifications of $P_0(x)$ and $Q_0(y)$ and suppose $\{p_1,\dots,p_{\mu}\}$ (resp.\ $\{q_1,\dots,q_{\nu}\}$) are the non-degenerate critical values of $P$ (resp.\ $Q$), where $\mu$ (resp.\ $\nu$) is the Milnor number of $P$ (resp.\ $Q$). Fix $p_*$ (resp.\ $q_*$) a regular value of $P$ (resp.\ $Q$) and $\{\gamma_1,\dots,\gamma_{\mu}\}$ (resp.\ $\{\delta_1,\dots,\delta_{\nu}\}$) a distinguished collection of vanishing paths. The critical values of $P(x)+Q(y)$ are $\{p_i+q_j\}$, for $i=1\dots,\mu$, $j=1,\dots,\nu$. We define a path $\epsilon_{ij}$ between the regular value $p_*+q_*$ and each $p_i+q_j$ as follows. First, note that the path $p_*+\delta_j$ is a path from $p_*+q_*$ to $p_*+q_j$. Second, the path $\gamma_i+q_j$ goes from $p_*+q_j$ to $p_i+q_j$. The path $\epsilon_{ij}$ is given by the smoothening of the concatenation of these two paths. See Figure~\ref{fig:pathsdecoupled}, where for $\mu=2$ and $\nu=3$, we have sketched the critical values and choices of vanishing paths for $P$ and $Q$, the paths $p_*+\delta_j$ (pink) and $\gamma_i+q_j$ (blue), and the smoothing $\epsilon_{ij}$ of their concatenations (black).

In general, we can assume sufficient genericity on the choices of Morsifications $P_0(x)$ and $Q_0(y)$ and on the associated vanishing paths. In particular, we can assume the critical values of $P(x)+Q(y)$ to be pairwise distinct, and we can construct the paths $\{\epsilon_{ij}\}_{i,j}$ to be pairwise disjoint away from $p_*+q_*$ \cite[Assumption~2.37]{Keating15}. If this is not the case, we can appropriately rescale $Q_0(y)$ and isotope the vanishing cycles $\{\delta_1,\dots,\delta_{\nu}\}$ \cite[Section~2.3.1]{Keating15}. Under these generic assumptions, the vanishing paths $\{\epsilon_{ij}\}_{i,j}$ make up a distinguished collection (as in \cite[Section~16d]{SeidelBk}), and they are ordered lexicographically, with $\epsilon_{ij}<\epsilon_{hk}$ exactly when $i<h$ or $i=h$ and $j<k$.

\begin{figure}[t]
   \centering
   \begin{minipage}{.1 \textwidth}
    \centering
   \begin{tikzpicture}[scale=.8,v/.style={draw,shape=circle, fill=black, minimum size=1mm, inner sep=0pt, outer sep=0pt},cross/.style={cross out,draw=black, minimum size=1mm},cross/.style={cross out, draw=black, minimum size=1mm, inner sep=0pt, outer sep=0pt},cross/.default={1pt}]
    \node[v] at (3,0) {};

    \node[cross,label=above:$p_1$] at (2.5,.5) {};
    \node[cross,label=above:$p_2$] at (3,.5) {};
   
    \draw[color=ballblue] (3,0) to (2.5,.5);
    \draw[color=ballblue] (3,0) to (3,.5);
  
\end{tikzpicture}
\end{minipage}
   \begin{minipage}{.4 \textwidth}
    \centering
   \begin{tikzpicture}[scale=.8,v/.style={draw,shape=circle, fill=black, minimum size=1mm, inner sep=0pt, outer sep=0pt},cross/.style={cross out,draw=black, minimum size=1mm},cross/.style={cross out, draw=black, minimum size=1mm, inner sep=0pt, outer sep=0pt},cross/.default={1pt}]
    \node[v] at (0,0) {};

    \node[cross,label=above:$q_1$] at (3,0) {};
    \node[cross,label=above:$q_2$] at (-1,2) {};
    \node[cross,label=above:$q_3$] at (-1,-2) {};
  
    \draw[color=frenchrose] (0,0) to (3,0);
    \draw[color=frenchrose] (0,0) to (-1,2);
    \draw[color=frenchrose] (0,0) to[in=0] (-1.5,3) to[out=180,in=110] (-2,1) to[out=290,in=90] (0,-1.5) to[out=270,in=-30] (-1,-2);
  
\end{tikzpicture}
\end{minipage}
   \begin{minipage}{.4 \textwidth}
    \centering
        \begin{tikzpicture}[scale=.8,v/.style={draw,shape=circle, fill=black, minimum size=1mm, inner sep=0pt, outer sep=0pt},cross/.style={cross out,draw=black, minimum size=1mm},cross/.style={cross out, draw=black, minimum size=1mm, inner sep=0pt, outer sep=0pt},cross/.default={1pt}]
            \node[v] at (0,0) {};

            \node[cross] at (2.5,.5) {};
            \node[cross] at (3,.5) {};

            \node[cross] at (-1.5,2.5) {};
            \node[cross] at (-1,2.5) {};

            \node[cross] at (-1.5,-1.5) {};
            \node[cross] at (-1,-1.5) {};

            \draw[frenchrose] (0,0) to (3,0);
            \draw[frenchrose] (0,0) to (-1,2);
            \draw[frenchrose] (0,0) to[in=0] (-1.5,3) to[out=180,in=110] (-2,1) to[out=290,in=90] (0,-1.5) to[out=270,in=-30] (-1,-2);

            \draw[ballblue] (3,0) to (2.5,.5);
            \draw[ballblue] (3,0) to (3,.5);

            \draw[ballblue] (-1,2) to (-1.5,2.5);
            \draw[ballblue] (-1,2) to (-1,2.5);

            \draw[ballblue] (-1,-2) to (-1.5,-1.5);
            \draw[ballblue] (-1,-2) to (-1,-1.5);

            \draw[dashed] (0,0) to (2.5,.5);
            \draw[dashed] (0,0) to[out=-5,in=190] (3.2,-.25) to[out=10] (3,.5);

            \draw[dashed] (0,0) to (-1.5,2.5);
            \draw[dashed] (0,0) to (-1,2.5);

            \draw[dashed] (0,0) to[in=0] (-1.5,2.9) to[out=180,in=110] (-1.8,1) to[out=290,in=90] (.2,-1.5) to[out=270,in=-30] (-1,-2.2) to[out=150,in=-90] (-1.5,-1.5) {};

            \draw[dashed] (0,0) to[in=0] (-1.5,3.1) to[out=180,in=110] (-2.2,1) to[out=290,in=90] (-.2,-1.5) to[out=270,in=-30] (-.8,-1.9) to[out=150,in=-90] (-1,-1.5);

        \end{tikzpicture}
    \end{minipage}
        \caption{Vanishing paths for $P$ (left), $Q$ (centre), and $P+Q$ (right), for $\mu=2$, $\nu=3$.}
        \label{fig:pathsdecoupled}
    \end{figure}

In addition to the above collection of vanishing paths, \cite[Section~2]{Gabrielov73} also provides an explicit description of the associated vanishing cycles and their pairwise intersection numbers using a second fibration on the fibre above $p_*+q_*$, given by projection to $Q(y)$.

When $m=1$ and $Q_0(y)=y^k$, for some $k\geq 2$, \cite[Section~2]{FU11} computes the endomorphism algebra of a preferred collection of vanishing cycles in the fibre above $p_*+q_*$. We review a slightly modified, but equivalent, version their construction.  We refer to \cite[Section~2.5.4]{Keating15} for a similar account.

Assume that $p_*$ is small on the negative imaginary axis of $\C$, that the critical values of $P$ are real, positive and sufficiently small, and that each $\gamma_i$ is the straight segment between $p_*$ and $p_i$; whenever this is not the case, an isotopy of the base of $P$ puts us in this setting. These are the only assumptions that differ from \cite{FU11}, where the authors assume to have chosen an isotopy on the base of $P$ so that its critical values are the set of $\mu^{\textnormal{th}}$ roots of unity, and the vanishing paths are the straight segments between them and the origin (their fixed critical value). Our setting can be isotoped to the one described by Futaki\textendash Ueda. Let $Q(y)=y^k-Ny$ be a Morsification of $Q_0$ for sufficiently big $N$, $q_*=0$, and each $\delta_i$ the straight segment between the origin and $q_i$. Label the critical values of $Q$ counter-clockwise, so that $q_1$ is on the positive real axis (see Figure~\ref{fig:FUpaths}, left).
\begin{remark}\label{rk:FUlabelling}
    An inconsequential choice one can make is to label the critical values of $Q$ clockwise or counter-clockwise. We choose to follow the latter convention in this section, as this matches Futaki--Ueda's \cite{FU11}, and the one implicitly used in \cite[Section~2.5.4]{Keating15}. For the sake of clarity, we will highlight where we use this convention.
\end{remark}
We consider the Lefschetz bifibration (as in \cite[Section~15e]{SeidelBk}) given by
\begin{equation*}
\begin{tikzcd}
    \C^{n+1} \arrow[r, "\varpi"] \arrow[rr, bend right, "P+Q"'] &\C^2 \arrow[r, "\psi"]& \C
\end{tikzcd}
\end{equation*}
where $\varpi(x,y)=(P(x)+Q(y),y)$ and $\psi$ is the projection to the first coordinate. The critical values of $P$ and $Q$ are depicted in Figure~\ref{fig:FUpaths} (left, $\mu=3$ and $k=5$ in figure) and the vanishing paths detailed at the beginning of the section are as in Figure~\ref{fig:FUpaths} (centre).

\begin{figure}[t]
   \centering
     \begin{minipage}{.31 \textwidth}
          \centering
        \begin{tikzpicture}[scale=.8,v/.style={draw,shape=circle, fill=black, minimum size=1mm, inner sep=0pt, outer sep=0pt},cross/.style={cross out,draw=black, minimum size=.25mm},cross/.default={1mm}]

                \node[v,color=ballblue,label=135:$p_1$] at (.75,0) {};
                \node[v,color=ballblue,label=above:$p_2$] at (1,0) {};
                \node[v,color=ballblue,label=45:$p_3$] at (1.25,0) {};

                \node[v,color=frenchrose,label=right:$q_1$] at (2,0) {};
                \node[v,color=frenchrose,label=above:$q_2$] at (0,2) {};
                \node[v,color=frenchrose,label=left:$q_3$] at (-2,0) {};
                \node[v,color=frenchrose,label=below:$q_4$] at (0,-2) {};
        \end{tikzpicture}
    \end{minipage}
     \begin{minipage}{.31 \textwidth}
          \centering
        \begin{tikzpicture}[scale=.8,v/.style={draw,shape=circle, fill=black, minimum size=1mm, inner sep=0pt, outer sep=0pt},]

            \node at (0,-.25) {$*$};
            \node[v] at (2.75,0) {};
            \node[v] at (3,0) {};
            \node[v] at (3.25,0) {};

            \node[v] at (.75,2) {};
            \node[v] at (1,2) {};
            \node[v] at (1.25,2) {};

            \node[v] at (-1.25,0) {};
            \node[v] at (-1,0) {};
            \node[v] at (-.75,0) {};

            \node[v] at (.75,-2) {};
            \node[v] at (1,-2) {};
            \node[v] at (1.25,-2) {};

            \draw[color=ForestGreen] (0,-.25) to[bend right=20] (2.75,0);
            \draw[color=red] (0,-.25) to[bend right=30] (3,0);
            \draw[color=blue] (0,-.25) to[bend right=40] (3.25,0);

            \draw[color=ForestGreen] (0,-.25) to[bend right=10] (.75,2);
            \draw[color=red] (0,-.25) to[bend right=15] (1,2);
            \draw[color=blue] (0,-.25) to[bend right=20] (1.25,2);

            \draw[color=ForestGreen] (0,-.25) to[bend left=40] (-1.25,0);
            \draw[color=red] (0,-.25) to[bend left=30] (-1,0);
            \draw[color=blue] (0,-.25) to[bend left=20] (-.75,0);

            \draw[color=ForestGreen] (0,-.25) to[bend right=10] (.75,-2);
            \draw[color=red] (0,-.25) to[out=270,in=250,looseness=1.5] (1,-2);
            \draw[color=blue] (0,-.25) to[out=250,in=260,looseness=2] (1.25,-2);

        \end{tikzpicture}
    \end{minipage}
    \begin{minipage}{.31 \textwidth}
        \centering
      \begin{tikzpicture}[scale=.8,v/.style={draw,shape=circle, fill=black, minimum size=1mm, inner sep=0pt, outer sep=0pt},cross/.style={cross out,draw=black, minimum size=.25mm},cross/.default={1mm}]


              \node[v] at (-.75,0) {};
              \node[v] at (-1,0) {};
              \node[v] at (-1.25,0) {};


              \node[v] at (2,0) {};
              \node[v] at (0,2) {};
              \node[v] at (-2,0) {};
              \node[v] at (0,-2) {};

              \draw[color=blue] (-2,0) to[bend right=40] (-1.25,0);
              \draw[color=red] (-2,0) to[bend right=40] (-1,0);
              \draw[color=ForestGreen] (-2,0) to[bend right=40] (-.75,0);

              \draw[color=blue] (2,0) to[bend left=40] (-1.25,0);
              \draw[color=red] (2,0) to[bend left=40] (-1,0);
              \draw[color=ForestGreen] (2,0) to[bend left=40] (-.75,0);

              \draw[color=blue] (0,-2) to[bend left=20] (-1.25,0);
              \draw[color=red] (0,-2) to[bend left=15] (-1,0);
              \draw[color=ForestGreen] (0,-2) to[bend left=10] (-.75,0);
             
              \draw[color=blue] (0,2) to[out=-70,in=-80,looseness=2] (-1.25,0);
             \draw[color=red] (0,2) to[out=270,in=-70,looseness=1.25] (-1,0);
              \draw[color=ForestGreen] (0,2) to[bend left=10] (-.75,0);

      \end{tikzpicture}
  \end{minipage}
        \caption{(left) Critical values of $P$ and $Q$. (centre) Vanishing paths on the base of $P+Q$. (right) Trajectories of the critical values of $\pi_t$ along vanishing paths, colour-coded so that vanishing cycles correspond to trajectories.}
        \label{fig:FUpaths}
    \end{figure}

The vanishing cycles on the regular fibre above $p_*+q_*$ associated to the chosen vanishing paths can be obtained from a \emph{matching sphere} construction, following \cite[Section~16i]{SeidelBk}. Specifically, for $t \in \C$ a regular value for $P+Q$, we consider the Lefschetz fibration
\begin{equation*}
    \varpi_t: (P+Q)^{-1}(t) \to \psi^{-1}(t), \quad (x,y) \mapsto (t,y)
\end{equation*}
with a collection of critical values denoted by $\textnormal{Critv}(\varpi_t)$ \cite[Section~15e]{SeidelBk}. By varying $t$ on a chosen vanishing path $\epsilon_{ij}$, the critical values $\textnormal{Critv}(\varpi_t)$ trace out trajectories in the base of $\varpi_t$, identified with the complex plane for each $t$. As one approaches $p_i+q_j$, exactly two points of $\textnormal{Critv}(\varpi)$ come together at a branch point of $\psi|_{\textnormal{Critv}(\varpi)}$ \cite[Lemma~15.8]{SeidelBk}. The trajectory with which these points came together in the base of $\varpi_{p_*+q_*}$ is called a \emph{matching path}; above this, there is a matching cycle, which Seidel proves to be isotopic (as a framed Lagrangian sphere in $(P+Q)^{-1}(p_*+q_*)$) to the vanishing cycle associated to the vanishing path $\epsilon_{ij}$ \cite[Lemma~16.15]{SeidelBk}. We denote each of these matching cycles by $C_{i,j}$, and we view it as an object of the Fukaya category of $(P+Q)^{-1}(p_*+q_*)$.

In general, the above construction is quite involved. Futaki\textendash Ueda approach this computation by considering the family of maps
    \begin{equation*}
        \pi_t:=Q \circ  \varpi_t: (P+Q)^{-1}(t) \to \C, \quad (x,y) \mapsto s:=Q(y)=t-P(x),
    \end{equation*}
    whose fibres above $s$ are given by
    \begin{equation*}
        \pi^{-1}_t(s)=P^{-1}(t-s) \times Q^{-1}(s).
    \end{equation*}
    This is singular whenever $s$ is either a critical value of $Q$ or $t$ minus a critical value of $P$. As one varies $t$ along $\epsilon_{ij}$ in Figure~\ref{fig:FUpaths} (centre), the critical values of minus $P$ move on the base of $Q$ until $-p_i$ hits $q_j$. The trajectories of the critical values of $\pi_t$ along vanishing paths are depicted in Figure~\ref{fig:FUpaths} (right), and they are translations and reflections of the vanishing paths $\{\epsilon_{ij}\}_{i,j}$. Each trajectory is the image by $Q$ of the matching path $\sigma_{ij}$ corresponding to $\epsilon_{ij}$. The matching paths (obtained as the inverse images of these trajectories by $Q$) on the base of $\varpi_{p_*+q_*}$ are depicted in Figure~\ref{fig:FUmatchingpaths} (left). See also \cite[Figure~5]{FU11} and \cite[Figure~9]{Keating15}.

\begin{figure}[t]
    \centering
    \begin{minipage}{.45 \textwidth}
        \centering
    \begin{tikzpicture}[scale=.7,v/.style={draw,shape=circle, fill=black, minimum size=1mm, inner sep=0pt, outer sep=0pt},]

        \node[v,label=90:$p_{11}$] at (0:3) {};
        \node[v,label={[label distance=.3cm]90:$p_{21}$}] at (0:3.5) {};
        \node[v,label=270:$p_{31}$] at (0:4) {};

        \node[v,label=90+72:$p_{12}$] at (72:3) {};
        \node[v,label=90+72:$p_{22}$] at (72:3.5) {};
        \node[v,label=90+72:$p_{32}$] at (72:4) {};

        \node[v,label={[label distance=.2cm]90+144+20:$p_{13}$}] at (144:3) {};
        \node[v,label=90+144:$p_{23}$] at (144:3.5) {};
        \node[v,label=180:$p_{33}$] at (144:4) {};

        \node at (-1.9,-2.3) {$p_{14}$};
        \node at (-2.3,-2.7) {$p_{24}$};
        \node at (-2.7,-3.1) {$p_{34}$};

        \node[v] at (216:3) {};
        \node[v] at (216:3.5) {};
        \node[v] at (216:4) {};

        \node[v,label=0:$p_{15}$] at (288:3) {};
        \node[v,label=0:$p_{25}$] at (288:3.5) {};
        \node[v,label=0:$p_{35}$] at (288:4) {};

        \draw[color=ForestGreen] (0:3) to (72:3);
        \draw[color=red] (0:3.5) to[out=230,in=310, looseness=1.2] (72:3.5);
        \draw[color=blue] (0:4) to[out=220,in=320,looseness=1.7] (72:4);

        \draw[color=ForestGreen] (72:3) to (144:3);
        \draw[color=red] (72:3.5) to[out=230+72,in=310+72, looseness=1.2] (144:3.5);
        \draw[color=blue] (72:4) to[out=220+72,in=320+72,looseness=1.7] (144:4);

        \draw[color=ForestGreen] (144:3) to (216:3);
        \draw[color=red] (144:3.5) to[out=230+144,in=310+144, looseness=1.2] (216:3.5);
        \draw[color=blue] (144:4) to[out=220+144,in=320+144,looseness=1.7] (216:4);

        \draw[color=ForestGreen] (216:3) to (288:3);
        \draw[color=red] (216:3.5) to[out=230+216,in=310+216, looseness=1.2] (288:3.5);
        \draw[color=blue] (216:4) to[out=220+216,in=320+216,looseness=1.7] (288:4);

    \end{tikzpicture}
    \end{minipage}
    \begin{minipage}{.45 \textwidth}
        \centering
        \begin{tikzpicture}[scale=.65,v/.style={draw,shape=circle, fill=black, minimum size=1mm, inner sep=0pt, outer sep=0pt},]

            \node[v,label=120:$p_{11}$] at (0:3) {};
            \node[v,label=90:$p_{21}$] at (0:3.5) {};
            \node[v,label=60:$p_{31}$] at (0:4) {};
    
            \node[v,label=90+72:$p_{12}$] at (72:3) {};
            \node[v,label=90+72:$p_{22}$] at (72:3.5) {};
            \node[v,label=90+72:$p_{32}$] at (72:4) {};

            \node[v,label=270:$p_{13}$] at (144:3) {};
            \node[v,label=90+144:$p_{23}$] at (144:3.5) {};
            \node[v,label=180:$p_{33}$] at (144:4) {};

            \node[v,label=0:$p_{14}$] at (216:3) {};
            \node[v,label=90+216:$p_{24}$] at (216:3.5) {};
            \node[v,label=270:$p_{34}$] at (216:4) {};

            \node[v,label=0:$p_{15}$] at (288:3) {};
            \node[v,label=0:$p_{25}$] at (288:3.5) {};
            \node[v,label=0:$p_{35}$] at (288:4) {};

            \draw (0,0) to[bend right=10] (0:3);
            \draw (0,0) to[bend right=15] (0:3.5);
            \draw (0,0) to[bend right=20] (0:4);

            \draw (0,0) to[bend right=10] (72:3);
            \draw (0,0) to[bend right=15] (72:3.5);
            \draw (0,0) to[bend right=20] (72:4);

            \draw (0,0) to[bend right=10] (144:3);
            \draw (0,0) to[bend right=15] (144:3.5);
            \draw (0,0) to[bend right=20] (144:4);

            \draw (0,0) to[bend right=10] (216:3);
            \draw (0,0) to[bend right=15] (216:3.5);
            \draw (0,0) to[bend right=20] (216:4);

            \draw (0,0) to[bend right=10] (288:3);
            \draw (0,0) to[bend right=15] (288:3.5);
            \draw (0,0) to[bend right=20] (288:4);
        \end{tikzpicture}
    \end{minipage}

        \caption{(left) Matching paths on the base of $\varpi_{p_*+q_*}$, for $\mu=3$ and $k=5$. (right) Vanishing paths for $\varpi_{p_*+q_*}$ associated to each vanishing cycle $V_{i,j}$.}
    \label{fig:FUmatchingpaths}
\end{figure}

The description of the matching paths provided by Futaki--Ueda is the following. Critical points of $\varpi_{p_*+q_*}$ are solutions of the following:
\begin{equation*}
    \begin{cases}
        P(x)+Q(y)=p_*+q_* \\
        dP(x)=0 \\
        dQ(y)=c
    \end{cases}
\end{equation*}
for some $c\in\C$; in particular, for each critical point $p_i$ of $P$, there are $k$ critical values $Q^{-1}(p_*+q_* -P(p_i))=\{p_{i1},\dots,p_{ik}\}$ of $\varpi_{p_*+q_*}$. We can arrange these to be positive rescalings of the $k^{\textnormal{th}}$ roots of unity, and the matching path $\sigma_{ij}$ to connect $p_{i(j+1)}$ with $p_{ij}$. As per Remark~\ref{rk:FUlabelling}, the critical values of $\varpi_{p_*+q_*}$ are labelled counter-clockwise for increasing $k$. Denote by $C_i \subset P^{-1}(p_*)$ the vanishing cycle associated to the vanishing path $\gamma_i$ on the base of $P$. The fibre $P^{-1}(p_*)$ can be naturally identified with $\varpi_{p_*+q_*}^{-1}(0)$, so that each vanishing cycle $V_{i,j}$ of $\varpi_{p_*+q_*}$ along the line from the origin to $p_{ij}$ in Figure~\ref{fig:FUmatchingpaths} (right) corresponds to $C_i$. This gives an isomorphism of $A_{\infty}$-algebras
 \begin{equation}\label{eq:FUcorrmor}
    \hom(V_{i,j},V_{i',j'}) \cong \hom (C_i,C_{i'})
\end{equation}
 computed in the Fukaya categories $\Fuk(\varpi_{p_*+q_*}^{-1}(0))$ and $\Fuk(P^{-1}(p_*))$ respectively. In particular, for each $i,j$ we consider the morphism $e_{i,j}: V_{i,j+1} \to V_{i,j}$ corresponding to the identity $\textnormal{id}_{C_{i}}$ under
 \begin{equation*}
     \hom(V_{i,j+1},V_{i,j}) \cong \hom (C_i,C_i).
 \end{equation*}
 For $i={1,\dots,\mu}$ and $j=1,\dots,k-1$, we denote by $S_{i,j}$ the following object:
\begin{equation*}
    S_{i,j} := \Cone(V_{i,j+1} \xrightarrow{e_{i,j}}V_{i,j})
\end{equation*} 
in the directed subcategory of $\Fuk(\varpi_{p_*+q_*}^{-1}(0))$ with respect to the order
\begin{equation}\label{eq:opplexicographicFU}
   (i,j) < (i',j') \quad \textnormal{ if } j > j' \textnormal{, or } j=j' \textnormal{ and }  i<i'.
  \end{equation}
The following theorem provides a description of the directed Fukaya category of \linebreak ${(P+Q)^{-1}(p_*+q_*)}$ in terms of the Fukaya category of $P^{-1}(p_*)$, following the identification of the Lagrangian spheres in $(P+Q)^{-1}(p_*+q_*)$ with the matching cycles $C_{i,j}$. It is due to Seidel, though we use the formulation of \cite[Theorem~2.1]{FU11}.

\begin{proposition}[{\cite[Proposition~18.21]{SeidelBk}}]\label{prop:Seidelmatching}
    The Fukaya category \linebreak ${\Fuk((P+Q)^{-1}(p_*+q_*))}$ consisting of cycles $C_{i,j}$ is quasi-equivalent to full subcategory of $\Fuk(\varpi_{p_*+q_*}^{-1}(0))$ consisting of the objects $S_{i,j}$.\qed
\end{proposition}
     
\subsubsection{A distinguished basis for Brieskorn\textendash Pham}\label{sec:FukBP}
Fix $n\geq 1$. We define a Morsification $\what_{n,1}$ of $\fhat_{n,1}$. First, we choose a preliminary non-generic Morsification $\whatcirc_{n,1}$ of $\fhat_{n,1}$ satisfying the following:
\begin{enumerate}[i.]
    \item $\whatcirc_{n,1}$ is a real polynomial deformation with real and non-degenerate critical points;
    \item Each critical value of $\whatcirc_{n,1}$ is either zero or positive.
\end{enumerate}

\begin{lemma}\label{lem:preliminarywn1}
    Define $\whatcirc_{n,1}$ to be
    \begin{equation*}
        \whatcirc_{n,1}(x) :=
        \begin{cases}
            \prod(x+j)^2 \quad &\text{for $n+1$ even} \\
            (x+ \frac{n+2}{2}) \prod (x+j)^2 \quad &\text{for $n+1$ odd}
        \end{cases}
    \end{equation*}
    where the product ranges over $j=1,\dots, \frac{n+1}{2}$ (resp.\ $j=1,\dots, \frac{n}{2}$) for $n+1$ even (resp.\ odd). This is a Lefschetz fibration, and a non-generic Morsification of $\fhat_{n,1}$. Moreover, fixing a certain regular value $c$, it admits a distinguished collection of Lefschetz thimbles which are entirely contained in the real line.
\end{lemma}

\begin{proof}
    The Milnor number of $ \fhat_{n,1}$, computed as the dimension of the Jacobi algebra, is $n$. We can temporarily restrict $\whatcirc_{n,1}$ to the real line, and view it as a real polynomial. The critical points of $\whatcirc_{n,1}$ coincide with local minima and maxima (and, \emph{a priori}, saddle points) of the graph of $\whatcirc_{n,1}$ on the real plane. In fact, the local minima (global if $n+1$ is even) are obtained at the first $\frac{n+1}{2}$ (resp.\ $\frac{n}{2}$) negative integers, and local maxima are obtained at negative real numbers between them. Since there are $n$ such points, each critical point is non-degenerate (also, $\whatcirc_{n,1}$ has no saddle point). Moreover, $\whatcirc_{n,1}$ vanishes at the local minima, and is positive at the local maxima. Fix $c$ to be a positive real number, which is a regular value of $\whatcirc_{n,1}$ strictly smaller than all positive critical ones. The regular fibre above $c$ is a finite collection of points ($n+1$ of them). We can label the points in $(\whatcirc_{n,1})^{-1}(c)$ by $c^{\circ}_1,\dots,c^{\circ}_{n+1}$, with $c^{\circ}_1>\dots>c^{\circ}_{n+1}$ (each point is a negative, real number). For $j=1,\dots,n$, denote by $\Dhat^{\circ}_j$ the shortest arc in $(\whatcirc_{n,1})^{-1}(\CR)$ from $c^{\circ}_{j}$ to $c^{\circ}_{j+1}$. 
    
    Returning to the complex picture, we note that $\whatcirc_{n,1}$ is a cover of the complex plane branched at $n$ points (with coinciding, non-degenerate branch values). The Fukaya\textendash Seidel categories of branched covers of surfaces are described in \cite[Section~2A]{Seidel01b}, and the case of the disc is explicitly explained in \cite[Section~2B]{Seidel01b}. Following Seidel's description, each arc $\Dhat^{\circ}_j$ represents a Lefschetz thimble for the Fukaya--Seidel category of $\whatcirc_{n,1}$.
\end{proof}

\begin{lemma}\label{lem:disjointnessDhat}
    Fix $c$ to be a regular value of $\whatcirc_{n,1}$ which is smaller than all positive critical values, as in the proof of Lemma~\ref{lem:preliminarywn1}. For $j=1,\dots,n$, let $\Dhat^{\circ}_j$ be the Lefschetz thimble associated to the shortest arc in $(\whatcirc_{n,1})^{-1}(\CR)$ from $c^{\circ}_{j}$ to $c^{\circ}_{j+1}$. Then $\Dhat^{\circ}_j$ and $\Dhat^{\circ}_{j'}$ are disjoint whenever $j,j'$ are both even or both odd.
\end{lemma}

\begin{proof}
    By construction in the proof of Lemma~\ref{lem:preliminarywn1}, $\Dhat^{\circ}_j$ and $\Dhat^{\circ}_{j'}$ are not disjoint only when $j'=j \pm 1$.
\end{proof}

By Lemma~\ref{lem:disjointnessDhat}, we can apply a further perturbation to $\whatcirc_{n,1}$ (see Remark~\ref{rk:furtherperturbation}) to separate the critical values. In particular, the relative order after perturbation of the critical values that are strictly positive before perturbing is inconsequential. Similarly, so is the relative order after perturbation of the critical values that vanish before perturbing.

\begin{definition}
    We denote by $\what_{n,1}$ the Morsification of $\fhat_{n,1}$ obtained by applying a small perturbation to $\whatcirc_{n,1}$ so that, for $j=1,\dots,n$, each critical value $\hat{w}_j$ is real and non-negative, and given by
    \begin{equation*}
        \hat{w}_j=
        \begin{cases*}
            \frac{j-1}{2}\epsilon &\quad \text{for $j$ odd} \\
            \frac{1}{2}+\floor{\frac{n-j}{2}}\epsilon &\quad \text{for $j$ even} \\ 
        \end{cases*}
    \end{equation*}
    for some small $0 < \epsilon < \frac{1}{n-2}$. With slight abuse of notation, we denote by $\Dhat^{\circ}_{j}$ the thimble associated to the critical value $\hat{w}_j$, obtained from that in the proof of Lemma~\ref{lem:preliminarywn1} under compactly supported Hamiltonian isotopy.
\end{definition}
We can now fix the regular value for $\what_{n,1}$ to be $*=-i$, and use symplectic parallel transport along the segment $c(1-t)-it$ ($t$ varying in the unit segment) on the thimbles $\{\Dhat^{\circ}_{j}\}_j$ constructed in the proof of Lemma~\ref{lem:preliminarywn1}.
    \begin{definition}\label{def:BFShatn1}
        For each $j=1,\dots,n$, we denote by $\hat{\lambda}_{\hat{w}_j}$ the straight segment from $-i$ to $\hat{w}_j$. We define $\lambdaboldhat:=\lambdaboldhat^{\bm{1}}$ to be the collection of such vanishing paths. We denote by $\SODhat$ the associated full collection of Lefschetz thimbles in $\F(\fhat_{n,1})$. For $h=0,1$, denote by $\SODhati_{h}$ the strictly full triangulated subcategories generated by the subcollections of objects in $\SODhat$ corresponding to the critical values $\hat{w}_j$, with $h \leq 2\hat{w} < h+1$. We write $\Dhat_j$ for the symplectic parallel transported $\Dhat^{\circ}_{j}$, and $x^{\textnormal{crit}}_{j}$ for the critical point of $\what_{n,1}$ corresponding to the thimble $\Dhat_j$. We write $c_j$ for the point in the regular fibre $\what_{n,1}^{-1}(-i)$ obtained from $c^{\circ}_j$ under compactly supported Hamiltonian isotopy and symplectic parallel transport. We define $\BFShat_{n,1}$ to be the endomorphism algebra of $\SODhat$ over $\kk$. 
    \end{definition}

\begin{proposition}\label{lemma:atomicMors}
    There is a suitable choice of gradings for the thimbles underlying $\SODhat$ such that $\BFShat_{n,1}$ is isomorphic (as a graded $\kk$-algebra) to $\Bhat_{n,1}$. In particular, $\SODhat$ is a strong full exceptional collection, and $\F(\fhat_{n,1})$ admits a semi-orthogonal decomposition ${\F(\fhat_{n,1})=\langle \SODhati_0,\SODhati_1 \rangle}$.
\end{proposition}

\begin{proof}
    $\BFShat_{n,1}$ is isomorphic to $\Ghat_{n,1}$ (as an ungraded $\kk$-algebra) by construction. Indeed, the isomorphism is given by the one-to-one correspondence between the thimbles $\Dhat_{j}$, for $j$ odd (resp.\ even), and the sources (resp.\ sinks) of the alternating $A_{n}$ quiver underlying $\Bhat_{n,1}$. Moreover, for $j\neq j'$, the morphism space $\hom_{\BFShat_{n,1}}(\Dhat_j,\Dhat_{j'})$ is trivial unless $j$ is odd and $j'=j\pm 1$, in which case it is one-dimensional and generated by the unique intersection $c_{j}\in \Dhat_j \cap \Dhat_{j-1}$ or $c_{j+1}\in \Dhat_j \cap \Dhat_{j+1}$. In particular, all differentials trivially vanish. This gives a description of $\Hom_{\BFShat_{n,1}}(\Dhat_j,\Dhat_{j'})$ at a cohomological level, which is itself trivial and one-dimensional exactly when $j$ is odd and $j'=j\pm 1$. Note that all compositions vanish, since there are no holomorphic triangles (we recall that in the dictionary provided in \cite[Section~2A]{Seidel01b}, a pair of pants product corresponds to a triple intersection). Equip the Lefschetz thimbles in $\SODhat$ with arbitrary grading structure. Since the quiver underlying their endomorphism algebra is a tree, we can perform a suitable shift to each object so that all morphisms lie in degree 0. Higher $A_{\infty}$-products in $\BFShat_{n,1}$ vanish for degree reasons.
\end{proof}

Fix $n,d \geq 1$. Following the constructions reviewed in Section~\ref{sec:Gabrielov}, we view  $\fhat_{n,d}$ as the sum of decoupled singularities:
\begin{equation*}
    \fhat_{n,d}(x_1,\dots,x_d)=\fhat_{n,d-1}(x_1,\dots,x_{d-1})+\fhat_{n,1}(x_d).
\end{equation*}

\begin{definition}\label{def:whatnd}
    Define $\what_{n,d}$ to be
    \begin{equation*}
        \what_{n,d}(x_1,\dots x_d):= \sum_{j=1}^d \what_{n,1}(x_j).
    \end{equation*}
   Its critical values are of the form $\hat{w}=\hat{w}_{i_1}+\dots+\hat{w}_{i_d}$. For each such $\hat{w}$, define $\ghat_{\hat{w}}$ the straight segment between $-id$ and $\hat{w}$. For $h=0,\dots,d$, define by $\lambdaboldhat_h$ the collection of those paths, for which the real endpoint satisfies $h\leq 2\hat{w}<h+1$. We denote by $\lambdaboldhat:=\lambdaboldhat^{\bm{d}}$ the collection of all these paths. 
\end{definition}

\begin{lemma}
    The map $\what_{n,d}$ is a non-generic Morsification of $\fhat_{n,d}$, and $\lambdaboldhat$ is a non-generic distinguished collection of vanishing paths.
\end{lemma}

\begin{proof}
    Critical points of $\what_{n,d}$ are of the form $(x^{\textnormal{crit}}_{i_1},\dots,x^{\textnormal{crit}}_{i_d})$, where each $x^{\textnormal{crit}}_{j}$ is a critical point of $\what_{n,1}$. They are non-degenerate, as $\what_{n,1}$ has non-degenerate singularities at those points. The collection $\lambdaboldhat$ is a non-generic one by Definition~\ref{def:nongenericness}.
\end{proof}

\begin{definition}\label{def:BFShatnd}
    Let $*=-id$ be fixed regular value for $\what_{n,d}$, one of the endpoint of each path underlying the collection $\lambdaboldhat$ in Definition~\ref{def:whatnd}. If $x^{\textnormal{crit}}_{i}$ is the critical point of $\what_{n,1}$ corresponding to the thimble $\Dhat_{i}$ ($i=1,\dots,n$), we denote by $\Dhat_I$ the Lefschetz thimble corresponding to the critical point $x_I:=(x^{\textnormal{crit}}_{i_1},\dots,x^{\textnormal{crit}}_{i_d})$ and critical value ${\hat{w}_I:=\hat{w}_{i_1}+\dots+\hat{w}_{i_d}}$. For $h=0,\dots,d$, we denote by $\SODhati_h$ the strictly full triangulated subcategory of $\F(\fhat_{n,d})$ generated by the subcollection of thimbles associated to $\lambdaboldhat_h$. We define $\SODhat$ to be the full collection of thimbles for $\F(\fhat_{n,d})$ associated to $\lambdaboldhat$ equipped with arbitrary grading structure, and we denote its endomorphism algebra over $\kk$ by $\BFShat_{n,d}$.
\end{definition}

The index $I$ appearing in Definition~\ref{def:BFShatnd} is indexed by the following set.

\begin{definition}\label{def:Ihatnd}
    For $n$ and $d$ positive integers, we define the poset
\begin{equation*}
    \Ihat_{n,d}:=\left\{ I=(i_1,\dots,i_d)\in \mathbb{N}^d \mid 1\leq i_1,\dots,i_d\leq n \right\},
\end{equation*}
ordered by $I < J$ if $\text{ev}_{I} < \text{ev}_{J}$, where $\text{ev}_{I}$ (resp.\ $\text{ev}_{J}$) denotes the number of even entries of the tuple $I$ (resp.\ $J$).
\end{definition}

We depicted the collection of vanishing paths in Definition~\ref{def:BFShatnd} in Figure \ref{fig:Gabrielovpaths}, for \linebreak ${d+1=n=3}$. In figure, $\{p_1=\hat{w}_1,p_2=\hat{w}_3,p_3=\hat{w}_2\}$ and $\{q_1=\hat{w}_1,q_2=\hat{w}_3,q_3=\hat{w}_2\}$ denote the critical values of $\what_{n,d-1}$ and $\what_{n,1}$ respectively; in blue and pink, we depicted the concatenation of the inductively-constructed vanishing paths $\lambdaboldhat^{\bm{d-1}}$ and $\lambdaboldhat^{\bm{1}}$, whose smoothening $\lambdaboldhat^{\bm{d}}$ are the black lines.

\begin{figure}[t]
    \centering
    \begin{tikzpicture}[v/.style={draw,shape=circle, fill=black, minimum size=1mm, inner sep=0pt, outer sep=0pt},cross/.style={cross out,draw=black, minimum size=1mm}]
        \node[v,label={180:$p_*+q_*$}] at (0,1) {};
            \node[v,label={180:$p_*+q_1$}] at (0,2) {};
            \node[v,label={0:$p_*+q_2$}] at (1,2) {};
            \node[v,label={-5:$p_*+q_3$}] at (5,2) {};
    
            \node[cross,label={135:$p_1+q_1$}] at (0,3) {};
            \node[cross,label={[align=center]above:$p_1+q_2$\\$p_2+q_1$}] at (1,3) {};
            \node[cross,label={45:$p_2+q_2$}] at (2,3) {};
            \node[cross,label={[align=center]above:$p_1+q_3$\\$p_3+q_1$}] at (5,3) {};
            \node[cross,label={[align=center]45:$p_2+q_3$\\$p_3+q_2$}] at (6,3) {};
            \node[cross,label={90:$p_3+q_3$}] at (10,3) {};
    
            \draw[frenchrose,line width=1.25pt] (0,1) to (0,2);
            \draw[frenchrose,line width=1.25pt] (0,1) to (1,2);
            \draw[frenchrose,line width=1.25pt] (0,1) to (5,2);
    
            \draw[ballblue,line width=1.25pt] (0,2) to (0,3);
            \draw[ballblue,line width=1.25pt] (0,2) to (1,3);
            \draw[ballblue,line width=1.25pt] (0,2) to (5,3);
            \draw[ballblue,line width=1.25pt] (1,2) to (1,3);
            \draw[ballblue,line width=1.25pt] (1,2) to (2,3);
            \draw[ballblue,line width=1.25pt] (1,2) to (6,3);
            \draw[ballblue,line width=1.25pt] (5,2) to (5,3);
            \draw[ballblue,line width=1.25pt] (5,2) to (6,3);
            \draw[ballblue,line width=1.25pt] (5,2) to (10,3);
    
            \draw[] (0,1) to (0,3);
            \draw[] (0,1) to (1,3);
            \draw[] (0,1) to (2,3);
            \draw[] (0,1) to (5,3);
            \draw[] (0,1) to (6,3);
            \draw[] (0,1) to (10,3);
    \end{tikzpicture}

    \caption{Non-generic vanishing paths (black) for Brieskorn\textendash Pham singularities, as smoothing of the blue segments, for $d+1=n=3$.}
    \label{fig:Gabrielovpaths}
\end{figure}

\begin{proposition}\label{prop:Gabungraded}
    There is an appropriate choice of gradings of the distinguished collection of objects underlying $\BFShat_{n,d}$ such that this is an $A_{\infty}$-algebra concentrated in degree 0, and it is isomorphic to $\Bhat_{n,d}$ as a graded $\kk$-algebra.
\end{proposition}

\begin{proof}
    Our proof is an adaptation of \cite[Section~3]{FU11}. We prove the claim inductively on $d$, assuming that $\BFShat_{n,d-1}$ and $\Bhat_{n,d-1}$ are isomorphic as graded $\kk$-algebras. For $I\in \Ihat_{n,d-1}$, denote by $C_I \subset \what_{n,d-1}^{-1}(-i(d-1))$ the vanishing cycle associated to the Lefschetz thimble $\Dhat_I$ in Definition~\ref{def:whatnd}. We view $\what_{n,d}$ as a sum of decoupled singularities ${\what_{n,d}(x,y)=\what_{n,d-1}(x)+\what_{n,1}(y)}$, and $\what_{n,1}$ as a branched cover of $\C$. The endomorphism algebra  $\BFShat_{n,d}$ is independent of choices of representative of isotopy classes of thimbles and of the corresponding vanishing cycles. We can therefore isotope the base of $\what_{n,1}$ and use results of Futaki\textendash Ueda recalled in Section~\ref{sec:Gabrielov}. We first give a description of a collection of vanishing cycles for $\what_{n,d}$.

    We isotope the base of $\what_{n,1}$ so that its critical values coincide with the $n^{\textnormal{th}}$ roots of unity, and the vanishing paths with the straight segments from these to the origin. Specifically, the isotopy rotates counter-clockwise and rescales the vanishing paths (and associated critical values), so that the new critical value associated to $\ghat_{\hat{w}_1}$ is the positive real root of unity. This is shown in Figure~\ref{fig:isotopybase} (for $n=6$ in figure), where we depicted in black the original critical values and vanishing paths for $\what_{n,1}$, in red the isotopy, and in blue the final values and paths. We also labelled the isotoped critical values by $q_i$, for $i=1,\dots,n$, which follows the labelling of the critical values of $Q=\what_{n,1}$ in Remark~\ref{rk:FUlabelling}. The endomorphism algebra of the thimbles associated to these paths is the alternating $A_n$ quiver by Proposition~\ref{lemma:atomicMors}. 

    \begin{figure}[t]
        \centering
        \centering
        \begin{minipage}{.45 \textwidth}
            \centering
            \begin{tikzpicture}[scale=.65,cross/.style={cross out,draw=black, minimum size=1mm},cross/.style={cross out, draw=black, minimum size=1.5mm, inner sep=0pt, outer sep=0pt},cross/.default={1pt}]
                \node[cross] at (-.5,0) {};
                \node at (-.7,.4){$\hat{w}_1$};
                \node[cross] at (0,0) {};
                \node at (0,.4){$\hat{w}_3$};
                \node[cross] at (.5,0) {};
                \node at (.7,.4){$\hat{w}_5$};

                \node[cross] at (2,0) {};
                \node at (1.8,.4){$\hat{w}_6$};
                \node[cross] at (2.5,0) {};
                \node at (2.5,.4){$\hat{w}_4$};
                \node[cross] at (3,0) {};
                \node at (3.2,.4){$\hat{w}_2$};
             
                \node at (-.5,-1) {$*$};

                \draw (-.5,-1) to (-.5,0);
                \draw (-.5,-1) to (0,0);
                \draw (-.5,-1) to (.5,0);
                \draw (-.5,-1) to (2,0);
                \draw (-.5,-1) to (2.5,0);
                \draw (-.5,-1) to (3,0);

                \node[cross,blue] at (0:4){};
                \node[cross,blue] at (60:4){};
                \node[cross,blue] at (120:4){};
                \node[cross,blue] at (180:4){};
                \node[cross,blue] at (240:4){};
                \node[cross,blue] at (300:4){};
    
                \draw[red,->,shorten >= 2mm,shorten <= .5mm](-.5,0) to[out=180,in=90] (-1,-.5) to[out=270,in=210] (0:4);
                \draw[red,->,shorten >= 2mm,shorten <= .5mm](0,0) to[out=100,in=90,looseness=1.5] (-1.5,0) to[out=270,in=150] (300:4);
                \draw[red,->,shorten >= 2mm,shorten <= .5mm](.5,0) to[out=100,in=90,looseness=1.5] (-2,0) to[out=270,in=70] (240:4);
                \draw[red,->,shorten >= 2mm,shorten <= .5mm](2,0) to[out=120,in=45,looseness=1] (180:4);
                \draw[red,->,shorten >= 2mm,shorten <= .5mm](2.5,0) to[out=100,in=0,looseness=1] (120:4);
                \draw[red,->,shorten >= 2mm,shorten <= .5mm](3,0) to[out=60,in=-45,looseness=1] (60:4);
            \end{tikzpicture}
        \end{minipage}
        \begin{minipage}{.45 \textwidth}
            \centering
            \begin{tikzpicture}[scale=.5,cross/.style={cross out,draw=black, minimum size=1mm},cross/.style={cross out, draw=black, minimum size=1.5mm, inner sep=0pt, outer sep=0pt},cross/.default={1pt}]
    
            \node[blue] at (0,0) {$*$};
            \node[cross,blue,label={[blue]above:$\hat{w}_1=q_1$}] at (0:4){};
            \node[cross,blue,label={[blue]above:$\hat{w}_2=q_2$}] at (60:4){};
            \node[cross,blue,label={[blue]above:$\hat{w}_4=q_3$}] at (120:4){};
            \node[cross,blue,label={[blue]above:$\hat{w}_6=q_4$}] at (180:4){};
            \node[cross,blue,label={[blue]below:$\hat{w}_5=q_5$}] at (240:4){};
            \node[cross,blue,label={[blue]below:$\hat{w}_3=q_6$}] at (300:4){};
    
            \draw[blue] (0,0) to  (0:4);
            \draw[blue] (0,0) to  (60:4);
            \draw[blue] (0,0) to (120:4);
            \draw[blue] (0,0) to (180:4);
            \draw[blue] (0,0) to (240:4);
            \draw[blue] (0,0) to (300:4);
        
            \end{tikzpicture}
        \end{minipage}
            \caption{Isotopy on the base of $\what_{n,1}$.}
        \label{fig:isotopybase}
    \end{figure}

    We apply the construction detailed in Section~\ref{sec:Gabrielov}, using the vanishing paths in Figure~\ref{fig:isotopybase} (right) for the Morsification of $\what_{n,1}$. We consider the Lefschetz fibration
    \begin{equation*}
        \varpi_t: \what_{n,d}^{-1}(t)\to \C, \quad (x_1,\dots,x_d)\mapsto x_d,
    \end{equation*}
    and we denote by $V_{I,j}$ the vanishing cycle in $\Fuk(\varpi_{-id}^{-1}(0))$ corresponding to the vanishing path $\ghat_{\hat{w}}$ of the critical value $\hat{w}=\what_{n,d}(x_{I,j})$. We denote by $\sigma_{I,j}$ its associated matching path. We depicted the matching paths in Figure~\ref{fig:FUalternatingmpaths}, for some integers $n,d$ so that $\Ihat_{n,d-1}$ has size 3: in figure, we drew each $p_{X,1}$ ($X\in \{I,J,K\}=\Ihat_{n,d-1}$) as the positive real rescaled root of unity, and the remaining $p_{X,j}$ are the other rescaled roots of unity, ordered counter-clockwise for increasing $j$. Each matching path $\sigma_{X,1}$ connects $p_{X,2}$ to $p_{X,1}$, while for $j$ even (resp.\ $j\neq1$ odd), each path $\sigma_{X,j}$ connects $p_{X,\frac{2n+4-j}{2}}$ to $p_{X,\frac{j+2}{2}}$ (resp.\ $p_{X,\frac{2n+5-j}{2}}$ to $p_{X,\frac{j+3}{2}}$).  

    \begin{figure}[t]
        \centering
         \begin{minipage}{.45 \textwidth}
            \begin{tikzpicture}[scale=.7,v/.style={draw,shape=circle, fill=black, minimum size=1mm, inner sep=0pt, outer sep=0pt},]
   \node[v] at (0:3) {};
                   \node[v] at (0:3.5) {};
                   \node[v] at (0:4) {};
           
                   \node[v] at (72:3) {};
                   \node[v] at (72:3.5) {};
                   \node[v] at (72:4) {};
           
                   \node[v] at (144:3) {};
                   \node[v] at (144:3.5) {};
                   \node[v] at (144:4) {};
           
                   \node[v] at (216:3) {};
                   \node[v] at (216:3.5) {};
                   \node[v] at (216:4) {};
           
                   \node[v] at (288:3) {};
                   \node[v] at (288:3.5) {};
                   \node[v] at (288:4) {};
           
                   \draw[color=ForestGreen] (0:3) to (72:3);
                   \draw[color=red] (0:3.5) to[out=230,in=310, looseness=1.2] (72:3.5);
                   \draw[color=blue] (0:4) to[out=220,in=320,looseness=1.7] (72:4);
           
                  \draw[color=ForestGreen] (72:3) to (288:3);
                  \draw[color=red] (72:3.5) to[out=290,in=140] (288:3.5);
                  \draw[color=blue] (72:4) to[out=300,in=160,looseness=1] (288:4);
           
                   \draw[color=ForestGreen] (288:3) to (144:3);
                   \draw[color=red] (288:3.5) to[out=160,in=20, looseness=1] (144:3.5);
                   \draw[color=blue] (288:4) to[out=180,in=30,looseness=1.2] (144:4);
           
                   \draw[color=ForestGreen] (144:3) to (216:3);
                   \draw[color=red] (144:3.5) to[out=0,in=110, looseness=1] (216:3.5);
                   \draw[color=blue] (144:4) to[out=20,in=30,looseness=1.1] (190:2.5) to[out=210,in=110,looseness=1.1] (216:4);
       
                   \node[color=ForestGreen] at (2.9,2) {$\sigma_{I,1}$};
                       \draw[->,densely dotted,thick,color=ForestGreen] (3,1.75) to (2.75,.5);
                   \node[color=red] at (3.5,2.5) {$\sigma_{J,1}$};
                       \draw[->,densely dotted,thick,color=red] (3.5,2.25) to (3.25,-.13);
                   \node[color=blue] at (4.25,2) {$\sigma_{K,1}$};
                       \draw[->,densely dotted,thick,color=blue] (4,1.75) to (3.75,-.13);

                   \node[color=ForestGreen] at (-.2,2) {$\sigma_{I,2}$};
                       \draw[->,densely dotted,thick,color=ForestGreen] (.2,1.8) to (.75,1);
                   \node[color=red] at (-.2,2.5) {$\sigma_{J,2}$};
                       \draw[->,densely dotted,thick,color=red] (.2,2.3) to (1,1);
                   \node[color=blue] at (-.2,3) {$\sigma_{K,2}$};
                       \draw[->,densely dotted,thick,color=blue] (.2,2.8) to (1.25,1.15);
       
                   \node[color=ForestGreen] at (1.75,-2) {$\sigma_{I,3}$};
                       \draw[->,densely dotted,thick,color=ForestGreen] (1.1,-1.9) to (.1,-1.5);
                   \node[color=red] at (1.75,-2.5) {$\sigma_{J,3}$};
                       \draw[->,densely dotted,thick,color=red] (1.1,-2.4) to (-.2,-1.8);
                   \node[color=blue] at (1.75,-3) {$\sigma_{K,3}$};
                       \draw[->,densely dotted,thick,color=blue] (1.1,-2.9) to (-.7,-1.9);
           
                   \node[color=ForestGreen] at (-3.5,1) {$\sigma_{I,4}$};
                       \draw[->,densely dotted,thick,color=ForestGreen] (-2.8,1) to (-2.5,1);
                   \node[color=red] at (-3.5,.5) {$\sigma_{J,4}$};
                       \draw[->,densely dotted,thick,color=red] (-2.8,.5) to (-2.6,.5);
                   \node[color=blue] at (-3.5,0) {$\sigma_{K,4}$};
                       \draw[->,densely dotted,thick,color=blue] (-2.8,0) to (-2.2,0);
            \end{tikzpicture}
             \end{minipage}
             \begin{minipage}{.45 \textwidth}
                \begin{tikzpicture}[scale=.6,v/.style={draw,shape=circle, fill=black, minimum size=1mm, inner sep=0pt, outer sep=0pt},]

                    \node[v] at (0:3) {};
                    \node[v] at (0:3.5) {};
                    \node[v] at (0:4) {};
            
                    \node[v] at (60:3) {};
                    \node[v] at (60:3.5) {};
                    \node[v] at (60:4) {};
            
                    \node[v] at (120:3) {};
                    \node[v] at (120:3.5) {};
                    \node[v] at (120:4) {};
    
                    \node[v] at (180:3) {};
                    \node[v] at (180:3.5) {};
                    \node[v] at (180:4) {};
            
                    \node[v] at (240:3) {};
                    \node[v] at (240:3.5) {};
                    \node[v] at (240:4) {};
    
                    \node[v] at (300:3) {};
                    \node[v] at (300:3.5) {};
                    \node[v] at (300:4) {};

                    \draw[color=ForestGreen] (0:3) to (60:3);
                    \draw[color=red] (0:3.5) to[out=230,in=-20, looseness=1] (60:3.5);
                    \draw[color=blue] (0:4) to[out=220,in=-10,looseness=1.2] (60:4);
            
                   \draw[color=ForestGreen] (60:3) to (300:3);
                   \draw[color=red] (60:3.5) to[out=-30,in=160] (300:3.5);
                   \draw[color=blue] (60:4) to[out=-20,in=170,looseness=1.2] (300:4);
            
                   \draw[color=ForestGreen] (300:3) to (120:3);
                   \draw[color=red] (300:3.5) to[out=180,in=60, looseness=1] (120:3.5);
                   \draw[color=blue] (300:4) to[out=190,in=-90,looseness=1] (0,0) to[out=90,in=70,looseness=1] (120:4);

                   \draw[color=ForestGreen] (120:3) to (240:3);
                   \draw[color=red] (120:3.5) to[out=0,in=110, looseness=1] (240:3.5);
                   \draw[color=blue] (120:4) to[out=10,in=120,looseness=1] (240:4);
    
                   \draw[color=ForestGreen] (240:3) to (180:3);
                   \draw[color=red] (240:3.5) to[out=130,in=-60, looseness=1] (-2.25,.25) to[out=120,in=90, looseness=1] (180:3.5);
                   \draw[color=blue] (240:4) to[out=150,in=-60,looseness=1] (-2.5,.25) to[out=120,in=60,looseness=1] (180:4);

                \node[color=ForestGreen] at (2.9,3.5) {$\sigma_{I,1}$};
                   \draw[->,densely dotted, thick, color=ForestGreen] (3,3.25) to (2.75,.5);
               \node[color=red] at (3.5,4) {$\sigma_{J,1}$};
                   \draw[->,densely dotted,thick,color=red] (3.5,3.75) to (3.25,-.13);
               \node[color=blue] at (4.5,3.5) {$\sigma_{K,1}$};
                   \draw[->,densely dotted,thick,color=blue] (4,3.25) to (3.75,-.13);
    
                \node[color=ForestGreen] at (.4,2) {$\sigma_{I,2}$};
                   \draw[->,densely dotted,thick,color=ForestGreen] (.8,1.8) to (1.4,1);
               \node[color=red] at (.4,2.5) {$\sigma_{J,2}$};
                   \draw[->,densely dotted,thick,color=red] (1,2.5) to (1.8,.5);
               \node[color=blue] at (.4,3) {$\sigma_{K,2}$};
                   \draw[->,densely dotted,thick,color=blue] (1.1,3) to (2.25,.5);
    
                   \node[color=ForestGreen] at (2.5,-1) {$\sigma_{I,3}$};
                   \draw[->,densely dotted,thick,color=ForestGreen] (1.8,-1) to (.65,-1);
               \node[color=red] at (2.5,-1.5) {$\sigma_{J,3}$};
                   \draw[->,densely dotted,thick,color=red] (1.8,-1.5) to (.3,-1.8);
               \node[color=blue] at (2.5,-2) {$\sigma_{K,3}$};
                   \draw[->,densely dotted,thick,color=blue] (1.8,-2) to (.45,-2.45);
    
                \node[color=ForestGreen] at (-2.5,2) {$\sigma_{I,4}$};
                   \draw[->,densely dotted,thick,color=ForestGreen] (-1.9,2) to (-1.65,1.5);
               \node[color=red] at (-2.6,1.5) {$\sigma_{J,4}$};
                    \draw[->,densely dotted,thick,color=red] (-1.9,1.5) to (-1.25,.75);
               \node[color=blue] at (-2.7,1) {$\sigma_{K,4}$};
                    \draw[->,densely dotted,thick,color=blue] (-2,1) to (-1.75,.75);

                \node[color=ForestGreen] at (-3.55,-1) {$\sigma_{I,5}$};
                    \draw[->,densely dotted,thick,color=ForestGreen] (-2.9,-1) to (-2.5,-1);
                \node[color=red] at (-3.55,-1.5) {$\sigma_{J,5}$};
                    \draw[->,densely dotted,thick,color=red] (-2.9,-1.5) to (-2.15,-1.75);
                \node[color=blue] at (-3.6,-2) {$\sigma_{K,5}$};
                \draw[->,densely dotted,thick,color=blue] (-2.9,-2) to (-2.55,-2);
                \end{tikzpicture}
                 \end{minipage}

            \caption{Matching paths for $\what_{n,d-1}(x)+\what_{n,1}(y)$, on the base of the second fibration (projection to $y$), for $n=4$ (left) and $n=5$ (right), where the labels for the critical values are the same as in Figure~\ref{fig:FUmatchingpaths}.
            }
        \label{fig:FUalternatingmpaths}
    \end{figure}

    For $I,j$ as above, define $S_{I,j}$ to be the following object:
    \begin{align*}
        S_{I,j}=
        \begin{cases}
            \Cone(V_{I,2}\to V_{I,1}) &\quad  \textnormal{ for $j=1$} \\
            \Cone(V_{I,\frac{2n+4-j}{2}}\to V_{I,\frac{j+2}{2}}) &\quad \textnormal{ for $j$ even} \\
            \Cone(V_{I,\frac{2n+5-j}{2}}\to V_{I,\frac{j+3}{2}}) &\quad \textnormal{ for $j\neq1$ odd} \\
        \end{cases}
    \end{align*}
    in the directed subcategory of $\Fuk(\varpi_{-id}^{-1}(0))$ with respect to \eqref{eq:opplexicographicFU}, where the cone is formed over the morphism $e:=e_{I,j}:V_{I,k}\to V_{I,j}$ (for appropriate $k$) corresponding to the identity $id_{I}:C_{I}\to C_{I}$ under \eqref{eq:FUcorrmor}.

    Denote by $\Fuk^{\rightarrow}$ the directed subcategory of $\Fuk(\varpi_{-id}^{-1}(0))$ consisting of the objects $\{ S_{I,j}\}_{I,j}$, with respect to the lexicographic order
    \begin{equation}\label{eq:lexicographicFU} 
        (I,j) < (I',j') \quad \textnormal{ if } j < j' \textnormal{, or } j=j' \textnormal{ and }  I<I'.
       \end{equation}
    Proposition~\ref{prop:Seidelmatching} shows that the directed Fukaya category of vanishing cycles, as a subcategory of $\Fuk(\what_{n,d}^{-1}(-id))$, is quasi-equivalent to $\Fuk^{\rightarrow}$.
    \begin{remark}
        We remark the difference between \eqref{eq:opplexicographicFU} and \eqref{eq:lexicographicFU}, which is the same difference between the unnumbered equations in \cite[pp.\ 444 and 445]{FU11}. The former is due to the order in Remark~\ref{rk:FUlabelling}, and prescribes which objects $S_{I,j}$ of the Fukaya category $\Fuk(\varpi_{-id}^{-1}(0))$ to consider. The latter follows the lexicographic order of the vanishing paths, and gives an order of the matching spheres as objects of the directed Fukaya subcategory of $\what_{n,d}^{-1}(-id)$.
    \end{remark}
    
    We claim that the isomorphism of graded $\kk$-algebras between $\Bhat_{n,d}$ and $\BFShat_{n,d}$ is realised by the one-to-one correspondence
    \begin{equation*}
        (I,j) \leftrightarrow S_{I,j}
    \end{equation*}
    between the vertices of the quiver underlying $\Bhat_{n,d}$ and objects of the directed subcategory $\Fuk^{\rightarrow}$. By the construction in Section~\ref{sec:Gabrielov}, each $S_{I,j}$ corresponds to a thimble underlying the full exceptional collection $\SODhat$. More precisely, $S_{I,j}$ corresponds to the critical point ${x_{I,j}=(x^{\textnormal{crit}}_{I},x^{\textnormal{crit}}_{j})=(x^{\textnormal{crit}}_{i_1},\dots,x^{\textnormal{crit}}_{i_{d-1}},x^{\textnormal{crit}}_{j})}$ of $\what_{n,d}$.

    We recall from Section~\ref{sec:mainresults} that given two vertices $(I,j) =(i_1,\dots,i_{d-1},j)$ and \linebreak ${(I',j') =(i'_1,\dots,i'_{d-1},j')}$ in the quiver with relations underlying $\Bhat_{n,d}$, there is a path from the former to the latter exactly when the following conditions hold:
    \begin{itemize}
        \item $\lvert i'_h - i_h \rvert \leq 1$ for each $h\in \{1,\dots, d-1\}$, and when equality holds the corresponding $i_h$ is odd;
        \item $\lvert j' - j \rvert \leq 1$, and when equality holds $j$ is odd.
    \end{itemize}

    For $S_{I,j}=\Cone(V_{I,y} \to V_{I,x})$ and $S_{I',j'}=\Cone(V_{I',y'} \to V_{i',x'})$ (each pair $(x,y)$ and $(x',y')$ determined by $j$ and $j'$ in the definition of $S_{I,j}$ and $S_{I',j'}$ above), one has
    \begin{equation}\label{eq:doublecpx}
        \hom_{\Fuk^{\rightarrow}}(S_{I,j},S'_{I',j'})=
        \left\{
         \begin{tikzcd}
            \hom(V_{I,x},V_{I',y'}) \arrow[r,"\bullet \circ e_{I,j}"]\arrow[d,"e_{I',j'} \circ \bullet"] &  \hom(V_{I,y},V_{I',y'})\arrow[d,"e_{I',j'} \circ \bullet"]  \\
             \hom(V_{I,x},V_{I',x'})\arrow[r,"\bullet \circ e_{I,j}"] &  \hom(V_{I,y},V_{I',x'})
        \end{tikzcd}
        \right\}
    \end{equation}
    where the right-hand side denotes the total complex of the double complex, and the morphisms are computed in the directed subcategory of $\Fuk(\varpi_{-id}^{-1}(0))$ with respect to \eqref{eq:opplexicographicFU}. If $j=j'$, the right-hand side of \eqref{eq:doublecpx} is given by
    \begin{equation*}
        \left\{
         \begin{tikzcd}
           0 \arrow[r]\arrow[d] &  \hom(C_I,C_{I'})\arrow[d,"id"]  \\
           \hom(C_I,C_{I'})\arrow[r,"id"] &  \hom(C_I,C_{I'})
        \end{tikzcd}
        \right\}
        \cong
        \hom(C_I,C_{I'})
    \end{equation*}
    with morphisms computed in the Fukaya category of the distinguished basis of vanishing cycles in $\what_{n,d-1}^{-1}(-i(d-1))$. Natural representatives of a basis of the cohomology group of this complex are given by
    \begin{equation*}
        \begin{tikzcd}
            V_{I,2} \arrow[r,"e"]\arrow[d,"x_{I,I'}"] &  V_{I,1}\arrow[d,"x_{I,I'}"]  \\
            V_{I',2}\arrow[r,"e"] & V_{I',1}
         \end{tikzcd}
         \quad \quad \quad
         \begin{tikzcd}
           V_{I,\frac{2n+4-j}{2}} \arrow[r,"e"]\arrow[d,"x_{I,I'}"] &  V_{I,\frac{j+2}{2}}\arrow[d,"x_{I,I'}"]  \\
           V_{I',\frac{2n+4-j}{2}}\arrow[r,"e"] & V_{I',\frac{j+2}{2}}
        \end{tikzcd}
        \quad \quad \quad
        \begin{tikzcd}
            V_{I,\frac{2n+5-j}{2}} \arrow[r,"e"]\arrow[d,"x_{I,I'}"] &  V_{I,\frac{j+3}{2}}\arrow[d,"x_{I,I'}"]  \\
            V_{I',\frac{2n+5-j}{2}}\arrow[r,"e"] & V_{I',\frac{j+3}{2}}
         \end{tikzcd}
    \end{equation*}
    for $j=1$, $j$ even, and $j\neq 1$ odd respectively, where $x_{I,I'}$ runs over a basis of $\hom(C_I,C_{I'})$. If $j=1$ and $j'=2$, the right-hand side of \eqref{eq:doublecpx} is given by
    \begin{equation*}
        \left\{
         \begin{tikzcd}
           0 \arrow[r]\arrow[d] &  0\arrow[d,]  \\
          0\arrow[r,] &  \hom(C_I,C_{I'})
        \end{tikzcd}
        \right\}
        \cong
        \hom(C_I,C_{I'})[-1]
    \end{equation*}
    and each
    \begin{equation*}
        \begin{tikzcd}
         & V_{I,2} \arrow[r,"e_{I,1}"]\arrow[d,"x_{I,I'}"] &  V_{I,1}  \\
          V_{I',n+1}\arrow[r,"e_{I',2}"] & V_{I',2} &
       \end{tikzcd}
    \end{equation*}
    is a representative of an element of a basis of the cohomology group. If $j\neq 1$ is odd and $j'=j\pm 1$, the right-hand side of \eqref{eq:doublecpx} is given by
    \begin{equation*}
        \left\{
         \begin{tikzcd}
           0 \arrow[r]\arrow[d] &  \hom(C_I,C_{I'})\arrow[d,"id"]  \\
           \hom(C_I,C_{I'})\arrow[r,"id"] &  \hom(C_I,C_{I'})
        \end{tikzcd}
        \right\}
        \cong
        \hom(C_I,C_{I'})
    \end{equation*}
    and each
    \begin{equation*}
        \begin{tikzcd}
          V_{I,\frac{2n+5-j}{2}} \arrow[r,"e"]\arrow[d,"x_{I,I'}"] &  V_{I,\frac{j+3}{2}}\arrow[d,"x_{I,I'}"]  \\
          V_{I',\frac{2n+4-j'}{2}}\arrow[r,"e"] & V_{I',\frac{j'+2}{2}}
       \end{tikzcd}
    \end{equation*}
    is a representative of an element of a basis of the cohomology group. In any other case, the total complex on the right-hand side of \eqref{eq:doublecpx} is contractible. The above computations show that $\hom(S_{I,j},S_{I',j'})$ is non-trivial and one-dimensional exactly when there is a path from $I$ to $I'$, and when either $j=j'$, or $\lvert j' - j \rvert = 1$ and $j$ is odd. In particular, all differentials trivially vanish, and the above provides a computation of $\Hom(S_{I,j},S_{I',j'})$ at a cohomological level, which coincides with the chain-level morphisms.
    
    The representatives of cohomology given above give a straightforward way to compute compositions of pairs of morphisms, which agree with the composition of paths in the quiver with relations underlying $\Bhat_{n,d}$. To check that the commutativity of the compositions hold, suppose $j=j'_1 \neq 1$ is odd (the case $j=j'_1 = 1$ is analogous), and $j''=j'_2=j\pm 1$. The two compositions
    \begin{equation*}
        \begin{tikzcd}
            V_{I,\frac{2n+5-j}{2}} \arrow[r,]\arrow[d] &  V_{I,\frac{j+3}{2}}\arrow[d]  \\
            V_{I'_2,\frac{2n+5-j'_1}{2}}\arrow[r]\arrow[d] & V_{I'_2,\frac{j'_1+3}{2}}\arrow[d]\\
            V_{I'',\frac{2n+4-j''}{2}}\arrow[r] & V_{I'',\frac{j''+2}{2}}
       \end{tikzcd}
       \quad \quad \quad
       \begin{tikzcd}
        V_{I,\frac{2n+5-j}{2}} \arrow[r]\arrow[d] &  V_{I,\frac{j+3}{2}}\arrow[d]  \\
        V_{I'_1,\frac{2n+4-j'_2}{2}}\arrow[r]\arrow[d]  & V_{I'_1,\frac{j'_2+2}{2}}\arrow[d] \\
        V_{I'',\frac{2n+4-j''}{2}}\arrow[r] & V_{I'',\frac{j''+2}{2}}
     \end{tikzcd}
    \end{equation*}
    give commutativity of the following square:
    \begin{equation*}
        \begin{tikzcd}
            S_{I'_1,j'_2} \arrow[r] &   S_{I'',j''}  \\
            S_{I,j}\arrow[r]\arrow[u] & S_{I'_2,j'_1}.\arrow[u]
         \end{tikzcd}
    \end{equation*}
The above computations also show that morphisms spaces are concentrated in degree 0. Indeed, a morphism between two objects listed above is given by a morphism of complexes between the two representatives of the respective bases of cohomology, and there is no such morphism in degree other than 0. Finally, higher $A_{\infty}$-products in $\BFShat_{n,d}$ vanish for degree reasons.
\end{proof}

\begin{remark}
    We note that for $d=2$ we recover the semi-orthogonal decomposition we constructed in \cite[Section~2.2]{DiDedda23}, despite the non-generic Morsification considered there not being decomposable in a sum of atomic ones.
\end{remark}

\begin{remark}
    The proof of Proposition~\ref{prop:Gabungraded} relies on the specific choices of vanishing paths made in the proof of Proposition~\ref{lemma:atomicMors}. In Figure~\ref{fig:nonvpaths} (centre), we have drawn a collection of paths arising from a different choice of distinguished basis for $\what_{n,1}$ (the one associated to the linearly oriented $A_n$ quiver in Figure~\ref{fig:nonvpaths}, left) that is not a distinguished collection, and not even a non-generic one. As pointed out in \cite[Section~2.5.2]{Keating15}, one could rescale the decoupled Morsifications (Figure~\ref{fig:nonvpaths}, right), but this would not be compatible with the computations we carry out at the end of this section.
\end{remark}

    \begin{figure}[t]
        \centering
        \begin{minipage}{.17 \textwidth}
            \centering
            \begin{tikzpicture}[scale=.65,v/.style={draw,shape=circle, fill=black, minimum size=1mm, inner sep=0pt, outer sep=0pt},cross/.style={cross out,draw=black, minimum size=1mm},cross/.style={cross out, draw=black, minimum size=1.5mm, inner sep=0pt, outer sep=0pt},cross/.default={1pt}]
                \node[v] at (0,0) {};

                \node[cross] at (0,1) {};
                \node[cross] at (.5,1) {};
                \node[cross] at (2,1) {};

                \draw (0,0) to (0,1);
                \draw (0,0) to (2,1);
                \draw (0,0) to[out=15,in=300] (2.2,1.1) to[out=120,in=20] (1.4,.9) to[out=200,in=-10] (.5,1);

            \end{tikzpicture}
        \end{minipage}
    \begin{minipage}{.4 \textwidth}
        \centering
        \begin{tikzpicture}[scale=.9,v/.style={draw,shape=circle, fill=black, minimum size=1mm, inner sep=0pt, outer sep=0pt},cross/.style={cross out,draw=black, minimum size=1mm},cross/.style={cross out, draw=black, minimum size=1.5mm, inner sep=0pt, outer sep=0pt},cross/.default={1pt}]
            \node[v] at (0,0) {};

                \node[cross] at (0,2) {};
                \node[cross] at (.5,2) {};
                \node[cross] at (2,2) {};

                \node[cross] at (1,2) {};
                \node[cross] at (2.5,2) {};
                \node[cross] at (4,2) {};

                \draw[] (0,0) to (0,2);
                \draw[] (0,0) to[out=70,in=240] (.5,1) to[out=60,in=300] (2.1,2.1) to[out=120,in=20] (1.4,1.9) to[out=200,in=-10] (.5,2);

                \draw[] (0,0) to (2,2);
                \draw[] (0,0) to (4,2);
                \draw[] (0,0) to[out=15,in=300] (4.2,2.1) to[out=120,in=20] (3.4,1.9) to[out=200,in=-10] (2.5,2);

                \draw[] (0,0) to (.5,2);
                \draw[] (0,0) to (2.5,2);

                \draw[] (0,0) to[out=35,in=200] (1.8,1.1) to[out=20,in=340]  (2.8,2.1) to[out=160,in=20] (1.9,1.9) to[out=200,in=-10] (1,2);            
        \end{tikzpicture}
    \end{minipage}
        \begin{minipage}{.4 \textwidth}
            \centering
            \begin{tikzpicture}[scale=.8,v/.style={draw,shape=circle, fill=black, minimum size=1mm, inner sep=0pt, outer sep=0pt},cross/.style={cross out,draw=black, minimum size=1mm},cross/.style={cross out, draw=black, minimum size=1.2mm, inner sep=0pt, outer sep=0pt},cross/.default={1pt}]
                \node[v] at (0,0) {};
    
                \node[v] at (0,2) {};
                \node[v] at (1,2) {};
                \node[v] at (4,2) {};

                \node[cross] at (0,2.28) {};
                \node[cross] at (.14,2.28) {};
                \node[cross] at (.56,2.28) {};

                \node[cross] at (1,2.28) {};
                \node[cross] at (1.14,2.28) {};
                \node[cross] at (1.56,2.28) {};

                \node[cross] at (4,2.28) {};
                \node[cross] at (4.14,2.28) {};
                \node[cross] at (4.56,2.28) {};

                \draw (0,0) to (0,2);
                \draw (0,0) to (4,2);
                \draw (0,0) to[out=15,in=300] (4.8,2.4) to[out=120,in=20] (2.8,1.8) to[out=200,in=-10] (1,2);

                \draw (0,2) to (0,2.28);
               \draw (0,2) to (.56,2.28);
                \draw (0,2) to[out=15,in=300] (.7,2.35) to[out=120,in=20] (.4,2.3) to[out=200,in=-10] (.14,2.28);

                \draw (1,2) to (1,2.28);
                \draw (1,2) to (1.56,2.28);
                 \draw (1,2) to[out=15,in=300] (1.7,2.35) to[out=120,in=20] (1.4,2.3) to[out=200,in=-10] (1.14,2.28);

                 \draw (4,2) to (4,2.28);
                 \draw (4,2) to (4.56,2.28);
                  \draw (4,2) to[out=15,in=300] (4.7,2.35) to[out=120,in=20] (4.4,2.3) to[out=200,in=-10] (4.14,2.28);
            \end{tikzpicture}
        \end{minipage}
        \caption{Vanishing paths for the linearly oriented $A_3 \otimes A_3$ quiver. (left) Choice of vanishing paths for which the endomorphism algebra of the collection of thimbles is the path algebra of the linearly oriented $A_3$ quiver. (centre) Non-example of vanishing paths constructed from vanishing paths on the left. (right) Admissible vanishing paths for which the endomorphism algebra of the collection of thimbles is the linearly oriented $A_3 \otimes A_3$.}
        \label{fig:nonvpaths}
    \end{figure}

    \begin{proposition}\label{prop:Gabrielovthimbles}
        $\SODhat$ is a full exceptional collection for $\F(\fhat_{n,d})$. In particular, $\F(\fhat_{n,d})$ admits a semi-orthogonal decomposition $\F(\fhat_{n,d})= \langle \SODhati_0,\dots, \SODhati_d \rangle$ of length $d+1$.
    \end{proposition}
    
    \begin{proof}
        It suffices to check that the morphism spaces between thimbles above coinciding paths vanish. We remark that this should hold in more generality for any non-generic collection of vanishing paths, see \cite[Remark 6.6]{Seidel15}. By  Proposition~\ref{prop:Gabungraded}, the isomorphism $\BFShat_{n,d}\cong \Bhat_{n,d}$ is realised by the one-to-one correspondence
        \begin{equation*}
            \Dhat_I \leftrightarrow I
        \end{equation*}
        between thimbles and vertices of the quiver underlying $\Bhat_{n,d}$ (indexed by lattice points). In particular, there is a morphism from $\Dhat_I$ to $ \Dhat_J$ exactly when there is a path from $I$ to $J$ in the quiver with relations underlying $\Bhat_{n,d}$. By Section~\ref{sec:mainresults}, this exists when $\lvert j_h - i_h \rvert \leq 1$, and equality holds for arbitrarily many $h\in \{1,\dots, d\}$, in which case the corresponding index $i_h$ is odd. Suppose there is at least one $h$ so that $i_h$ is odd, and $j_h=i_h\pm 1$. In this case, $\what_{n,1}(x^{\textnormal{crit}}_{j_h})>\what_{n,1}(x^{\textnormal{crit}}_{i_h})$, so $\what_{n,d}(x^{\textnormal{crit}}_{i_1},\dots,x^{\textnormal{crit}}_{i_d})<\what_{n,d}(x^{\textnormal{crit}}_{j_1},\dots,x^{\textnormal{crit}}_{j_d})$, and the thimbles $\Dhat_I$ to $ \Dhat_J$ lie above different vanishing paths.
    \end{proof}

\begin{theorem}\label{thm:mainthm1}
    There is a quasi-equivalence of triangulated $A_{\infty}$-categories
    \begin{equation*}
        \F(\fhat_{n,d}) \xrightarrow{\text{ $\simeq$ }} \perf (\Ghat_{n,d}),
    \end{equation*}
    induced by the isomorphism of $A_{\infty}$-algebras:
    \begin{equation}\label{eq:mainthm1}
        \BFShat_{n,d} \xrightarrow{\text{ $\cong$ }}  \Ghat_{n,d}.
    \end{equation}
\end{theorem}

\begin{proof}
    The two triangulated categories are generated by the constructed collection of thimbles and by projective $\Bhat_{n,d}$-modules respectively. $\F(\fhat_{n,d})$ admits a tilting object by Proposition~\ref{prop:Gabungraded}. In particular, $\BFShat_{n,d}$ and $\Bhat_{n,d}$ are isomorphic as $A_{\infty}$-algebras (concentrated in degree 0). The claim follows a general result of Keller \cite[Theorem~3.8]{Keller06}. We recall that this states the following. Let $\mathcal{T}$ be a given $\kk$-linear triangulated category (which is \emph{algebraic}, see \cite[Section~3.6]{Keller06}), and suppose $\mathcal{G}$ is a full subcategory. Then there is a differential graded category $\mathcal{A}$ whose cohomology is isomorphic to $\mathcal{G}$, and a functor $ \mathcal{T} \to D(\mathcal{A})$ to the derived category $D(\mathcal{A})$. Moreover, this induces an equivalence $ \mathcal{T} \to \perf(\mathcal{A})$ whenever $\mathcal{T}=\textnormal{add}(\mathcal{G})$.
\end{proof}

We note some favourable properties of the thimbles underlying \eqref{eq:mainthm1}.

\begin{definition}\label{def:index}
    We say a vertex $I$ of the quiver underlying $\Bhat_{n,d}$ in Definition~\ref{def:Bhatnd} has \emph{index} $h$ if $\text{ev}_I=h$. Under the bijection $\Dhat_I \leftrightarrow I$ underlying \eqref{eq:mainthm1}, we say $\Dhat_I$ has \emph{index} $h$ if $I$ has index $h$.
 \end{definition}

 The following is a direct consequence of the explicit description of the full collection $\SODhat$ of $\F(\fhat_{n,d})$ we constructed.
    
\begin{lemma}\label{lem:indexinvariant}
    For $I\in \Ihat_{n,d}$, $\Dhat_I$ in $\SODhat$ has index $h$ exactly when it lies in $\SODhati_h$.
\end{lemma}

\begin{proof}
    This follows from Proposition~\ref{prop:Gabrielovthimbles} and Theorem~\ref{thm:mainthm1}.
\end{proof}

The following proposition is the foundation on which Section~\ref{section:quotientFC} is built.

\begin{proposition}\label{prop:actiononthimbles}
    The symmetric group $\S_d$ acts on the collection $\SODhat$ by
    \begin{equation*}
        \sigma \cdot \Dhat_{I} = \Dhat_{\sigma \cdot I}
    \end{equation*}
    for $\sigma \in \S_d$. In particular, the $\S_d$-action is closed under restriction to each component $\SODhati_h$, and is free away from the diagonals of $\Ihat_{n,d}$.
\end{proposition}

\begin{proof}
    Since the non-generic Morsification $\what_{n,d}$ of $\fhat_{n,d}$ is invariant with respect to the natural action of $\S_d$ on $\C^d$, the latter restricts to the set of critical points. By Definition~\ref{def:BFShatnd}, this is given by
    \begin{equation*}
        \sigma \cdot x_{I} = x_{\sigma \cdot I},
    \end{equation*}
    whenever $x_I$ is a critical point, $\sigma \in \S_d$, and $I \in  \Ihat_{n,d}$. Moreover, the index of the Lefschetz thimbles in Definition~\ref{def:index} is invariant under the action of the symmetric group, fact which follows from Lemma~\ref{lem:indexinvariant} and from Definition~\ref{def:index}. Finally, the vanishing paths associated to the critical values are by construction pointwise fixed by the action of the symmetric group. The claim follows from the fact that each Lefschetz thimble is uniquely determined by symplectic parallel transport from its corresponding critical value along the given vanishing path.
\end{proof}

\subsection{The Fukaya--Seidel category \texorpdfstring{$\F(\f_{n,d})$}{Ffnd}}{\label{section:quotientFC}}
The main results of this section are Proposition~\ref{prop:exccollectionf} and Theorem~\ref{thm:mainthm2}. The significance of the former consists in the construction of a distinguished collection of generators of $\F(\f_{n,d})$ from a given collection of generators of $\F(\fhat_{n,d})$. The latter concerns the computation of the Fukaya\textendash Seidel category generated by such collection.

\begin{definition}\label{def:Bnd}
    We define $\B_{n,d}$ to be the path algebra over $\kk$ of the quiver whose set of vertices is indexed by a tuple $(i_1,\dots,i_d)$ of \emph{strictly increasing} integers between $1$ and $n$, with arrows $(i_1,\dots,i_h,\dots,i_d) \to (i_1,\dots,i_h\pm 1,\dots,i_d)$ whenever $i_h$ is odd, and which is equipped with all possible commutativity relations.
\end{definition}

Throughout this section, fix positive integers $n$ and $d$, with $n \geq d$. Denote by \linebreak ${\bpi:\C^d \to \Sym^d(\C)}$ the branched covering map, and by $\f_{n,d}$ the map satisfying \eqref{eq:liftfnd}.

\begin{lemma}
    Under the identification $\varphi: \Sym^d(\C) \xrightarrow{\cong} \C^d$ in \eqref{eq:iso}, $\f_{n,d}\circ \varphi^{-1}$ defines a complex isolated hypersurface singularity.
\end{lemma}

\begin{proof}
    By the fundamental theorem of symmetric polynomials, $\fhat_{n,d}$ has a unique representation as a polynomial in the $d$ variables $\{e_i(\mathbf{x})\}_i$. Under \eqref{eq:iso}, $\f_{n,d} \circ \varphi^{-1}$ defines a polynomial in $d$ variables. The fact that it has an isolated singularity at the origin follows from chain rule for the differential applied to \eqref{eq:liftfnd}, which in particular implies that the critical locus of $\f_{n,d}\circ \varphi^{-1}$ is strictly contained in the critical locus of $\fhat_{n,d}$.
\end{proof}

 Note that when $n<d$, $\f_{n,d}\circ \varphi^{-1}$ is smooth, and its Fukaya\textendash Seidel category is trivial. Define $\w_{n,d}$ to be the perturbation of $\f_{n,d}$ satisfying $\what_{n,d}= \w_{n,d} \circ \bpi$.  Note also that $\F(\f_{n,d})$ and $\F(\f_{n,d}\circ \varphi^{-1})$ are equivalent under \eqref{eq:iso}.

\begin{lemma}\label{lem:Milnornumber}
    The Milnor number of $\f_{n,d}$ is $ n \choose d$, $\w_{n,d}$ is a (non-generic) Morsification of $\f_{n,d}$, and its critical values are (possibly a strict subset of) critical values of $\what_{n,d}$.
\end{lemma}

\begin{proof}
    Let $\what:=\what_{n,d}$ and $\w:=\w_{n,d}$. Applying the chain rule to $\what= \w \circ \bpi$, we note that the critical locus of $\w$ is strictly contained in the image under $\bpi$ of the critical locus of $\what$, and therefore so is the set of critical values. Moreover, the Hessian matrix $\textnormal{Hess}_{\what}$ satisfies
    \begin{equation*}
        \textnormal{Hess}_{\what}(\mathbf{x})=J_{\bpi}^{\top}(\mathbf{x})\textnormal{Hess}_{\w}(\bpi(\mathbf{x}))J_{\bpi}(\mathbf{x})+\sum_{i=1}^d \frac{\partial \w}{\partial u_i} (\bpi(\mathbf{x}))\textnormal{Hess}_{\bpi}^i(\mathbf{x})
    \end{equation*}
    where $\w=\w(u_1,\dots,u_d)$, $J_{\bpi}$ denotes the Jacobian matrix of $\bpi$, $J_{\bpi}^{\top}$ its transpose, and $\textnormal{Hess}_{\bpi}(\mathbf{x})=(\textnormal{Hess}^1_{\bpi}(\mathbf{x}),\dots,\textnormal{Hess}^d_{\bpi}(\mathbf{x}))$ is the vector of Hessian matrices of each component of $\bpi$. It follows that critical points of $\w$ cannot be degenerate, since otherwise they would be images under $\bpi$ of degenerate critical points of $\what$, which is Morse. From the same formula we observe that critical points of $\what$ cannot simultaneously be singular points of $\bpi$, and have images under $\bpi$ in singular points of $\w_{n,d}$. Since points of the form $(...,\alpha,...,\alpha,...)$ are branch points of $\bpi$, their images under $\bpi$ are not critical points of $\w_{n,d}$. Away from these points, the Jacobian matrix is invertible, and the claim follows from the chain rule.
\end{proof}

We now define a full exceptional collection of thimbles for $\F(\f_{n,d})$, associated to our preferred choice of (non-generic) Morsification. Let $\hat{\delta}:=\{I\in \Ihat_{n,d}\mid i_h=i_k \text{ for some }h\neq k\}$ denote the set of diagonals of the poset $\Ihat_{n,d}$, and set $\SODhat^{\circ}:=\SODhat \setminus \{\Dhat_I \mid I\in \hat{\delta}\}$. For $h=0,\dots,d$, further denote by $\SODhati^{\circ}_h$ the subcollection of thimbles in $\F(\fhat_{n,d})$ given by ${\SODhati^{\circ}_h:=\{\Dhat_I\in \SODhati_h\} \setminus \{\Dhat_I\in \SODhati_h \mid I\in \hat{\delta}\}}$. We introduce the following notation.

\begin{definition}\label{def:Ind}
    For $n\geq d$, we define the following set.
        \begin{equation*}
            \I_{n,d}:=\left\{I=(i_1,\dots,i_d) \in \mathbb{N}^d \mid 1\leq i_1<\dots<i_d\leq n \right\}.
        \end{equation*}
    We equip $\I_{n,d}$ with a poset structure by $I < J$ if $\text{ev}_{I} < \text{ev}_{J}$ and $\text{ev}_{I}$ (resp.\ $\text{ev}_{J}$) denotes the number of even entries of the tuple $I$ (resp.\ $J$). 
    \end{definition}
    
We note that, though $\I_{n,d}$ and $\N_{n,d}$ are set-wise identical, they have different poset structures. The poset in Definition~\ref{def:Ind} indexes a collection of objects in $\F(\f_{n,d})$. More precisely, we have the following.

\begin{definition}\label{def:BSFnd}
    For $I \in \I_{n,d}$, denote by $\Dt_{I}$ the geometric object representative of the equivalence class of $\Dhat_I$ in the quotient space $\Sym^d(\C)$, which satisfies $\Dt_I=\bpi(\Dhat_I)$. We define $\SOD$ to be the collection of such objects, equipped with arbitrary brane structures, and we denote its endomorphism algebra over $\kk$ by $\BFS_{n,d}$.
\end{definition}

In order to check that $\SOD$ in Definition~\ref{def:BSFnd} forms a full exceptional collection for $\F(\f_{n,d})$, we describe the vanishing paths associated to the objects underlying $\SOD$, and check that these form a (non-generic) distinguished collection. Fix the regular value for $\w_{n,d}$ to be $*=-id$.

\begin{definition}
    We denote by $\lambdabold$ the collection of straight segments between the (real, non-negative) critical values of $\what$ and $-id$.
\end{definition}

\begin{proposition}\label{prop:exccollectionf}
    The collection $\lambdabold$ is one of non-generic vanishing paths for the Lefschetz fibration $\w_{n,d}$, and the corresponding Lefschetz thimbles (equipped with arbitrary brane structures) form a full exceptional collection in $\F(\f_{n,d})$. In particular, $\SOD$ is a full exceptional collection in $\F(\f_{n,d})$, and $\F(\f_{n,d})$ admits a semi-orthogonal decomposition $\F(\f_{n,d})=\langle \SODi_0,\dots,\SODi_d \rangle$, where $\SODi_h$ is the strictly full triangulated subcategory generated by $\SODhati^{\circ}_h /\S_d$.
\end{proposition}

\begin{proof}
    The fact that $\lambdabold$ is a non-generic collection follows from Lemma~\ref{lem:Milnornumber}. The vanishing paths associated to $\SOD$ are strictly contained in the collection $\lambdaboldhat$, since the latter is fixed by the action of the symmetric group on $\SODhat$. In particular, the vanishing paths associated to $\SOD$ coincide with the collection $\lambdabold$.
\end{proof}

\begin{lemma}\label{lem:injectvspaces}
    Let $I,J \in \I_{n,d}$, and $\sigma \in \S_{d}$. If the morphism space
    \begin{equation*}
        \hom_{\F(\fhat_{n,d})}(\Dhat_I,\sigma \Dhat_J)
    \end{equation*}
    is non-trivial, $\sigma$ is the identity element in the symmetric group.
\end{lemma}

\begin{proof}
Fix $I,J \in \I_{n,d}$, and suppose $\hom_{\F(\fhat_{n,d})}(\Dhat_I,\sigma \Dhat_J)$ is non-trivial, for some non-trivial element $\sigma$ of the symmetric group. Let $\sigma J=(j'_1,\dots,j'_d)$. Since $\sigma J \neq J$, $\sigma J$ cannot be in $\I_{n,d}$, so there is an index $h\in\{1,\dots,d-1\}$, such that $j'_{h}>j'_{h+1}$. By Proposition~\ref{prop:Gabungraded}, we have a combinatorial condition for the non-vanishing of the morphism space in the statement. Namely, the following two conditions must hold:
    \begin{enumerate}[i.]
        \item $j'_{h}=i_{h}$, or $j'_{h}=i_{h}-1$, or $j'_{h}=i_{h}+1$;
        \item $j'_{h+1}=i_{h+1}$, or $j'_{h+1}=i_{h+1}-1$, or $j'_{h+1}=i_{h+1}+1$.
    \end{enumerate}
If, for example, $j'_{h}=i_{h}$ and $j'_{h+1}=i_{h+1}$, we cannot have $j'_{h}>j'_{h+1}$ since $i_{h}<i_{h+1}$ by assumption. The only non-trivial observation should be made in the following case. Suppose $j'_{h}=i_{h}+1$ and $j'_{h+1}=i_{h+1}-1$, and $i_{h}<i_{h+1}$ by initial assumption. In particular, again by Proposition~\ref{prop:Gabungraded}, both $i_{h}$ and $i_{h+1}$ are odd, so in fact one has $i_{h}+2\leq i_{h+1}$. It follows that in this case one also has $j'_{h}\leq j'_{h+1}$, which is a contradiction. Every other case is similarly verified.
\end{proof}

The goal for the remainder of the section is to compute the endomorphism algebra $\BFS_{n,d}$ in $\F(\f_{n,d})$. In order to do so, we use the computation for $\BFShat_{n,d}$ carried out in Section~\ref{sec:FukBP}, and we study how the holomorphic curves contributing to the $A_{\infty}$-products in $\BFShat_{n,d}$ interact with the diagonal locus $\Delta\subset \C^d$. 

We restrict the Morsification $\what_{n,1}$ to the preimage of the standard disc, and we fix $-i$ to be the regular value on the imaginary line.

\begin{remark}
    If $(x^{\textnormal{crit}}_1,\dots,x^{\textnormal{crit}}_d)$ is a critical point of $\what_{n,d}$, $\hat{w}$ its critical value, $\ghat_{\hat{w}}$ the straight line from $\hat{w}$ to $-id$, then $\ghat_I:=\ghat_{\hat{w}}=\ghat_{\hat{w}_1}+\dots+\ghat_{\hat{w}_d}$.
\end{remark}

For the remainder of the section, we fix $I\in\I_{n,d}$, i.e.\ we only consider the Lefschetz thimbles $\Dhat_I$ associated to the critical values $\hat{w}=\hat{w}_{i_1}+\dots+\hat{w}_{i_d}$, where for $i_h<i_k$ and $n$ even (resp.\ odd) we have the following:
\begin{equation*}
    \begin{cases*}
        \hat{w}_{i_h}<\hat{w}_{i_k} \quad &\text{if $i_h$ and $i_k$ are odd} \\
        \hat{w}_{i_h}>\hat{w}_{i_k} \quad &\text{if $i_h$ and $i_k$ are even} \\
        \hat{w}_{i_h} \leq \frac{1}{2}+\frac{n-2}{2}\epsilon-\hat{w}_{i_k} \text{ (resp.\ $\frac{1}{2}+\frac{n-3}{2}\epsilon-\hat{w}_{i_k}$)} \quad &\text{if $i_h$ is odd, $i_k$ even} \\
        \hat{w}_{i_h} > \frac{1}{2}+\frac{n-2}{2}\epsilon-\hat{w}_{i_k} \text{ (resp.\ $\frac{1}{2}+\frac{n-3}{2}\epsilon-\hat{w}_{i_k}$)} \quad &\text{if $i_h$ is even and $i_k$ odd}.
    \end{cases*}
\end{equation*}

The following proposition realises the thimbles in $\what_{n,d}^{-1}(\D_d)$ (where $\D_d$ is the disc of radius $d$) described above as products of thimbles in $\what_{n,1}^{-1}(\D)$, and its proof is similar to that of \cite[Lemma~2.2]{Auroux2}.

\begin{proposition}\label{prop:productsofthimbles}
    For each $I\in \I_{n,d}$, $\Dhat_I$ is isotopic to $\prod_{h=1}^d \Dhat_{i_h}$ as exact Lagrangians.
\end{proposition}

\begin{proof}
    Fix $I=(i_1,\dots,i_d)\in \I_{n,d}$, and let $\hat{w}=\hat{w}_{i_1}+\dots+\hat{w}_{i_d}$ be the critical value associated to $\Dhat_I$. Denote by $\nabla \Real \what$ and $\nabla \Imaginary\what$ the gradient vector fields for the real and imaginary part of $\what:=\what_{n,d}$ with respect to the K\"ahler metric $g=\omega(\cdot,J \cdot)$, where $J$ is the product complex structure on $\C^d$. Since $\what$ is holomorphic, the horizontal distribution on $\C^d$ giving symplectic parallel transport is spanned by $\nabla \Real \what$ and $\nabla \Imaginary\what$. Moreover, since $\ghat_{\hat{w}}$ is a straight segment from the real value $\hat{w}$ to $-id$, parallel transport along it is given by the gradient flow of
    \begin{equation*}
    -(\hat{w} \nabla \Real \what + d\nabla \Imaginary\what)
    \end{equation*}
     with respect to $g$. As $g$ is a product metric, the gradient vector at a point $(x_1,\dots,x_d)$ decomposes into
     \begin{equation*}
         -(\hat{w}_{i_1} \nabla \Real \what_{n,1}(x_1) + \nabla \Imaginary\what_{n,1}(x_1)),\quad \ldots \quad,-(\hat{w}_{i_d} \nabla \Real \what_{n,1}(x_d) + \nabla\Imaginary \what_{n,1}(x_d)).
     \end{equation*}
    Analogously, if $x_i \in \what_{n,1}^{-1}(\ghat_{i_h})$ for some $h$, the parallel transport along $\ghat_{i_h}$ (which is a straight segment from $\hat{w}_{i_h}$ to $-i$) is given by
    \begin{equation*}
        -(\hat{w}_{i_h} \nabla  \Real \what_{n,1}(x_i) + \nabla \Imaginary\what_{n,1}(x_i)).
    \end{equation*}
    Therefore, parallel transport along $\ghat_{\hat{w}}$ decomposes into the product of the parallel transports along each $\ghat_{\hat{w}_{i_h}}$, and $\Dhat_I$ is exact Lagrangian isotopic to the product of Lagrangian thimbles in $\what^{-1}(\D_d)$.
\end{proof}

In analogy with \cite{Auroux2}, we will need to consider perturbed versions of the thimbles in Proposition~\ref{prop:productsofthimbles}. Fix a small number $\rho$ with positive real part, such that $-d(i\pm\rho)$ is on the boundary of the disc of radius $d$ in $\C$,  and consider the straight segments $\ghat_{I}^{\pm}$, connecting the critical values of $\what_{n,d}$ to $-d(i\pm\rho)$. This is depicted in Figure~\ref{fig:perturbedproduct} (right), see also \cite[Figures~1 and 3]{Auroux2}.

\begin{figure}[t]
    \centering
    \begin{minipage}{.45 \textwidth}
        \centering
    \begin{tikzpicture}[scale=.7,v/.style={draw,shape=circle, fill=black, minimum size=1.2mm, inner sep=0pt, outer sep=0pt}]
        \draw (1,0) ellipse (4 and 2);
    
        \node at (-.1,.5){$\ghat_1$};
        \draw[] (.1,.1) -- (-.1,-.1);
        \draw[] (-.1,.1) -- (.1,-.1);

        \node at (.7,0){$\dots$};

        \node at (.1+.6+.8,.5){$\ghat_{n-1}$};
        \draw[] (.1+.6+.8,.1) -- (.8-.1+.6,-.1);
        \draw[] (.8-.1+.6,.1) -- (.8+.1+.6,-.1);

        \node at (3,.5){$\ghat_{n}$};
        \draw[] (3+.1,.1) -- (3+-.1,-.1);
        \draw[] (3+-.1,.1) -- (3+.1,-.1);

        \node at (3.65,0){$\dots$};

        \node at (1+3+.3,.5){$\ghat_{2}$};
        \draw[] (1+3+.1+.3,.1) -- (1+3+-.1+.3,-.1);
        \draw[] (1+3+-.1+.3,.1) -- (1+3+.1+.3,-.1);

        \node[] (star) at (0,-2+.05) {$*$};
        \node[] at (0,-2.5) {$-i$};

        \node[color=red] (star) at (0-.5,-2+.15) {$*$};
        \node[color=red] at (-1.25,-2.5) {$-i-\rho$};

        \node[color=blue] (star) at (0+.5,-2) {$*$};
        \node[color=blue] at (1.25,-2.5) {$-i+\rho$};
        \node[color=red] at (-.5,-.5){$\ghat_1^{+}$};
        \node[color=blue] at (.5,-.5){$\ghat_1^{-}$};
        \draw[] (0,0) -- (0,-2+.05);
        \draw[] (1.4,0) -- (0,-2+.05);

        \draw[] (3,0) -- (0,-2+.05);
        \draw[] (4.3,0) -- (0,-2+.05);

        \draw[color=red] (0,0) -- (-.5,-2+.15);
        \draw[color=red] (1.4,0) -- (-.5,-2+.15);

        \draw[color=red] (3,0) -- (-.5,-2+.15);
        \draw[color=red] (4.3,0) -- (-.5,-2+.15);

        \draw[color=blue] (0,0) -- (.5,-2);
        \draw[color=blue] (1.4,0) -- (.5,-2);

        \draw[color=blue] (3,0) -- (.5,-2);
        \draw[color=blue] (4.3,0) -- (.5,-2);
    \end{tikzpicture}
\end{minipage}
    \begin{minipage}{.45 \textwidth}
        \centering
        \begin{tikzpicture}[scale=.7,v/.style={draw,shape=circle, fill=black, minimum size=1.2mm, inner sep=0pt, outer sep=0pt}]
            \draw (1,0) ellipse (4 and 2);
                
                    \node at (-.1,.5){$\ghat_I$};
                    \draw[] (.1,.1) -- (-.1,-.1);
                    \draw[] (-.1,.1) -- (.1,-.1);
            
                    \node at (.1+.6+.8,.5){$\ghat_J$};
                    \draw[] (.1+.6+.8,.1) -- (.8-.1+.6,-.1);
                    \draw[] (.8-.1+.6,.1) -- (.8+.1+.6,-.1);
            
                    \node at (3,.5){$\ghat_K$};
                    \draw[] (3+.1,.1) -- (3+-.1,-.1);
                    \draw[] (3+-.1,.1) -- (3+.1,-.1);
            
                    \node[] (star) at (0,-2+.05) {$*$};
                    \node[] at (0,-2.5) {$-id$};
            
                    \node[color=red] (star) at (0-.5,-2+.15) {$*$};
                    \node[color=red] at (-1.75,-2.5) {$-d(i+\rho)$};
            
                    \node[color=blue] (star) at (0+.5,-2) {$*$};
                    \node[color=blue] at (1.75,-2.5) {$-d(i-\rho)$};
                    \node[color=red] at (-.5,-.5){$\ghat_I^{+}$};
                    \node[color=blue] at (.5,-.5){$\ghat_I^{-}$};
                    \draw[] (0,0) -- (0,-2+.05);
                    \draw[] (1.4,0) -- (0,-2+.05);
            
                    \draw[] (3,0) -- (0,-2+.05);
            
                    \draw[color=red] (0,0) -- (-.5,-2+.15);
                    \draw[color=red] (1.4,0) -- (-.5,-2+.15);
            
                    \draw[color=red] (3,0) -- (-.5,-2+.15);

                    \draw[color=blue] (0,0) -- (.5,-2);
                    \draw[color=blue] (1.4,0) -- (.5,-2);
            
                    \draw[color=blue] (3,0) -- (.5,-2);
               
                \end{tikzpicture}
    \end{minipage}
\caption{The triple diagrams $(\D, \ghat_{\hat{w}_j}^{-}, \ghat_{\hat{w}_j}, \ghat_{\hat{w}_j}^{+})$ (left) and $(\D_d, \ghat_I^{-}, \ghat_I, \ghat_I^{+})$ (right).}
\label{fig:perturbedproduct}
\end{figure}

\begin{definition}\label{def:pushoffthimbles}
    For $I\in \I_{n,d}$, define $\Dhat_I^{\pm}$ to be the thimble associated to the straight segment $\ghat_{I}^{\pm}$, viewed as a vanishing path for the Fukaya--Seidel category of $\what_{n,d}$ with fixed regular value $*=-d(i\pm\rho)$.
\end{definition}

The following is analogous to Proposition~\ref{prop:productsofthimbles} for the perturbed versions of the thimbles.

\begin{lemma}\label{lem:productsofperturbedthimbles}
    For $I\in\I_{n,d}$, each $\Dhat_I^{\pm}$ is isotopic to $\prod_{h=1}^d \Dhat_{i_h}^{\pm}$ as exact Lagrangians.
\end{lemma}

\begin{proof}
    This is analogous to Proposition~\ref{prop:productsofthimbles}, noting that each $\ghat_{\hat{w}_{i_h}}^{\pm}$ is the straight segment from $\hat{w}_{i_h}$ to $-i\mp \rho$. This is shown in Figure~\ref{fig:perturbedproduct} (left), where in blue, black and red we depicted the vanishing paths $\ghat_{\hat{w}_j}^{-}, \ghat_{\hat{w}_j}, \ghat_{\hat{w}_j}^{+}$.
\end{proof}

The thimbles described in Definition~\ref{def:pushoffthimbles} define objects of the Fukaya--Seidel category of $\what_{n,d}$, where we have chosen different base points. Nonetheless, we can realise each holomorphic triangle contributing to a composition
\begin{equation*}
    \hom_{\F(\what_{n,d})}(\Dhat_J,\Dhat_K) \otimes \hom_{\F(\what_{n,d})}(\Dhat_I,\Dhat_J) \to \hom_{\F(\what_{n,d})}(\Dhat_I,\Dhat_K)
\end{equation*}
as a holomorphic triangle in $\what_{n,d}^{-1}(\D_d)$, with boundary conditions on $\Dhat_I^{-}\cup\Dhat_J\cup\Dhat_K^{+}$. More precisely, we have the following.

\begin{lemma}\label{lem:invarianceMaslov}
    Fix $I,J,K\in \I_{n,d}$. The Maslov index of the holomorphic triangle
    \begin{equation*}
        (u,\partial u): D \to (\what_{n,d}^{-1}(\D_d),\Dhat_I\cup\Dhat_J\cup\Dhat_K)
    \end{equation*}
    (when this exists) is equal to the Maslov index of the holomorphic triangle
    \begin{equation*}
        (u,\partial u): D \to (\what_{n,d}^{-1}(\D_d),\Dhat_I^{-}\cup\Dhat_J\cup\Dhat_K^{+}).
    \end{equation*}
\end{lemma}

\begin{proof}
    Each $\Dhat_I^{-}$ is Hamiltonian isotopic to $\Dhat_I$, i.e.\ it is the image under the time one flow of the vector field of a non-negative Hamiltonian function $H^{-}: \C^d \to \CR$, which is linear on the positive symplectisation of $\what_{n,d}^{-1}(\D_d)$. Similarly, each $\Dhat_I$ is Hamiltonian isotopic to $\Dhat_I^{+}$, for a non-negative Hamiltonian function $H^{+}: \C^d \to \CR$ satisfying the same conditions. The claim follows from the fact that the Maslov index is invariant under Hamiltonian isotopy.
\end{proof}

Proposition~\ref{prop:productsofthimbles}, Lemma~\ref{lem:productsofperturbedthimbles} and Lemma~\ref{lem:invarianceMaslov} allow us to compute morphisms and compositions by looking at intersection points and two-chains in Figure~\ref{fig:perturbedproduct} (left), using a combinatorial description of holomorphic triangle count similar to that defined for Heegaard diagrams \cite{LOT18}. Note that the diagram in Figure~\ref{fig:perturbedproduct} can have triple intersections; when this is the case, we can perturb them (in particular, we can perturb the triple intersection at the real branch values).

For the remainder of this section, fix a triple $(I,J,K)$, so that there is a Maslov index 0 triangle, called $\phi$, with boundary conditions on $\Dhat_I^{-}\cup\Dhat_J\cup\Dhat_K^{+}$, contributing to the product
\begin{equation*}
    \hom_{\F(\what_{n,d})}(\Dhat_J,\Dhat_K) \otimes \hom_{\F(\what_{n,d})}(\Dhat_I,\Dhat_J) \to \hom_{\F(\what_{n,d})}(\Dhat_I,\Dhat_K).
\end{equation*}
Sufficient and necessary conditions for this to exist are given in Proposition~\ref{prop:Gabungraded}, which we summarise in Definition~\ref{def:Maslovcondition}. We will prove in Proposition~\ref{prop:disjoindiagonal} that $\phi$ is disjoint from the diagonal locus $\Delta \subset \C^d$. We first realise $\phi$ as a $2$-chain in $\C$.

\begin{definition}\label{def:Maslovcondition}
    For $I,J,K \in \I_{n,d}$, we say the triple $(\Dhat_I^{-},\Dhat_J,\Dhat_K^{+})$ is \emph{composable} if for each $h=1,\dots,d$, we have that $i_h=j_h=k_h$, or that $i_h=j_h=k_h\pm 1$ and $i_h$ is odd, or that $k_h=j_h=i_h\pm 1$ and $i_h$ is odd. When this holds, we pick generators $x\in \Dhat_I^{-}\cap \Dhat_J$, $y\in \Dhat_J\cap \Dhat_K^{+}$, and $z\in \Dhat_I^{-}\cap \Dhat_K^{+}$, which are realised as the $d$-tuples $x=(x_1,\dots,x_d)$, $y=(y_1,\dots,y_d)$,  $z=(z_1,\dots,z_d)$, for $x_{h}\in \Dhat_{i_h}^{-}\cap \Dhat_{j_h}$, $y_{h}\in \Dhat_{j_h}\cap \Dhat_{k_h}^{+}$, $z_{h}\in \Dhat_{i_h}^{-}\cap \Dhat_{k_h}^{+}$.
\end{definition}

\begin{definition}\label{def:phi'}
    Fix $I,J,K \in \I_{n,d}$ so that $(\Dhat_I^{-},\Dhat_J,\Dhat_K^{+})$ is composable. Denote by $\phi'$ the 2-chain in $\D\subset \C$ given by the union over $h$ of $\phi'_h$ in $\D\subset \C$, where each $\phi'_h$ is the holomorphic triangle with boundary conditions on $\Dhat_{i_h}^{-}\cup\Dhat_{j_h}\cup\Dhat_{k_h}^{+}$. More specifically, each $\phi'_h$ denotes the image of a holomorphic map ${u_h: (D,\partial D) \to (\what_{n,1}^{-1}(\D), \Dhat_{i_h}^{-}\cup\Dhat_{j_h}\cup\Dhat_{k_h}^{+})}$.
\end{definition}
We can depict $\phi'$ in Definition~\ref{def:phi'} by its projection under $\what_{n,d}$, which realises it as a 2-chain in the base of $\what_{n,1}$, with boundary on the arcs of the diagram. See Figure~\ref{fig:tupleoftrigs} (left), where for $n=7$ and $d=3$ we depicted the 2-chain $\phi'$, with boundary conditions on $\Dhat_{1,5,7}^{-}\cup\Dhat_{2,4,7}\cup\Dhat_{2,4,6}^{+}$. The three components of $\phi'$ are $\phi'_{1},\phi'_{2}$ and $\phi'_{3}$, here depicted as $\what_{n,1}(\phi'_{1})$, $\what_{n,1}(\phi'_{2})$, and $\what_{n,1}(\phi'_{3})$. These have, respectively, boundary conditions on $\ghat_{1}^{-}\cup\ghat_{2}\cup\ghat_{2}^{+}$, $\ghat_{5}^{-}\cup\ghat_{4}\cup\ghat_{4}^{+}$, and $\ghat_{5}^{-}\cup\ghat_{4}\cup\ghat_{4}^{+}$, which are the shaded triangles with vertex in $\hat{w}_2$, $\hat{w}_4$, and $\hat{w}_7$ respectively. The following definition and proposition give a necessary condition which determines which 2-chains can occur.

\begin{definition}[{\cite[Section~2.3]{LMW08}}]
    We say that a pair of embedded triangles with boundary on $\ghat_{i_h}^{-}\cup\ghat_{j_h}\cup\ghat_{k_h}^{+}$ overlap ``\emph{head-to-tail}'' if they overlap in the configuration depicted in Figure~\ref{fig:tupleoftrigs} (right), where parallel sides are ``of the same type'', i.e.\ any of $\ghat,\ghat^{-},\ghat^{+}$.
\end{definition}

\begin{figure}[t]
    \centering
    \begin{minipage}{.6 \textwidth}
        \centering
    \begin{tikzpicture}[scale=.8,v/.style={draw,shape=circle, fill=black, minimum size=1.2mm, inner sep=0pt, outer sep=0pt}]
        \draw (1,0) ellipse (4.5 and 2);

        \node at (.1-1+.8-.1,.5){$\hat{w}_{1}$};
        \draw[] (.1-1+.8,.1) -- (.8-.1-1,-.1);
        \draw[] (.8-.1-1,.1) -- (.8+.1-1,-.1);

        \node at (.1-.5+.8,.5){$\hat{w}_{3}$};
        \draw[] (.1-.5+.8,.1) -- (.8-.1-.5,-.1);
        \draw[] (.8-.1-.5,.1) -- (.8+.1-.5,-.1);

        \node at (.1+.8,.5){$\hat{w}_{5}$};
        \draw[] (.1+.8,.1) -- (.8-.1,-.1);
        \draw[] (.8-.1,.1) -- (.8+.1,-.1);

        \node at (.1+.5+.8+.1,.5){$\hat{w}_{7}$};
        \draw[] (.1+.5+.8,.1) -- (.8-.1+.5,-.1);
        \draw[] (.8-.1+.5,.1) -- (.8+.1+.5,-.1);

        \node at (3-.5-.1,.5){$\hat{w}_{6}$};
        \draw[] (3+.1-.5,.1) -- (3+-.1-.5,-.1);
        \draw[] (3+-.1-.5,.1) -- (3+.1-.5,-.1);
        
        \node at (3,.5){$\hat{w}_{4}$};
        \draw[] (3+.1,.1) -- (3+-.1,-.1);
        \draw[] (3+-.1,.1) -- (3+.1,-.1);

        \node at (3+.5+.1,.5){$\hat{w}_{2}$};
        \draw[] (3+.1+.5,.1) -- (3+-.1+.5,-.1);
        \draw[] (3+-.1+.5,.1) -- (3+.1+.5,-.1);
 
        \node[] (star) at (1,-2) {$*$};
        \node[color=red] (star) at (1-.5,-2) {$*$};
        \node[color=blue] (star) at (1+.5,-2) {$*$};
        \draw[color=blue] (.8-1,0) -- (1+.5,-2);
        \draw[color=blue] (.8,0) -- (1+.5,-2);
        \draw[color=blue] (.8+.5,0) -- (1+.5,-2);

        \draw[] (.8+.5,0) -- (1,-2);
        \draw[] (3,0) -- (1,-2);        
        \draw[] (3+.5,0) -- (1,-2);

        \draw[color=red] (3-.5,0) -- (.5,-2);
        \draw[color=red] (3+.5,0) -- (.5,-2);
        \draw[color=red] (3,0) -- (.5,-2);

        \draw[fill=gray,draw=none, fill opacity=.3] (1.3,0) to  (1.4,-1.1) to  (1.1,-1.4);
        \draw[fill=gray,draw=none, fill opacity=.3] (3,0) to (1.37,-1.63) to  (1.28,-1.37);
        \draw[fill=gray,draw=none, fill opacity=.3] (3.5,0) to (1.3,-1.75) to  (1.15,-1.57);
    \end{tikzpicture}
\end{minipage}
    \begin{minipage}{.3 \textwidth}
        \centering
        \begin{tikzpicture}[scale=.95,v/.style={draw,shape=circle, fill=black, minimum size=1.2mm, inner sep=0pt, outer sep=0pt}]
            \draw (-1,1) to (-1,-1);
            \draw (-1,-1) to (1,0);
            \draw (-1,1) to (1,0);
        
            \draw (-1+1,1) to (-1+1,-1);
            \draw (-1+1,-1) to (1+1,0);
            \draw (-1+1,1) to (1+1,0);
        
            \end{tikzpicture}
    \end{minipage}
    \caption{(left) The 2-chain $\phi'$ as product of triangles. (right) ``\emph{Head-to-tail}'' overlap.}
    \label{fig:tupleoftrigs}
\end{figure}

\begin{proposition}\label{prop:headtail}
    Fix $I,J,K \in \I_{n,d}$ so that $(\Dhat_I^{-},\Dhat_J,\Dhat_K^{+})$ is composable, and let $\phi'$ be as in Definition~\ref{def:phi'}. The image of $\phi'$ under $\what_{n,d}$ is the union of embedded triangles which are either disjoint or overlap pairwise head-to-tail. 
 \end{proposition}

 \begin{remark}
    We remark the similarities with  \cite[Lemma~2.6]{LMW08}, which states that when we have a \emph{nice} diagram in the sense of Sarkar \cite{Sarkar}, if a rigid holomorphic triangle is represented by the union of embedded holomorphic triangles, these pairwise are either disjoint or overlap head-to-tail. In our case, the diagram in not nice. This is similar to the proof of \cite[Proposition~3.5]{Auroux2}, in which setting the Floer product counts unions of embedded triangles which either are disjoint or overlap head-to-tail, despite the diagram not being nice.
\end{remark}

\begin{proof}
    By Definition~\ref{def:Maslovcondition}, there are 3 possible configurations for each triangle $\phi'_h$:
    \begin{enumerate}
        \item If $\hat{w}_{i_h}=\hat{w}_{j_h} = \hat{w}_{k_h}$, the image under $\what_{n,1}$ of $\phi'_h$ is a ``degenerate'' triangle, with vertices at the branch value $\hat{w}_{i_h}$ of $\what_{n,1}$;
        \item If $\hat{w}_{i_h}=\hat{w}_{j_h} \neq \hat{w}_{k_h}$, the image under $\what_{n,1}$ of $\phi'_h$ has one vertex at the branch value $\hat{w}_{i_h}$, as depicted in Figure~\ref{fig:possiblephiprime} (left);
        \item If $\hat{w}_{i_h}\neq\hat{w}_{j_h} = \hat{w}_{k_h}$, the image under $\what_{n,1}$ of $\phi'_h$ has one vertex at the branch value $\hat{w}_{j_h}$, and is depicted in Figure~\ref{fig:possiblephiprime} (right).
    \end{enumerate}

    \begin{figure}[t]
        \centering
        \begin{minipage}{.45 \textwidth}
            \centering
        \begin{tikzpicture}[v/.style={draw,shape=circle, fill=black, minimum size=1.2mm, inner sep=0pt, outer sep=0pt}]
    
            \node at (.1+.6+.8,.5){\color{blue}$\ghat_{i_h}^{-}$\color{black}, $\ghat_{j_h}$};
            \draw[] (.1+.6+.8,.1) -- (.8-.1+.6,-.1);
            \draw[] (.8-.1+.6,.1) -- (.8+.1+.6,-.1);
    
            \node at (3,.5){\color{red}$\ghat_{k_h}^{+}$\color{black}};
            \draw[] (3+.1,.1) -- (3+-.1,-.1);
            \draw[] (3+-.1,.1) -- (3+.1,-.1);
     
            \node[] (star) at (0,-2) {$*$};
            \node[color=red] (star) at (0-.5,-2) {$*$};
            \node[color=blue] (star) at (0+.5,-2) {$*$};
          
            \draw[] (1.4,0) -- (0,-2);        
            \draw[color=red] (3,0) -- (-.5,-2);
            \draw[color=blue] (1.4,0) -- (.5,-2);
       
            \draw[fill=gray,draw=none, fill opacity=.3] (1.4,0) to (.35,-1.5) to  (.85,-1.25);
        \end{tikzpicture}
    \end{minipage}
    \begin{minipage}{.45 \textwidth}
        \centering
        \begin{tikzpicture}[v/.style={draw,shape=circle, fill=black, minimum size=1.2mm, inner sep=0pt, outer sep=0pt}]
    
            \node at (.1+.6+.8,.5){\color{blue}$\ghat^{-}_{i_h}$\color{black}};
            \draw[] (.1+.6+.8,.1) -- (.8-.1+.6,-.1);
            \draw[] (.8-.1+.6,.1) -- (.8+.1+.6,-.1);
    
            \node at (3,.5){$\ghat_{j_h}$, \color{red}$\ghat_{k_h}^{+}$\color{black}};
            \draw[] (3+.1,.1) -- (3+-.1,-.1);
            \draw[] (3+-.1,.1) -- (3+.1,-.1);
     
            \node[] (star) at (0,-2) {$*$};
            \node[color=red] (star) at (0-.5,-2) {$*$};
            \node[color=blue] (star) at (0+.5,-2) {$*$};
          
            \draw[] (3,0) -- (0,-2);        
            \draw[color=red] (3,0) -- (-.5,-2);
            \draw[color=blue] (1.4,0) -- (.5,-2);
       
            \draw[fill=gray,draw=none, fill opacity=.3] (3,0) to (.85,-1.23) to  (.72,-1.52);
        \end{tikzpicture}        
    \end{minipage}
    \caption{Possible configurations for $\phi'_h$.}
    \label{fig:possiblephiprime}
    \end{figure}

Now consider two embedded triangles $\phi'_h$, $\phi'_{\ell}$. If they are both of type $1$, they are clearly disjoint (as, in this case, $\hat{w}_{i_h}=\hat{w}_{j_h} = \hat{w}_{k_h} \neq \hat{w}_{i_{\ell}}=\hat{w}_{j_{\ell}} = \hat{w}_{k_{\ell}}$). If they are both of type $2$ (or, analogously, if they are both of type $3$), their images under $\what_{n,1}$ intersect head-to-tail, as depicted in the two left-most pictures in Figure~\ref{fig:headtailalltypes} ($\phi'_h$ and $\phi'_{\ell}$ are shaded in orange and blue respectively). Moreover, if $\phi'_h$ is of type $1$, and $\phi'_{\ell}$ is of any other type, the triangles are disjoint (this is because $\phi'_{\ell}$ does not have a vertex on the branch value $\hat{w}_{i_h}=\hat{w}_{j_h}=\hat{w}_{k_h}$). Finally, if $\phi'_h$ is of type $2$ and $\phi'_{\ell}$ is of type $3$, their images under $\what_{n,1}$ are disjoint, as depicted in the two right-most pictures in Figure~\ref{fig:headtailalltypes} ($\phi'_h$ and $\phi'_{\ell}$ are shaded in orange and blue respectively).\qedhere

\begin{figure}[t]
    \centering
    \begin{minipage}{.24 \textwidth}
        \centering
    \begin{tikzpicture}[scale=.7,v/.style={draw,shape=circle, fill=black, minimum size=1.2mm, inner sep=0pt, outer sep=0pt}]

        \node at (.1+.6,.5){$\ghat_{i_h}$};
        \draw[] (.1+.8,.1) -- (.8-.1,-.1);
        \draw[] (.8-.1,.1) -- (.8+.1,-.1);

        \node at (.1+.6+.8,.5){$\ghat_{i_{\ell}}$};
        \draw[] (.1+.6+.8,.1) -- (.8-.1+.6,-.1);
        \draw[] (.8-.1+.6,.1) -- (.8+.1+.6,-.1);

        \node at (3-.1,.5){$\ghat_{k_{\ell}}$};
        \draw[] (3+.1,.1) -- (3+-.1,-.1);
        \draw[] (3+-.1,.1) -- (3+.1,-.1);

        \node at (3+.75,.5){$\ghat_{k_{h}}$};
        \draw[] (3+.1+.5,.1) -- (3+-.1+.5,-.1);
        \draw[] (3+-.1+.5,.1) -- (3+.1+.5,-.1);
 
        \node[] (star) at (0,-2) {$*$};
        \node[color=red] (star) at (0-.5,-2) {$*$};
        \node[color=blue] (star) at (0+.5,-2) {$*$};
        \draw[] (.8,0) -- (0,-2);     
        \draw[color=blue] (.8,0) -- (.5,-2);
        \draw[color=red] (3+.5,0) -- (-.5,-2);

        \draw[] (1.4,0) -- (0,-2);        
        \draw[color=red] (3,0) -- (-.5,-2);
        \draw[color=blue] (1.4,0) -- (.5,-2);
   
        \draw[fill=ballblue,draw=none, fill opacity=.3] (1.4,0) to  (.32,-1.53) to  (.83,-1.25);
        \draw[fill=orange,draw=none, fill opacity=.3] (.8,0) to (.13,-1.7) to  (.58,-1.45);
    \end{tikzpicture}
\end{minipage}
\begin{minipage}{.24 \textwidth}
    \centering
    \begin{tikzpicture}[scale=.7,v/.style={draw,shape=circle, fill=black, minimum size=1.2mm, inner sep=0pt, outer sep=0pt}]

        \node at (.1-1+.8,.5){$\ghat_{i_h}$};
        \draw[] (.1-1+.8,.1) -- (.8-.1-1,-.1);
        \draw[] (.8-.1-1,.1) -- (.8+.1-1,-.1);

        \node at (.1+.8,.5){$\ghat_{i_{\ell}}$};
        \draw[] (.1+.8,.1) -- (.8-.1,-.1);
        \draw[] (.8-.1,.1) -- (.8+.1,-.1);

        \node at (3-.1,.5){$\ghat_{k_{\ell}}$};
        \draw[] (3+.1,.1) -- (3+-.1,-.1);
        \draw[] (3+-.1,.1) -- (3+.1,-.1);

        \node at (3+.75,.5){$\ghat_{k_{h}}$};
        \draw[] (3+.1+.5,.1) -- (3+-.1+.5,-.1);
        \draw[] (3+-.1+.5,.1) -- (3+.1+.5,-.1);
 
        \node[] (star) at (1,-2) {$*$};
        \node[color=red] (star) at (1-.5,-2) {$*$};
        \node[color=blue] (star) at (1+.5,-2) {$*$};
        \draw[color=blue] (.8-1,0) -- (1+.5,-2);
        \draw[color=blue] (.8,0) -- (1+.5,-2);

        \draw[] (3,0) -- (1,-2);        
        \draw[] (3+.5,0) -- (1,-2);        
        \draw[color=red] (3+.5,0) -- (.5,-2);
        \draw[color=red] (3,0) -- (.5,-2);

        \draw[fill=ballblue,draw=none, fill opacity=.3] (3,0) to (1.37,-1.63) to  (1.28,-1.37);
        \draw[fill=orange,draw=none, fill opacity=.3] (3.5,0) to (1.3,-1.75) to  (1.15,-1.57);
    \end{tikzpicture}
\end{minipage}
\begin{minipage}{.24\textwidth}
    \centering
\begin{tikzpicture}[scale=.7,v/.style={draw,shape=circle, fill=black, minimum size=1.2mm, inner sep=0pt, outer sep=0pt}]

    \node at (.1+.6,.5){$\ghat_{i_h}$};
    \draw[] (.1+.8,.1) -- (.8-.1,-.1);
    \draw[] (.8-.1,.1) -- (.8+.1,-.1);

    \node at (.1+.6+.8,.5){$\ghat_{i_{\ell}}$};
    \draw[] (.1+.6+.8,.1) -- (.8-.1+.6,-.1);
    \draw[] (.8-.1+.6,.1) -- (.8+.1+.6,-.1);

    \node at (3-.1,.5){$\ghat_{k_{\ell}}$};
    \draw[] (3+.1,.1) -- (3+-.1,-.1);
    \draw[] (3+-.1,.1) -- (3+.1,-.1);

    \node at (3+.75,.5){$\ghat_{k_{h}}$};
    \draw[] (3+.1+.5,.1) -- (3+-.1+.5,-.1);
    \draw[] (3+-.1+.5,.1) -- (3+.1+.5,-.1);

    \node[] (star) at (0,-2) {$*$};
    \node[color=red] (star) at (0-.5,-2) {$*$};
    \node[color=blue] (star) at (0+.5,-2) {$*$};
    \draw[] (.8,0) -- (0,-2);     
    \draw[color=blue] (.8,0) -- (.5,-2);
    \draw[color=red] (3+.5,0) -- (-.5,-2);

    \draw[color=blue] (1.4,0) -- (.5,-2);
    \draw[color=red] (3,0) -- (-.5,-2);

    \draw[] (3,0) -- (0,-2);
    \draw[fill=orange,draw=none, fill opacity=.3] (.8,0) to (.12,-1.7) to  (.58,-1.45);
    \draw[fill=ballblue,draw=none, fill opacity=.3] (3,0) to (.83,-1.25) to  (.7,-1.52);
\end{tikzpicture}
\end{minipage}
\begin{minipage}{.24\textwidth}
\centering
\begin{tikzpicture}[scale=.7,v/.style={draw,shape=circle, fill=black, minimum size=1.2mm, inner sep=0pt, outer sep=0pt}]

    \node at (.1+.6,.5){$\ghat_{i_{\ell}}$};
    \draw[] (.1+.8,.1) -- (.8-.1,-.1);
    \draw[] (.8-.1,.1) -- (.8+.1,-.1);

    \node at (.1+.6+.8,.5){$\ghat_{i_h}$};
    \draw[] (.1+.6+.8,.1) -- (.8-.1+.6,-.1);
    \draw[] (.8-.1+.6,.1) -- (.8+.1+.6,-.1);

    \node at (3-.1,.5){$\ghat_{k_h}$};
    \draw[] (3+.1,.1) -- (3+-.1,-.1);
    \draw[] (3+-.1,.1) -- (3+.1,-.1);

    \node at (3+.75,.5){$\ghat_{k_{\ell}}$};
    \draw[] (3+.1+.5,.1) -- (3+-.1+.5,-.1);
    \draw[] (3+-.1+.5,.1) -- (3+.1+.5,-.1);

    \node[] (star) at (0,-2) {$*$};
    \node[color=red] (star) at (0-.5,-2) {$*$};
    \node[color=blue] (star) at (0+.5,-2) {$*$};
    \draw[] (1.4,0) -- (0,-2);     
    \draw[color=blue] (.8,0) -- (.5,-2);
    \draw[color=red] (3+.5,0) -- (-.5,-2);

    \draw[color=blue] (1.4,0) -- (.5,-2);
    \draw[color=red] (3,0) -- (-.5,-2);

    \draw[] (3.5,0) -- (0,-2);
    \draw[fill=orange,draw=none, fill opacity=.3] (1.4,0) to (.3,-1.55) to  (.85,-1.22);
    \draw[fill=ballblue,draw=none, fill opacity=.3] (3.5,0) to (.57,-1.45) to  (.55,-1.68);
\end{tikzpicture}
\end{minipage}
\caption{All possible head-to-tail or empty overlaps.}
\label{fig:headtailalltypes}
\end{figure}

\end{proof}

The following two propositions identify $\phi'$ with $\phi$.

\begin{proposition}\label{prop:phiprimedisjoindiagonal}
    Fix $I,J,K \in \I_{n,d}$ so that $(\Dhat_I^{-},\Dhat_J,\Dhat_K^{+})$ is composable, and let $\phi'$ be as in Definition~\ref{def:phi'}. Then $\phi'$ is disjoint from the diagonal locus $\Delta$ in $\C^d$.
 \end{proposition}

 \begin{proof}
    The proof is similar to part of that of \cite[Proposition~3.5]{Auroux2}, which concerns the computation of the Maslov index of certain holomorphic discs in  $\Sym^d(\C)$, and their intersection number with the diagonal in $\Sym^d(\C)$. We first prove that the algebraic intersection number $\#(\phi'\cap \Delta)$ is zero. Note that each thimble $\Dhat_I^{-}, \Dhat_J, \Dhat_K^{+}$ is disjoint from $\Delta$, so the boundary of $\phi'$ is disjoint from $\Delta$, and $\#(\phi'\cap\Delta)$ is well-defined. Since $\phi'$ is union of triangles in $\what_{n,1}^{-1}(\D)$, this is equivalent to proving that the intersection number $\#(\phi_h'\cap \phi_{\ell}')$ vanishes for any pair $(\phi_h',\phi_{\ell}')$. By Proposition~\ref{prop:headtail}, there are two possibilities: either $\phi'_h$ and $\phi'_{\ell}$ are disjoint, or their images under $\what_{n,1}$ overlap head-to-tail. In the first case, $\#(\phi'_h\cap\phi'_{\ell})=0$. On the other hand, if $\phi'_h$ and $\phi'_{\ell}$ overlap head-to-tail, denote by $(u_{h},u_{\ell})$ the product map, such that the image of the unit disc under $u_{h}$ and $u_{\ell}$ is $\phi'_h$ and $\phi'_{\ell}$ respectively. The intersection number $\#(\phi'_h\cap\phi'_{\ell})$ can be evaluated by considering the rotation number of the boundaries $\partial\phi'_h,\partial \phi'_{\ell}$ around each other. More precisely, the restriction of the map $u_{\ell}-u_{h}$ to its boundary can be seen as a map defining a loop in $\C^\times$, and its degree is $\#(\phi'_h\cap\phi'_{\ell})$. Note that this restriction avoids the origin in $\C$ precisely because the head-to-tail configuration does not allow for sides ``of the same type'' to intersect. In the case of head-to-tail overlap, the degree of the map is 0. By positivity of intersections \cite[Appendix~B]{Wendl20}, 
    $\phi'_h$ and $\phi'_{\ell}$ are disjoint. It follows that $\phi'$ and $\Delta$ are disjoint.
 \end{proof}

 \begin{proposition}\label{prop:vanishingMaslov}
    Fix $I,J,K \in \I_{n,d}$ so that $(\Dhat_I^{-},\Dhat_J,\Dhat_K^{+})$ is composable, and let $\phi'$ be as in Definition~\ref{def:phi'}. Then $\phi'$ is a holomorphic triangle in $\what_{n,d}^{-1}(\D_d)$ with Maslov index 0.
\end{proposition}

\begin{proof}
    We view $\phi'$ as the image of a holomorphic map $ (u,\partial u): D \to \what_{n,d}^{-1}(\D_d)$, with boundary conditions on $\prod_{h}\Dhat_{i_h}^{-}\cup\prod_{h}\Dhat_{j_h}\cup\prod_{h}\Dhat_{k_h}^{+}$. This satisfies $\phi'_{h}=\textnormal{pr}_{h}(u(D))$, where $\textnormal{pr}_{h}$ denotes the projection to the $h^{\textnormal{th}}$ factor. We compute the Maslov index of $\phi'$, and verify that this is 0. Indeed, the Maslov index of the disjoint union of discs is the sum of the individual Maslov indices. By Proposition~\ref{prop:phiprimedisjoindiagonal}, the triangles $\{\phi'_h\}_h$ are pairwise disjoint. Moreover, each individual triangle is rigid, and has Maslov index 0. We again remark the similarity with \cite[Lemma~2.6]{LMW08}.
\end{proof}

\begin{corollary}\label{cor:phiphiprime}
    Fix $I,J,K \in \I_{n,d}$ so that $(\Dhat_I^{-},\Dhat_J,\Dhat_K^{+})$ is composable, and let $\phi'$ be as in Definition~\ref{def:phi'}. Then $\phi=\phi'$.
\end{corollary}

\begin{proof}
    The holomorphic triangle $\phi$ is the unique (by Proposition~\ref{prop:Gabungraded}) rigid triangle with boundary conditions on $\Dhat_{I}^{-}\cup\Dhat_{J}\cup\Dhat_{K}^{+}$. These are exact Lagrangian isotopic to \linebreak ${\prod_{h}\Dhat_{i_h}^{-}\cup\prod_{h}\Dhat_{j_h}\cup\prod_{h}\Dhat_{k_h}^{+}}$ by Lemma~\ref{lem:productsofperturbedthimbles}. The claim follows from the vanishing of the Maslov index of $\phi'$ in Proposition~\ref{prop:vanishingMaslov}.
\end{proof}

We conclude that holomorphic triangles with boundary on $\Dhat_{I}\cup\Dhat_{J}\cup\Dhat_{K}$, for \linebreak ${I,J,K\in\I_{n,d}}$, are disjoint from the diagonal locus $\Delta$ in $\C^d$.

\begin{proposition}\label{prop:disjoindiagonal}
     $\phi$ is disjoint from $\Delta \subset \C^d$.
\end{proposition}

\begin{proof}
    This follows from Proposition~\ref{prop:phiprimedisjoindiagonal} and Corollary~\ref{cor:phiphiprime}.
\end{proof}

The following proposition allows us to compute the endomorphism algebra $\BFS_{n,d}$ in $\F(\f_{n,d})$.

\begin{proposition}\label{lem:quotientalgebra}
    Let $I,J \in \I_{n,d}$. For an appropriate choice of grading of the objects $\{\Dt_I\}_I$, there is an isomorphism of $A_{\infty}$-algebras:
    \begin{equation}\label{eq:isoalgebrasquotient}
        \bigoplus_{I,J}\hom_{\F(\fhat_{n,d})}(\Dhat_I,\Dhat_J) \xrightarrow{\text{ $\cong$ }} \bigoplus_{I,J}\hom_{\F(\f_{n,d})}(\Dt_I,\Dt_J).
    \end{equation}
\end{proposition}

\begin{proof}
    The obstructions to the existence and uniqueness of a canonical grading on $\Sym^d(\C)$ both vanish. By Proposition~\ref{prop:productsofthimbles}, each thimble $\Dhat_I$ in $\C^d$ is exact Lagrangian isotopic to a product of thimbles in $\C$. It follows that $\Dt_I=\bpi(\Dhat_I)$ is the (unordered) product of thimbles in $\C$. Each thimble is contractible, so each $\Dt_I$ admits a choice of a grading structure, which is unique up to suitable shifts. Equip the objects $\{\Dt_I\}_I$ with such grading. 

    For $I\in\I_{n,d}$, each $\Dhat_I$ is disjoint from the diagonal locus $\Delta \subset \C^d$ by Proposition~\ref{prop:actiononthimbles}. By Proposition~\ref{prop:disjoindiagonal} the holomorphic triangles contributing to the compositions in the \linebreak $A_{\infty}$-algebra structure of the left-hand side of \eqref{eq:isoalgebrasquotient} are disjoint from the diagonal locus in $\C^d$, i.e.\ from the branch locus of the branched cover $\bpi$. By Proposition~\ref{prop:Gabungraded}, these are the only polygons contributing to the $A_{\infty}$-products on the left-hand side of the equation. We can now restrict $\bpi$ to the complement of the diagonal locus, where it is a $d!$-to-$1$ cover. Since images under $\bpi$ of pseudo-holomorphic polygons are pseudo-holomorphic, the images of polygons contributing to the $A_{\infty}$-products on the left-hand side of the equation contribute to the $A_{\infty}$-products on the right-hand side.
    
    Conversely, suppose there is a pseudo-holomorphic polygon contributing to the \linebreak $A_{\infty}$-products on the right-hand side of \eqref{eq:isoalgebrasquotient}. Since  each $\Dhat_I$ is disjoint from the diagonal locus in $\C^d$, each image $\Dt_I$ under $\bpi$ is itself is disjoint from the diagonal locus of $\Sym^d(\C)$. Following Seidel's construction on finite coverings in \cite[Section~8b]{Seidel15}, the right-hand side of \eqref{eq:isoalgebrasquotient} admits a splitting
    \begin{equation*}
        \bigoplus_{I,J}\hom_{\F(\f_{n,d})}(\Dt_I,\Dt_J) \cong \bigoplus_{I,J} \bigoplus_{\sigma \in \S_d}\hom_{\F(\f_{n,d})}(\Dt_I,\Dt_J)_{\sigma}
    \end{equation*}
    where each $\sigma$-summand consists of those $x \in \Dt_I \cap \Dt_J$, whose unique lift $\hat{x}$ satisfies ${\hat{x} \in \Dhat_I \cap \sigma \Dhat_J}$. Fixing $I$ and $J$,
    \begin{equation}\label{eq:Seidelsplitting}
        \bigoplus_{\sigma \in \S_d}\hom_{\F(\f_{n,d})}(\Dt_I,\Dt_J)_{\sigma} \xrightarrow{\cong} \bigoplus_{\sigma \in \S_d}\hom_{\F(\fhat_{n,d})}(\Dhat_I,\sigma \Dhat_J), \qquad x \mapsto \hat{x}
    \end{equation}
    defines an isomorphism of vector spaces, and by Lemma~\ref{lem:injectvspaces}, the right-hand side of \linebreak \eqref{eq:Seidelsplitting} is isomorphic to $\hom_{\F(\fhat_{n,d})}(\Dhat_I, \Dhat_J)$. In particular, given any generator $x$ of \linebreak $\hom_{\F(\f_{n,d})}(\Dt_I,\Dt_J)$, we can uniquely lift it to a generator $\hat{x}$ of $\hom_{\F(\fhat_{n,d})}(\Dhat_I,\Dhat_J)$. Moreover, given any pseudo-holomorphic polygon contributing to the $A_{\infty}$-products on the right-hand side of \eqref{eq:isoalgebrasquotient}, this can be lifted to a pseudo-holomorphic polygon contributing to the $A_{\infty}$-products on the left-hand side of the same equation \cite[Section~8b]{Seidel15}.
\end{proof}

\begin{theorem}\label{thm:mainthm2}
  There is a quasi-equivalence of triangulated $A_{\infty}$-categories
    \begin{equation*}
        \F(\f_{n,d}) \xrightarrow{\text{ $\simeq$ }} \perf (\B_{n,d}),
    \end{equation*}
    induced by the isomorphism of $A_{\infty}$-algebras:
    \begin{equation*}
        \BFS_{n,d} \xrightarrow{\text{ $\cong$ }}  \B_{n,d}.
    \end{equation*}
\end{theorem}

\begin{proof}
    For $I\in \I_{n,d}$, denote by $P_I$ the indecomposable projective $\B_{n,d}$-module associated to the vertex $I$ of the quiver underlying the $\kk$-algebra. Proposition~\ref{lem:quotientalgebra}, together with Proposition~\ref{prop:Gabungraded}, gives a chain-level description of the morphisms. Since for any pair of objects the morphism complexes are either one or zero-dimensional, all differential trivially vanish, which gives a description of the morphisms at the cohomological-level. The bijection
    \begin{equation*}
        \Dt_I \leftrightarrow P_I,
    \end{equation*}
    between the distinguished collection of generators, together with Proposition~\ref{lem:quotientalgebra}, gives the desired isomorphism.
    
    As in the proof of Theorem~\ref{thm:mainthm1}, the two triangulated categories are generated by the constructed collection of thimbles and by projective $\B_{n,d}$-modules. $\F(\f_{n,d})$ admits a tilting object by a combination of Propositions~\ref{prop:exccollectionf}, \ref{lem:quotientalgebra} and \ref{prop:Gabungraded}. In particular, $\BFS_{n,d}$ and $\B_{n,d}$ are isomorphic as $A_{\infty}$-algebras concentrated in degree 0. The claim follows from Keller's result \cite[Theorem~3.8]{Keller06}.
\end{proof}

\subsection{The derived equivalence \texorpdfstring{$\F(\f_{n,d}) \simeq \W_n^{d}$}{Ffnd Wnd}}\label{section:F-W}
The main result of this section is Theorem~\ref{thm:mainthm3}. We highlight the fact that this is a very natural generalisation of \cite[Theorem~3.4]{DiDedda23} and indeed the proof adapts in a straightforward way.

We introduce a full exceptional collection of objects in the partially wrapped Fukaya category $\W_{n}^d$ defined by Auroux \cite{Auroux1}. Let $(\D,\Lambda_n)$ be the disc with a set of $n+1$ points on its boundary. Following the clockwise orientation of the boundary, label each component of $\D\setminus \Lambda_n$ with integers between $0$ and $n$.

\begin{definition}\label{def:arcsLI}
    For $\alpha \in \{0,\dots,\frac{n-2}{2}\}$ (resp.\ $\alpha \in \{0,\dots,\frac{n-1}{2}\}$) for $n+1$ odd (resp.\ even) and $\beta \in \{n-\alpha, n-\alpha-1\}$ (resp.\ $\beta \in \{n-\alpha, n-\alpha+1\}$), denote by $\L_i$ the isotopy class of (exact) Lagrangian arcs with endpoints on the boundary components labelled $\alpha$ and $\beta$, satisfying $i=\beta-\alpha$. For $I\in \I_{n,d}$, denote by $\L_I:=\L_{i_1}\times \dots \times \L_{i_d}$ the product of these arcs.
\end{definition}

We depicted the Lagrangian arcs in Definition~\ref{def:arcsLI} in Figure~\ref{fig:swingingarcs}, for $n=4,5$. Equipped with arbitrary brane structure, $\L_i$ defines an object of the partially wrapped Fukaya category $\W_n^1$. Similarly, equipped with arbitrary brane structure, $\L_I$ is an object of the partially wrapped Fukaya category $\W_n^d$. As the collection of arcs $\{\L_i\}_i$ satisfies the hypotheses of Theorem~\ref{thm:Auroux}, $\{\L_I\}_I$ is a full collection for $\W_n^d$. 

\begin{definition}\label{def:BWnd}
    Define $\BW_{n,d}$ to be the endomorphism algebra (in $\W_n^d$) over $\kk$ of the collection of arcs in Definition~\ref{def:arcsLI}, for $I\in \I_{n,d}$.
\end{definition}

\begin{notation*}
    Throughout this section, we write $\L$ for one of the arcs (or products of arcs) defined above, and $L$ for an unspecified arc (or product of arcs) in $\D\setminus\Lambda_n$.
\end{notation*}

\begin{figure}[t]
    \begin{minipage}{.45 \textwidth}
        \centering
        \begin{tikzpicture}[scale=.8,align=center, v/.style={draw,shape=circle, fill=black, minimum size=1.2mm, inner sep=0pt, outer sep=0pt}, font=\small, label distance=1pt,
            ]
            \draw (0,0) circle (3cm);
    
            \node[v] at(18:3cm){};
            \node[v] at(90:3cm){};
            \node[v] at(162:3cm){};
            \node[v] at(234:3cm){};
            \node[v] at(306:3cm){};

            \node at (54:3.4cm){0};
            \node at (342:3.4cm){1};
            \node at (270:3.4cm){$2=\frac{n}{2}$};
            \node at (198:3.4cm){$3$};
            \node at (126:3.5cm){$4=n$};
        
            \draw[color=blue] (55:3cm)[inner sep=0pt] to[bend left=30] node[pos=.5, fill=white] {$\L_4$} (125:3cm);
             \draw[color=blue] (36:3cm)[inner sep=0pt] to[bend left=10] node[pos=.5, fill=white] {$\L_3$} (195:3cm);
             \draw[color=blue] (210:3cm)[inner sep=0pt] to[bend left=10]node[pos=.5, fill=white] {$\L_2$} (-15:3cm);
             \draw[color=blue] (270:3cm)[inner sep=0pt] to[bend left=30] node[pos=.5, fill=white] {$\L_1$} (-30:3cm);
    
            \end{tikzpicture}
            \end{minipage}
            \qquad
            \begin{minipage}{.45 \textwidth}
                \centering
                \begin{tikzpicture}[scale=.8,align=center, v/.style={draw,shape=circle, fill=black, minimum size=1.2mm, inner sep=0pt, outer sep=0pt}, font=\small, label distance=1pt, every loop/.style={distance=1cm, label=right:}
                    ]
        
                    \draw (0,0) circle (3cm);
                    \node[v] at(30:3cm){};
                    \node[v] at(90:3cm){};
                    \node[v] at(150:3cm){};
                    \node[v] at(210:3cm){};
                    \node[v] at (270:3cm){};
                    \node[v] at(330:3cm){};

                    \node at (60:3.4cm){0};
                    \node at (0:3.4cm){1};
                    \node at (-60:3.6cm){2=$\frac{n-1}{2}$};
                    \node at (-120:3.6cm){3=$\frac{n+1}{2}$};
                    \node at (180:3.6cm){$4$};
                    \node at (120:3.5cm){$5=n$};
                        
                    \draw[color=blue] (60:3cm)[inner sep=0pt] to[bend left=30] node[pos=.5, fill=white] {$\L_5$} (120:3cm);
                    \draw[color=blue] (130:3cm)[inner sep=0pt] to[bend right=15] node[pos=.5, fill=white] {$\L_4$} (10:3cm);
                    \draw[color=blue] (180:3cm)[inner sep=0pt] to[bend right=0] node[pos=.5, fill=white] {$\L_3$} (00:3cm);
                    \draw[color=blue] (190:3cm)[inner sep=0pt] to[bend left=15] node[pos=.5, fill=white] {$\L_2$} (-40:3cm);
                    \draw[color=blue] (-60:3cm)[inner sep=0pt] to[bend right=30] node[pos=.5, fill=white] {$\L_1$} (-120:3cm);
    
    \end{tikzpicture}
    \end{minipage}
\caption{The generators $\L_{1},\dots,\L_{n}$ of $\W_n^{1}$, for $n+1=5$ (left) and $n+1=6$ (right).}
\label{fig:swingingarcs}
\end{figure}

\begin{proposition}\label{lem:dgaBW}
    $\BW_{n,d}$ is a differential graded $\kk$-algebra with vanishing differential. Moreover, there is an appropriate choice of gradings of the objects underlying $\BW_{n,d}$ such that this is concentrated in degree 0.
\end{proposition}

\begin{proof}
    Fix arbitrary gradings for each morphism (choice unique up to global shift). The collection $\{\L_i\}_i$ satisfies the assumption that each disc in $\D \setminus (\cup \L_i)$ contains at least one stop, so by \cite[Proposition~3.6]{Auroux2} the $A_{\infty}$-algebra has vanishing higher $A_{\infty}$-products. The fact that the differential vanishes follows from the following, more general claim: given an arbitrary collection $\{L_i\}_i$ of arcs in $(\D,\Lambda_n)$ such that each component of ${\partial \D\setminus \Lambda_n}$ contains at most two arc endpoints, and given $L=\prod_{h=1}^d L_{i_h}$ and $L'=\prod_{h=1}^d L_{i'_h}$, $\hom_{\W_n^d}(L,L')$ has vanishing differential. The proof of the claim is as follows. By \cite[Proposition~11]{Auroux1}, $\hom_{\W_n^d}(L,L')$ admits a basis indexed by bijective maps from ${S=\{i_1,\dots,i_d\}}$ to ${T=\{i'_1,\dots,i'_d\}}$, i.e.\ by appropriate strands diagrams with $n$ strands in $d$ places (as in \eqref{eq:strandshoms}). Suppose that two of these generators $x,y\in\hom_{\W_n^d}(L,L')$ are such that $d(x)=\sum_i y_i$, and $y=y_i$ for some $i$. Equivalently, we assume that there are two bijections $\phi_1, \phi_2:S \to T$ (indexing $x$ and $y$ respectively) and some $i < i'$ with
    \begin{equation*}
        \begin{cases*}
            \phi_1(i)>\phi_1(i') \\ 
            \phi_1(j)=\phi_2(j') \quad &\text{ for $(j,j')=(i,i')$ or $(j,j')=(i',i)$} \\
            \phi_1(j)=\phi_2(j) \quad & \text{ otherwise},
        \end{cases*}
    \end{equation*}
    such that $\hom_{\W_n^1}(L_{i}, L_{\phi_h(i)})$ and $\hom_{\W_n^1}(L_{i'}, L_{\phi_h(i')})$ do not vanish for $h=1,2$. In particular, there are three distinct arc endpoints contained in the same component of $(\D,\Lambda_n)$, which are the endpoints of $L_i,L_{\phi_1(i')}, L_{\phi_1(i)}$, or $L_{i'},L_{\phi_1(i)}, L_{\phi_1(i')}$ depending respectively on whether $\phi_1(i')\neq i$ or $\phi_1(i)\neq i'$.

    Finally, the grading of each morphism is determined by the number of crossings in the strands algebra isomorphic to it. By the computations above, no such crossing can take place when there are at most two arc endpoints in each component of $\D\setminus \Lambda_n$, so $\BW_{n,d}$ is concentrated in degree 0.
\end{proof}

\begin{theorem}\label{thm:mainthm3}
    There is a quasi-equivalence of triangulated $A_{\infty}$-categories
    \begin{equation*}
        \F(\f_{n,d}) \xrightarrow{\text{ $\simeq$ }} \W_n^{d},
    \end{equation*}
    induced by the isomorphism of $A_{\infty}$-algebras concentrated in degree 0:
    \begin{equation*}
        \BFS_{n,d} \xrightarrow{\text{ $\cong$ }} \BW_{n,d}.
    \end{equation*}
\end{theorem}

\begin{proof}
    For $I\in \I_{n,d}$, consider the bijection
    \begin{equation*}
        D_I \leftrightarrow \L_I
    \end{equation*}
    between the full exceptional collections in the respective categories. By Proposition~\ref{lem:dgaBW}, $\BW_{n,d}$ has vanishing differential. In particular, this gives a cohomological description of the endomorphism algebra of the collection of Lagrangians. In order to prove that $\BFS_{n,d}$ and $\BW_{n,d}$ are isomorphic as graded $\kk$-algebras, it suffices to show that there is a morphism from $\L_I$ to $\L_J$ whenever $I=(i_1,\dots,i_d)\leq J=(j_1,\dots,j_d)$, $\lvert j_h - i_h \rvert \leq 1$, and equality holds for arbitrarily many $h\in \{1,\dots, d\}$, in which case the corresponding index $i_h$ is odd. Indeed, if $\mathbf{h}=(h_1,\dots,h_k)$ is a subset of $\{1,\dots,d\}$ of size $k$, $\L_{i_{\mathbf{h}}}:=\L_{i_{h_1}}\times \dots \times \L_{i_{h_k}}$ denotes an object of $\W_n^{k}$, and $\L_{I\setminus i_{\mathbf{h}}}$ is the object of $\W_n^{d-k}$ such that $\L_I=\L_{i_{\mathbf{h}}} \times \L_{I\setminus i_{\mathbf{h}}}$, then there is a Reeb chord giving rise to a morphism
    \begin{equation*}
        \L_{I\setminus i_{\mathbf{h}}} \times \L_{i_{h_1}} \times \dots \times \L_{i_{h_k}} \to  \L_{I\setminus i_{\mathbf{h}}} \times \L_{i_{h_1}\pm 1} \times \dots \times \L_{i_{h_k}\pm 1}
    \end{equation*}
    exactly when the latter is a non-zero object of $\W_n^d$, and $h_{\ell}$ is odd for all $\ell=1,\dots,k$. This is given by the product of all such $x_{i_h}$ (each a Reeb chord from $\L_{i_{h}}$ to $\L_{i_h \pm 1}$), and the identity morphisms on the remaining components. These morphisms are in one-to-one correspondence with paths of length $k$ in the quiver with relations underlying $\B_{n,d}$, and the latter are in bijection with morphisms between indecomposable projective \linebreak $\B_{n,d}$-modules. That is, there is a morphism $P_I\to P_J$ exactly when there is a path from  the vertex labelled $I$ to $J$. By Theorem~\ref{thm:mainthm2}, the last statement is true precisely when there is a morphism from $D_I$ to $D_J$ in $\F(\f_{n,d})$. One can similarly check compositions, by varying $k$ in $\{1,\dots,d\}$.
\end{proof}

\begin{remark}
    Theorem~\ref{thm:mainthm3} is an algebraic proof of a quasi-equivalence between two geometrically defined $A_{\infty}$-categories, namely Seidel's Fukaya--Seidel category and Auroux' partially wrapped Fukaya category. We will provide in Section~\ref{sec:appendix} a heuristic but more geometric justification for this equivalence.
\end{remark}

Corollary~\ref{cor:maincor} is an immediate consequence of Theorem~\ref{thm:mainthm3} and \cite[Theorem~1]{DJL}. We also have the following.

\begin{proposition}
    There is a series of iterated Hurwitz moves on the collection of vanishing paths $\lambdabold$ described in Proposition~\ref{prop:exccollectionf} giving rise to a new distinguished collection of thimbles, whose endomorphism algebra is quasi-isomorphic (as an $A_{\infty}$-algebra concentrated in degree 0) to $\A_{n,d}$.
\end{proposition}

\begin{proof}
   This can be proven by induction on $d$ (inductively on $n$ for a fixed $d$), and \cite[Proposition~3.15]{DiDedda23} constitutes the base step. The claim follows from a careful identification of the subcollection $\{\lambda_I \mid I\in \I_{n-1,d}\}$ (resp.\ $\{\lambda_I \mid I\in \I_{n,d}, i_d=n\}$) of $\lambdabold$ with the vanishing paths associated to $\w_{n-1,d}$ (resp.\ $\w_{n-1,d-1}$) for $n$ odd (resp.\ even), and by an induction argument on both subcollections.
\end{proof}

\subsection{A tilting object for \texorpdfstring{$\A_{n,d}$}{And}}\label{sec:tiltingobject}
Fix $n \geq d$. We consider the full exceptional collection in $\W_n^d$ given by the $d$-fold product of the collection of arcs in \cite[Figure~1]{DJL}, and we denote by $\AW_{n,d}$ its endomorphism algebra over $\kk$. A useful application of Corollary~\ref{cor:maincor} is the following. We can use the derived equivalence of the $A_{\infty}$-algebras $\BW_{n,d}$ and $\AW_{n,d}$ and the triangles provided by Auroux in Lemma~\ref{lem:Aurouxtriangles} to produce a (triangulated) resolution of the former by the latter. Under the equivalences $\perf(\B_{n,d}) \simeq \W_n^d$ (a combination of Theorems~\ref{thm:mainthm2} and \ref{thm:mainthm3}) and $\perf(\A_{n,d}) \simeq \W_n^d$ (\cite[Theorem~1]{DJL}), such resolution gives rise to a \emph{tilting object} for $\A_{n,d}$, in the sense of \cite[Theorem~6.4(e)]{Rickard}. We give this as follows.

\begin{proposition}\label{prop:triangulatedres}
    For fixed $n,d$, $\A_{n,d}$ admits a tilting object $T=\bigoplus_{I\in \I_{n,d}}\K_I$, whose endomorphism algebra is isomorphic to $\B_{n,d}$. Moreover, each complex of projective $\A_{n,d}$-modules $\K_I$ admits a triangulated resolution as in Figure~\ref{eq:bigtriangulatedresolution}, where each $\alpha_h<\beta_h$ is uniquely determined by $\beta_h-\alpha_h=i_h$ and $\beta_h \in \{n-\alpha_h, n-\alpha_h-1\}$ for $n+1$ odd (resp.\ ${\beta_h \in \{n-\alpha_h, n-\alpha_h+1\}}$ for $n+1$ even), for all $h=1,\dots,d$, with additional convention that $\P_{j_1,\dots,j_d}=0$ if any of the two indexes coincide. The triangle has side length $2^d$, and the $r^{\textnormal{th}}$ entry in the bottom row is $\P_{x_1,\dots,x_d}$, where $x_h=\alpha_h$ exactly when $r\equiv 1,2,3,\dots,2^{h-1}\mod 2^h$ (otherwise $x_h=\beta_h$), with the additional (and only local) convention that $\P_{\dots,y,\dots,x,\dots}:=\P_{\dots,x,\dots,y,\dots}$ if $y>x$.
\begin{figure}[t]
    \centering
        \begin{tikzpicture}
    \node at(-5.4,0){$\P_{\alpha_d,\dots,\alpha_1}$};
        \draw[->,shorten >= .8cm,shorten <= .8cm] (-5.4,0) to (-3,0);
    \node at(-3,0){$\P_{\alpha_d,\dots,\beta_1}$};
    \node at (-1,0) {$\cdots$};
    \node at (0,0) {$\cdots$};
    \node at (1,0) {$\cdots$};
    \node at(3,0){$\P_{\alpha_1,\dots,\beta_{d}}$};
        \draw[->,shorten >= .8cm,shorten <= .8cm] (3,0) to (5.4,0);
    \node at(5.4,0){$ \P_{\beta_1,\dots,\beta_{d}}$};

        \draw[->,shorten >= .6cm,shorten <= .cm] (-4,1) to (-5.4,0);
        \draw[->,shorten >= .cm,shorten <= .3cm] (-3,0) to (-3,1);
            \node at (-2.75,.625) {$1$};
    \node at (-1,1) {$\cdots$};
    \node at (0,1) {$\cdots$};
    \node at (1,1) {$\cdots$};
        \draw[<-,shorten >= .cm,shorten <= .3cm] (3,0) to (3,1);
        \draw[->,shorten >= .cm,shorten <= .6cm] (5.4,0) to (4,1);
            \node at (4.6,.9) {$1$};

    \node at (-3,2) {\reflectbox{$\ddots$}};
    \node at (-1.5,2) {\reflectbox{$\ddots$}};
    \node at (0,2) {$\cdots$};
    \node at (1.5,2) {$\ddots$};
    \node at (3,2) {$\ddots$};

        \draw[->,shorten >= .6cm,shorten <= .cm] (-1.6,3) to (-3,2);
        \draw[->,shorten >= .cm,shorten <= .3cm] (-1,2) to (-1,3);
            \node at (-.75,2.625) {$1$};
        \draw[<-,shorten >= .cm,shorten <= .3cm] (1,2) to (1,3);
        \draw[->,shorten >= .cm,shorten <= .6cm] (3,2) to (1.6,3);
            \node at (2.2,2.9) {$1$};
    
    \node at(-1.2,3.4){$X_1$};
      \draw[->,shorten >= .4cm,shorten <= .4cm] (-1.2,3.4) to (1.2,3.4);
    \node at(1.2,3.4){$X_2$};

    \draw[->,shorten >= .4cm,shorten <= .4cm] (0,4.3) to (-1.2,3.4);
    \node at(0,4.3){$\K_I $};
    \draw[->,shorten >= .4cm,shorten <= .4cm] (1.2,3.4) to (0,4.3);
        \node at (.75,4.05) {$1$};
    \end{tikzpicture}
        \caption{Triangulated resolution of $\K_I$ by projective $\A_{n,d}$-modules.}
        \label{eq:bigtriangulatedresolution}
    \end{figure}
\end{proposition}

\begin{proof}
    Fix $I\in \I_{n,d}$, and let $K_I$ be the object in $\perf(A_{n,d})$ defined by the triangulated resolution in Figure~\ref{eq:bigtriangulatedresolution}. The claim follows from the identification, under the derived equivalence $\perf(\A_{n,d})\xrightarrow{\simeq} \W_n^{d}$ in \cite[Theorem~1]{DJL}, between $\K_I$ and an object in $\W_n^{d}$ which is quasi-isomorphic to $\L_I$. The quasi-isomorphism makes use of Lemma~\ref{lem:Aurouxtriangles}, and follows from the isomorphism $\B_{n,d}\cong \BW_{n,d}$ (a combination of Theorems~\ref{thm:mainthm2} and \ref{thm:mainthm3}).
\end{proof}

We note that Proposition~\ref{prop:triangulatedres} recovers the content of \cite[Section~3.2.2]{DiDedda23} for $d=2$.

\section{Type A symplectic higher Auslander correspondence}\label{section:3}
The aim of this section is to present and motivate what we will call \emph{symplectic higher Auslander correspondence} for the linearly oriented $A_n$ quiver. This construction mimics the inductive presentation of Morita equivalence classes of higher Auslander algebras of type A detailed in Section~\ref{sec:Auslbackground}.

We recall from Section~\ref{sec:Auslbackground} that higher Auslander algebras of type A are inductively constructed using the higher Auslander correspondence in the following way. Given the $\kk$-algebra $\A_{n,d}$, this admits a basic $d$-cluster-tilting module $M$. Its degree 0 endomorphism algebra $\End_{\A_{n,d}}(M)$ is isomorphic to $\A_{n+1,d+1}$, and the higher Auslander correspondence is realised by
\begin{equation}\label{eq:ch5typeAhAc}
    \A_{n,d} \mapsto \End_{\A_{n,d}}(M)=\A_{n+1,d+1}.
\end{equation}
In this section, we geometrically replicate this construction, using Fukaya\textendash Seidel categories and full exceptional collections of thimbles as our objects. Given a Fukaya\textendash Seidel category $\F(\Phi_{n,d})$, and a strong full exceptional collection $\Dbold$ of thimbles, whose endomorphism algebra is isomorphic to $\A_{n,d}$ (see Setup~\ref{setup:sympdata} below), we construct a new collection of objects $\Dbold^{+}$ of $\F(\Phi_{n,d})$. The degree 0 endomorphism algebra of $\Dbold^{+}$ is isomorphic to $\A_{n+1,d+1}$, and we think of the direct sum of the objects in $\Dbold^{+}$ as corresponding to the basic module $M$ from \eqref{eq:ch5typeAhAc}. By attaching handles along the boundaries of this new collection, we obtain a new Weinstein domain $\Sigma=\Sigma_{n+1,d+1}$, equipped with a collection $\bm{\Upsilon}$ of Lagrangian spheres. This is now the generic fibre of a new Lefschetz fibration $\Phi_{n+1,d+1}$, equipped with a strong full collection $\Deltabold$ of thimbles. The upshot of this construction is that the Fukaya\textendash Seidel categories $\F(\Phi_{n+1,d+1})$ and $\F(\f_{n+1,d+1})$ are derived equivalent, and the degree 0 endomorphism algebra of $\Deltabold$ in $\F(\Phi_{n+1,d+1})$ is isomorphic to the degree 0 endomorphism algebra of $\Dbold^{+}$ in $\F(\Phi_{n,d})$. What we refer to \emph{symplectic} higher Auslander correspondence is the inductive construction
\begin{equation*}
    \Dbold \mapsto \End^0_{\F(\Phi_{n,d})}(\Dbold^{+}) \cong \End_{\F(\Phi_{n+1,d+1})}(\Deltabold).
\end{equation*}
We expect, but do not prove at this time, the pair $(\X_{n+1,d+1},\Phi_{n+1,d+1})$ in Theorem~\ref{thm:introthm4} to be Weinstein deformation equivalent to $(\Sym^{d+1}(\C),\f_{n+1,d+1})$. Though this geometric identification is beyond the scope of this article, we provide partial justification of this expectation in Section~\ref{sec:appendix}.

Throughout this section, we denote by $\N_{n,d}^{+1}$ the subset of $\N_{n+1,d+1}$ given by
    \begin{equation*}
        \N_{n,d}^{+1}:=\{(1,I+1)=(1,i_1+1,\dots,i_d+1)\in \N_{n+1,d+1}\mid I=(i_1,\dots,i_d) \in \N_{n,d}\}.
    \end{equation*}

\begin{definition}\label{def:Kvect}
    Let $I \in \N_{n+1,d+1}$, and fix some $h \in \{1,\dots,d\}$. We define $K_{\vec{i_h}}$ to be the following complex of projective $\A_{n,d}$-modules:
    \begin{equation*}
       K_{\vec{i_h}}:= P_{\widehat{i_h}} \to P_{\widehat{i_{h-1}}} \to \dots \to P_{\widehat{i_2}} \to P_{\widehat{i_1}},
    \end{equation*}
    where the morphisms are the unique ones determined by the intertwining conditions \eqref{eq:intertwining}, and the last term is fixed to be in degree 0.
\end{definition}

For $I \in \N_{n+1,d+1}$ and $J \in \N_{n,d}^{+1}$ respectively, let $P_I$ be the associated indecomposable $d$-cluster tilting $\A_{n,d}$-module, and $\P_J$ be the indecomposable projective $\A_{n,d}$-module, as in Section~\ref{sec:Auslbackground}. Recall the total order of the index set $\N_{n,d}$ in Definition~\ref{def:Nnd}.

\begin{lemma}\label{prop:vanishingDehntwists}
    Fix $I \in \N_{n+1,d+1}$, and fix some $h \in \{2,\dots,d\}$. If $K_{\vec{i_h}}$ is as in Definition~\ref{def:Kvect} and viewing $\P_J$ as a complex of projective $\A_{n,d}$-modules concentrated in degree 0, the following holds for any $J \in \N_{n,d}^{+1}$ such that $\widehat{i_{h+1}} < J < \widehat{i_{h}}$:
     \begin{equation*}
        \Hom_{\perf(\A_{n,d})}(\P_J,K_{\vec{i_h}})=0.
    \end{equation*}
\end{lemma}

\begin{proof}
    This is a direct computation relying on the total order of the indecomposable projective $\A_{n,d}$-modules and the intertwining conditions \eqref{eq:intertwining} established in Section~\ref{sec:Auslbackground}, and it amounts to checking that $\Hom_{\A_{n,d}}(\P_J,\P_{\widehat{i_h}})$ vanishes exactly when $\Hom_{\A_{n,d}}(\P_J,\P_{\widehat{i_{h-1}}})$ does. We recall that the total order of indecomposable projective $\A_{n,d}$-modules (given in Section~\ref{sec:Auslbackground}) is so that $\P_I < \P_J$ when $I<J$ in $\N_{n,d}$, as per Definition~\ref{def:Nnd}. Indeed, the condition $\P_{\widehat{i_{h+1}}} < \P_J$ in particular implies that
  \begin{equation*}
    i_1 < j_1, \quad \text{ or } \quad  i_1=j_1 \text{ and } (i_2,\dots,i_{h-1},i_h,i_{h+2},\dots,i_{d+1})<(j_2,\dots,j_{d}),
  \end{equation*}
   while $\P_J < \P_{\widehat{i_h}}$ implies that
   \begin{equation*}
    j_1 < i_1, \quad \text{ or } \quad  j_1=i_1 \text{ and } (j_2,\dots,j_{d})<(i_2,\dots,i_{h-1},i_{h+1},i_{h+2},\dots,i_{d+1}).
  \end{equation*}
   Of these combinations, the only compatible ones are the latter two. Iterating, the condition $\P_{\widehat{i_{h+1}}} < \P_J < \P_{\widehat{i_h}}$ in particular implies that $j_{\ell}=i_{\ell}$ for $\ell=1,\dots,h-1$. Suppose first $\Hom_{\A_{n,d}}(\P_J,\P_{\widehat{i_{h-1}}})$ does not vanish, i.e.\ that the intertwining conditions for the pair $(\P_J,\P_{\widehat{i_{h-1}}})$ hold: 
   \begin{equation}\label{eq:staggeredintertwining}
    \underbracket{  j_1 \leq i_1 < \dots < j_{h-1} \leq   \color{black}  i_{h}}_{(*)}<  \underbracket{ j_h\leq i_{h+1} <\dots <j_d \leq i_{d+1} }_{(**)}.
  \end{equation}
   The ($*$) part of \eqref{eq:staggeredintertwining} always holds by $j_{\ell}=i_{\ell}$ for $\ell=1,\dots,h-1$. By $i_{h-1}<i_h<j_{h}$ and the ($**$) part of \eqref{eq:staggeredintertwining}, the pair $(\P_J,\P_{\widehat{i_h}})$ satisfies the intertwining conditions. The converse claim is analogous.
\end{proof}

We define the following setup. Let $\A_{n,d}$ be the $d$-dimensional Auslander algebra of type $A_n$.

\begin{setup}(Symplectic data associated to $\A_{n,d}$)\label{setup:sympdata}
    We define the \emph{symplectic data associated to $\A_{n,d}$} to be a tuple $(\X_{n,d},\Phi_{n,d},\Sigma_{n,d},\Vbold)$ satisfying the following:
    \begin{enumerate}[i.]
        \item $\X_{n,d}$ is a contractible symplectic manifold with boundary of (real) dimension $2d$;
        \item $\Phi_{n,d}:\X_{n,d} \to \D$ is an exact symplectic Lefschetz fibration to the standard disc $\D \subset \C$, with generic fibre $\Sigma_{n,d}$ above a fixed point on $\partial \D$ (which we can take to be $-i \in \partial \D$);
        \item $\Vbold:=\{V_I\}_I$, indexed and ordered by $I\in \N_{n,d}$, is a distinguished collection of vanishing cycles on $\Sigma_{n,d}$, whose associated Lefschetz thimbles $\Dbold=\{D_I\}_I$ (equipped with appropriate brane structures) form a strong full exceptional collection in $\F(f)$;
        \item $\Vbold$ is in one-to-one correspondence with indecomposable projective $\A_{n,d}$-modules, in such a way that for each $I,J \in \N_{n,d}$, there is an isomorphism of $A_{\infty}$-algebras
        \begin{equation*}
            \Hom_{\F(\Phi_{n,d})}(D_I,D_J) \xrightarrow{\cong}  \Hom_{\perf(\A_{n,d})}(P_I,P_J)
        \end{equation*}
        and an isomorphism (also of $A_{\infty}$-algebras)
        \begin{equation}\label{eq:symplecticdataiso}
            \End_{\F(\Phi_{n,d})}(\Dbold) \xrightarrow{\cong} \A_{n,d}
        \end{equation}
        underlying the derived equivalence
        \begin{equation}\label{eq:symplecticdataequiv}
            \F(\Phi_{n,d}) \xrightarrow{\simeq} \perf(\A_{n,d}).
        \end{equation}
    \end{enumerate}
\end{setup}

\begin{remark}
    As per Section~\ref{sec:FSbackground} and after taking its Liouville completion, $\Phi_{n,d}$ in Setup~\ref{setup:sympdata} can be extended to a map from an open symplectic manifold to $\C$. In this setting, it is more natural to consider $\X_{n,d}$ to have non-empty boundary.
\end{remark}

In this section, we assume the existence of Setup~\ref{setup:sympdata} associated to $\A_{n,d}$, and we construct symplectic data $(\X_{n+1,d+1},\Phi_{n+1,d+1},\Sigma_{n+1,d+1},\Upsilonbold)$ associated to $\A_{n+1,d+1}$. Of course, the existence of such data was proven explicitly in Section~\ref{section:2}. We give an alternative construction, in such a way that $\bm{\Upsilon}$ is obtained as a collection of \emph{matching spheres} constructed from $\Vbold$, and $\Sigma_{n+1,d+1}$ is obtained from a handle attachment procedure to $X_{n,d}$. We first define a collection of vanishing spheres on $\Sigma_{n,d}$.

\begin{definition}\label{def:halfmatchingspheres}
    For $I \in \N_{n,d}$, denote by $W_{1,I+1}$ the following vanishing sphere:
    \begin{equation*}
        W_{1,I+1}:=V_I.
    \end{equation*}
    For $J\in \N_{n+1,d+1}\setminus \N^{+1}_{n,d}$, define the following Lagrangian sphere $W_{J}$ in $\Sigma_{n,d}$:
    \begin{equation}\label{eq:sphereWJcone}
        W_{J}:=\tau_{V_{\widehat{j_{d+1}}}}\tau_{V_{\widehat{j_d}}}\cdots\tau_{V_{\widehat{j_2}}} V_{\widehat{j_1}}
    \end{equation}
    where $\tau$ denotes the symplectic Dehn twist, and let $\mathbf{W}$ be the collection of above constructed Lagrangian spheres in $\Sigma_{n,d}$:
    \begin{equation*}
        \mathbf{W}:=\{W_K\mid K \in \N_{n+1,d+1}\}= \{W_{1,I+1} \mid I\in \N_{n,d}\} \cup \{W_{J} \mid J\in \N_{n+1,d+1}\setminus \N^{+1}_{n,d}\}.
    \end{equation*}
\end{definition}

We claim that Lagrangian spheres $\Sigma_{n,d}$ appearing in Definition~\ref{def:halfmatchingspheres} arise from a collection of vanishing spheres. We now describe the collection of vanishing paths underlying these spheres, which we depicted in Figure~\ref{fig:twistsbyprojres} for $n+1=4$ and $d+1=2$. In figure, the straight segments $\gamma_1$, $\gamma_2$ and $\gamma_3$ are the vanishing paths associated to the spheres $W_{12}=V_1$, $W_{13}=V_2$ and $W_{14}=V_3$ respectively. The remaining paths $\gamma_{23}$, $\gamma_{24}$ and $\gamma_{34}$ are those associated to $W_{23}=\tau_{V_1}V_2$, $W_{24}=\tau_{V_1}V_3$ and $W_{34}=\tau_{V_2}V_3$ respectively. Let $\{\gamma_{I}\}_I$, indexed by $I \in \N_{n,d}$, be the distinguished collection of vanishing paths associated to $\Dbold$.

\begin{figure}[t]
    \centering
    \begin{tikzpicture}[scale=1,v/.style={draw,shape=circle, fill=black, minimum size=1mm, inner sep=0pt, outer sep=0pt},cross/.style={cross out,draw=black, minimum size=1mm},cross/.style={cross out, draw=black, minimum size=1.5mm, inner sep=0pt, outer sep=0pt},cross/.default={1pt}]
        \draw (0,1.5) ellipse (5cm and 1.5cm);
        \node[v] at (0,0) {};
        \node at (5,1.5) {$*$};
        \node[cross] at (-2,1) {};
        \node[cross] at (0,1) {};
        \node[cross] at (2,1) {};

        \draw (0,0) to (-2,1);
        \draw (0,0) to (0,1);
        \draw (0,0) to (2,1);
        \draw (0,0) to[out=150,in=-90] (-2.3,1.1) to[out=90,in=180] (-1.9,1.4) to[out=0,in=180] (0,.6) to[out=0,in=200] (2,1);
        \draw (0,0) to[out=160,in=-90] (-2.5,1) to[out=90,in=180] (-1.7,1.6) to[out=0,in=140] (0,1);
        \draw (0,0) to[out=170,in=-90] (-2.7,1) to[out=90,in=160] (.5,1.6) to[out=-20,in=160] (2,1);
        
        \node at (-1.5,1) {$\gamma_1$};
        \node at (.3,1) {$\gamma_2$};
        \node at (2.3,1) {$\gamma_3$};

        \draw (-6.5,1.4) edge[->] (-2.65,1.4);
        \draw (-6.35,1) edge[->] (-2.5,1);
        \draw (-6,.6) edge[->] (-2.15,.6);
        \node at (-6.9,1.4) {$\gamma_{34}$};
        \node at (-6.75,1) {$\gamma_{23}$};
        \node at (-6.4,.6) {$\gamma_{24}$};
\end{tikzpicture}
    \caption{The collection of vanishing paths for the vanishing spheres $\{W_K \mid K \in \N_{4,2}\}$, where $\gamma_1,\gamma_2,\gamma_3$ are the straight segments.}
    \label{fig:twistsbyprojres}
\end{figure}

\begin{proposition}
    The vanishing path $\gamma_J$ associated to each Lagrangian sphere $W_J$, for $J \in~\N_{n+1,d+1}\setminus \N^{+1}_{n,d}$, is the following:
    \begin{equation}\label{eq:inductiveclusterpaths}
        \gamma_J:=\left(\prod \tau_{\gamma_L}\right) \gamma_{\widehat{j_1}},
    \end{equation}
    where the product of Dehn twists ranges over all $L \in \N_{n,d}$ such that $L \leq \widehat{j_2}$.
\end{proposition}

\begin{proof}
    By Definition~\ref{def:halfmatchingspheres}, it suffices to prove that, for $\widehat{j_{h+1}} < L < \widehat{j_h}$ the vanishing cycles $V_L$ and $\tau_{V_{\widehat{j_h}}}\cdots \tau_{V_{\widehat{j_2}}} V_{\widehat{j_{1}}}$ are disjoint. This directly follows from Lemma~\ref{prop:vanishingDehntwists}, under the equivalence \eqref{eq:symplecticdataequiv} and the isomorphism \eqref{eq:Seidelisom} constructed under the restriction functor \eqref{eq:restrictionfunctor}.
\end{proof}

We choose isotopy classes of thimbles in $X_{n,d}$ that are in \emph{minimal position}, i.e.\ that have the minimal number of intersection points. This means taking isotopy classes of paths in $\D\setminus \{*\}$, where $*=1 \in \partial \D$, that intersect minimally. See Figure~\ref{fig:minimalposition}, for $n=3$ and $d=1$ in figure. Essentially by construction, we have the following.

    \begin{figure}[t]
        \centering
        \begin{tikzpicture}[scale=1,v/.style={draw,shape=circle, fill=black, minimum size=1mm, inner sep=0pt, outer sep=0pt},cross/.style={cross out,draw=black, minimum size=1mm},cross/.style={cross out, draw=black, minimum size=1.5mm, inner sep=0pt, outer sep=0pt},cross/.default={1pt}]
            \draw (0,1.5) ellipse (5cm and 1.5cm);
            \node at (5,1.5) {$*$};
            \node[cross] at (-2,1) {};
            \node[cross] at (0,1) {};
            \node[cross] at (2,1) {};
    
            \draw ($(0,1.5)+(-113.5:5 and 1.5)$) to (-2,1);
            \draw ($(0,1.5)+(-90:5 and 1.5)$) to (0,1);
            \draw ($(0,1.5)+(-66.5:5 and 1.5)$) to (2,1);
            \draw ($(0,1.5)+(-117.5:5 and 1.5)$) to[out=90,in=-90] (-2.3,1.1) to[out=90,in=180] (-1.9,1.4) to[out=0,in=180] (0,.6) to[out=0,in=200] (2,1);
            \draw ($(0,1.5)+(-120:5 and 1.5)$) to[out=90,in=-90] (-2.5,1) to[out=90,in=180] (-1.7,1.6) to[out=0,in=140] (0,1);
            \draw ($(0,1.5)+(-122.5:5 and 1.5)$) to[out=90,in=-90] (-2.7,1) to[out=90,in=160] (.5,1.6) to[out=-20,in=160] (2,1);
            
            \node at (-1.65,1) {$\gamma_1$};
            \node at (.3,1) {$\gamma_2$};
            \node at (2.3,1) {$\gamma_3$};
    
            \draw (-6.5,1.4) edge[->] (-2.65,1.4);
            \draw (-6.35,1) edge[->] (-2.55,1);
            \draw (-6,.6) edge[->] (-2.35,.6);
    
            \node at (-6.9,1.4) {$\gamma_{34}$};
            \node at (-6.75,1) {$\gamma_{23}$};
            \node at (-6.4,.6) {$\gamma_{24}$};
    \end{tikzpicture}
        \caption{Vanishing paths in minimal position, for $n=3$ and $d=1$.}
        \label{fig:minimalposition}
    \end{figure}

\begin{lemma}\label{lem:correspondenceCTthimbles}
    The derived equivalence \eqref{eq:symplecticdataequiv} determines an object of the Fukaya--Seidel category $\F(f)$ for each indecomposable direct summand of the $d$-cluster-tilting $\A_{n,d}$-module. Each of these objects is quasi-isomorphic to one of the thimbles constructed above.
\end{lemma}

\begin{proof}
    The thimbles $\{D_I\}_I$, $I\in \N_{n,d}^{+1}$, are determined by the indecomposable projective $\A_{n,d}$-modules by assumption. Proposition~\ref{prop:projresCT} realises each non-projective indecomposable direct summand $\CTmod_J$ of the $d$-cluster-tilting $\A_{n,d}$-module $M$ as an iterated cone of projective $\A_{n,d}$-modules in $\perf(\A_{n,d})$, via the resolution
    \begin{equation*}
        \P_{\widehat{j_{d+1}}}\to  \P_{\widehat{j_d}}\to \dots \to  \P_{\widehat{j_1}} \to  \CTmod_J \to 0.
    \end{equation*}
    On the other hand, each thimble $D_J$ is realised as the iterated twist of thimbles (in the sense of \cite[Section~5j]{SeidelBk}), explicitly given by \eqref{eq:inductiveclusterpaths} and \eqref{eq:sphereWJcone} as 
    \begin{equation*}
        \Cone\left(D_{\widehat{j_{d+1}}}\to \Cone\left(D_{\widehat{j_d}}\to\Cone(\dots \to \Cone(D_{\widehat{j_3}}\to \Cone(D_{\widehat{j_2}}\to D_{\widehat{j_1}}))\cdots)\right)\right)
    \end{equation*}
    where each cone is formed over the unique morphism corresponding (under the equivalence \eqref{eq:symplecticdataequiv}) to each unique morphism determined by the intertwining conditions, as in Proposition~\ref{prop:projresCT}. As mapping cones are determined up to quasi-isomorphism, the claim follows.
\end{proof}

\begin{proposition}\label{prop:unperturbedendalgebra}
    The following is an isomorphism of $A_{\infty}$-algebras:
    \begin{equation}\label{eq:unperturbedExts}
       \Hom^{0}_{\F(\Phi_{n,d})}(\oplus D_I)\cong\Ext_{\mathcal{U}^{(d)}}^{0}(\oplus \P_I)
    \end{equation}
    where $\Hom^{0}$ and $\Ext^{0}$ denote the degree zero morphisms in the appropriate categories, each $\P_I$ is an indecomposable projective $\A_{n+1,d+1}$-module, and each sum ranges over $\N_{n+1,d+1}$.
\end{proposition}

\begin{proof}
    Choose the unique gradings of the thimbles determined by \eqref{eq:symplecticdataiso}. The claim follows from Lemma~\ref{lem:correspondenceCTthimbles}. 
\end{proof}

There are two kinds of morphisms on the left-hand side of \eqref{eq:unperturbedExts}: those arising from (transverse) intersection points, and those arising from counter-clockwise Reeb chords along the boundary of $\D$ (stopped at $*=1\in \partial \D$).

\begin{proposition}\label{prop:characterisationmorphs}
    Fix $I,J \in \N_{n+1,d+1}$. The following statements hold.
    \begin{enumerate}[i.]
        \item A morphism $\hom(D_I,D_J)$ in $\F(\Phi_{n,d})$ arising from an intersection point is concentrated in degree 0 exactly when $I<J$.
        \item A morphism $\hom(D_I,D_J)$ in $\F(\Phi_{n,d})$ arising from a Reeb chord is concentrated in degree 0 exactly when $I<J$.
    \end{enumerate}
\end{proposition}

\begin{proof}
    The two cases are mutually exclusive by minimality of the number of intersections. To prove the first claim, we note that, fixing $I<J$ in $\N_{n+1,d+1} \setminus \N_{n,d}^{+1}$, there is an intersection point between $D_I$ and $D_J$ exactly when
    \begin{equation*}
        \P_{\widehat{i_2}} < \P_{\widehat{j_2}} < \P_{\widehat{i_1}} \leq \P_{\widehat{j_1}}
    \end{equation*}
    and the intertwining conditions \eqref{eq:intertwining} hold, in which case there is an obvious morphism between $\P_I$ and $\P_J$ in degree 0 in $\perf(\A_{n,d})$. The claim follows from Proposition~\ref{prop:unperturbedendalgebra}. The remaining intersection points occur when $I \in \N_{n,d}^{+1}$ and $J \notin \N_{n,d}^{+1}$. In this case, there is an intersection point between $D_I$ and $D_J$ exactly when $\P_{I} \leq \P_{\widehat{j_1}}$ and the intertwining conditions hold, and the corresponding morphism between $\P_I$ and $\P_J$ in $\perf(\A_{n,d})$ is in degree 0. 

    Fixing $I<J$ in $\N_{n+1,d+1} \setminus \N_{n,d}^{+1}$, the second case only occurs when
    \begin{equation*}
        \P_{\widehat{i_2}} = \P_{\widehat{j_2}} < \P_{\widehat{i_1}} < \P_{\widehat{j_1}}
    \end{equation*}
    and the intertwining conditions hold. In this case, there is an obvious morphism between $\P_I$ and $\P_J$ in degree 0 in $\perf(\A_{n,d})$. Finally, if $I \in \N_{n,d}^{+1}$ there is no Reeb chord from $\P_I$ and $\P_J$ unless $J \in \N_{n,d}^{+1}$, in which case the morphism is in degree 0 by Lemma~\ref{prop:unperturbedendalgebra}.
\end{proof}

    \begin{figure}[t]
        \centering
        \begin{tikzpicture}[scale=1,v/.style={draw,shape=circle, fill=black, minimum size=.5mm, inner sep=0pt, outer sep=0pt},cross/.style={cross out,draw=black, minimum size=1mm},cross/.style={cross out, draw=black, minimum size=1.5mm, inner sep=0pt, outer sep=0pt},cross/.default={1pt}]
            \draw (0,1.5) ellipse (5cm and 1.5cm);
            \node at (5,1.5) {$*$};
            \node[cross] at (-2,1) {};
            \node[cross] at (0,1) {};
            \node[cross] at (2,1) {};

            \draw ($(0,1.5)+(-75:5 and 1.5)$) to (-2,1);
            \draw ($(0,1.5)+(-85:5 and 1.5)$) to (0,1);
            \draw ($(0,1.5)+(-90:5 and 1.5)$) to (2,1);
            \draw ($(0,1.5)+(-100:5 and 1.5)$) to[out=150,in=-90] (-2.3,1.1) to[out=90,in=180] (-1.9,1.4) to[out=0,in=180] (0,.6) to[out=0,in=200] (2,1);
            \draw ($(0,1.5)+(-95:5 and 1.5)$) to[out=160,in=-90] (-2.5,1) to[out=90,in=180] (-1.7,1.6) to[out=0,in=140] (0,1);
            \draw ($(0,1.5)+(-105:5 and 1.5)$) to[out=170,in=-90] (-2.7,1) to[out=90,in=160] (.5,1.6) to[out=-20,in=160] (2,1);

            \node at (-1.5,1) {$\gamma_1$};
            \node at (.3,1) {$\gamma_2$};
            \node at (2.3,1) {$\gamma_3$};

            \draw (-6.5,1.4) edge[->] (-2.65,1.4);
            \draw (-6.35,1) edge[->] (-2.5,1);
            \draw (-6,.6) edge[->] (-2.15,.6);

            \node at (-6.9,1.4) {$\gamma_{34}$};
            \node at (-6.75,1) {$\gamma_{23}$};
            \node at (-6.4,.6) {$\gamma_{24}$};
    \end{tikzpicture}
        \caption{Positive perturbations on $\{W_K \mid K \in \N_{4,2}\}$.}
        \label{fig:positiveperturbations}
    \end{figure}

We define a collection of representatives of the isotopy classes of $\Dbold$. Fix $\gammabold$ to be the collection of paths \eqref{eq:inductiveclusterpaths} in minimal position.

    \begin{definition}[Positive perturbations, following {\cite[Definition~7]{Auroux1}}]\label{def:positiveperturbations}
        For each \linebreak ${I \in \N_{n+1,d+1}}$, choose representative of the isotopy class of $\gamma_I$ such that the endpoint of $\gamma_I$ on $\partial \D$ lies before that of any other $\gamma_J$ (following the counter-clockwise orientation of the boundary) exactly when $I>J$. Denote by $\gamma^{+}_I$ these choices of representatives, by $W_I^{+}$ the vanishing sphere above the endpoint of $\gamma^{+}_I$ on $\partial \D$, and by $D^{+}_I$ the corresponding representatives of the isotopy class of thimbles.
    \end{definition}

    Definition~\ref{def:positiveperturbations}, together with Propositions~\ref{prop:unperturbedendalgebra} and \ref{prop:characterisationmorphs}, allows us to realise all degree 0 morphisms in $\F(\Phi_{n,d})$ as intersection points between thimbles. See Figure~\ref{fig:positiveperturbations} for a depiction of the positive perturbations on $\gammabold=\{\gamma_{I}\}_{I\in \N_{4,2}}$, where $\gamma_{1i}:=\gamma_i$ for $i=1,2,3$. Our goal now is to construct a symplectic manifold (with boundary) $\Sigma_{n+1,d+1}$, from the pair $(\X_{n,d},\gammabold^{+})$.

\begin{definition}\label{def:Sigma}
    Define $\Sigma_{n+1,d+1}$ to be the symplectic manifold (with boundary) obtained by attaching ${n+1} \choose {d+1}$ handles of dimension $2d$ and index $d$ to $\X_{n,d}$ along the attaching spheres $\{W_I^{+} \mid I \in \N_{n+1,d+1}\}$. Denote by $D_I^{-}$ the core of each handle, by $W_I^{-}$ the attaching sphere, and by $\Upsilon_I$ each resulting matching sphere:
    \begin{equation*}
        \Upsilon_I := D_I^{+} \cup D_I^{-}.
    \end{equation*}
    The collection $\Upsilonbold$ of matching spheres is ordered by $\Upsilon_I<\Upsilon_J$ whenever $I<J$ in $\N_{n+1,d+1}$.
\end{definition}

We remark that the setting of matching spheres construction in Definition~\ref{def:Sigma} is that of the \emph{naive matching cycle construction} of Seidel (\cite[Section~16g]{SeidelBk}). We depicted the handle attachment procedure in Figure~\ref{fig:symplAuslfibre}: in figure, we projected $\Sigma_{n+1,d+1}$ to $\D$, for $n+1=4$ and $d+1=2$.

    \begin{figure}[t]
        \centering
        \begin{tikzpicture}[scale=1,v/.style={draw,shape=circle, fill=black, minimum size=.5mm, inner sep=0pt, outer sep=0pt},cross/.style={cross out,draw=black, minimum size=1mm},cross/.style={cross out, draw=black, minimum size=1.5mm, inner sep=0pt, outer sep=0pt},cross/.default={1pt}]
            \draw (0,1.5) ellipse (5cm and 1.5cm);
  
            \node[cross] at (-2,1) {};
            \node[cross] at (0,1) {};
            \node[cross] at (2,1) {};

            \draw ($(0,1.5)+(-75:5 and 1.5)$) to (-2,1);
            \draw ($(0,1.5)+(-85:5 and 1.5)$) to (0,1);
            \draw ($(0,1.5)+(-90:5 and 1.5)$) to (2,1);
            \draw ($(0,1.5)+(-100:5 and 1.5)$) to[out=150,in=-90] (-2.3,1.1) to[out=90,in=180] (-1.9,1.4) to[out=0,in=180] (0,.6) to[out=0,in=200] (2,1);
            \draw ($(0,1.5)+(-95:5 and 1.5)$) to[out=160,in=-90] (-2.5,1) to[out=90,in=180] (-1.7,1.6) to[out=0,in=140] (0,1);
            \draw ($(0,1.5)+(-105:5 and 1.5)$) to[out=170,in=-90] (-2.7,1) to[out=90,in=160] (.5,1.6) to[out=-20,in=160] (2,1);

            \draw ($(0,1.5)+(-75:5 and 1.5)$) to ($(0,1.5)+(-75:5 and 1.5)+(0,-1)$);
            \draw ($(0,1.5)+(-85:5 and 1.5)$) to ($(0,1.5)+(-85:5 and 1.5)+(0,-1)$);
            \draw ($(0,1.5)+(-90:5 and 1.5)$) to ($(0,1.5)+(-90:5 and 1.5)+(0,-1)$);
            \draw ($(0,1.5)+(-100:5 and 1.5)$) to  ($(0,1.5)+(-100:5 and 1.5)+(0,-1)$);
            \draw ($(0,1.5)+(-95:5 and 1.5)$) to ($(0,1.5)+(-95:5 and 1.5)+(0,-1)$) ;
            \draw ($(0,1.5)+(-105:5 and 1.5)$) to ($(0,1.5)+(-105:5 and 1.5)+(0,-1)$);
            
            \node[cross] at($(0,1.5)+(-75:5 and 1.5)+(0,-1)$) {};
            \node[cross] at ($(0,1.5)+(-85:5 and 1.5)+(0,-1)$) {};
            \node[cross] at ($(0,1.5)+(-90:5 and 1.5)+(0,-1)$) {};
            \node[cross] at ($(0,1.5)+(-95:5 and 1.5)+(0,-1)$){};
            \node[cross] at ($(0,1.5)+(-100:5 and 1.5)+(0,-1)$){};
            \node[cross] at ($(0,1.5)+(-105:5 and 1.5)+(0,-1)$) {};

            \draw ($(0,1.5)+(-107:5 and 1.5)$) to[out=-90,in=90] ($(0,1.5)+(-105:5 and 1.5)+(-.2,-.8)$) to[out=-90,in=180] ($(0,1.5)+(-105:5 and 1.5)+(0,-1.2)$) to[out=0,in=-90]  ($(0,1.5)+(-105:5 and 1.5)+(.2,-.8)$) to[out=90,in=-90] ($(0,1.5)+(-103:5 and 1.5)$);

            \draw ($(0,1.5)+(-102:5 and 1.5)$) to[out=-90,in=90] ($(0,1.5)+(-100:5 and 1.5)+(-.2,-.8)$) to[out=-90,in=180] ($(0,1.5)+(-100:5 and 1.5)+(0,-1.2)$) to[out=0,in=-90]  ($(0,1.5)+(-100:5 and 1.5)+(.2,-.8)$) to[out=90,in=-90] ($(0,1.5)+(-98:5 and 1.5)$);

            \draw ($(0,1.5)+(-97:5 and 1.5)$) to[out=-90,in=90] ($(0,1.5)+(-95:5 and 1.5)+(-.2,-.8)$) to[out=-90,in=180] ($(0,1.5)+(-95:5 and 1.5)+(0,-1.2)$) to[out=0,in=-90]  ($(0,1.5)+(-95:5 and 1.5)+(.2,-.8)$) to[out=90,in=-90] ($(0,1.5)+(-93:5 and 1.5)$);

            \draw ($(0,1.5)+(-92:5 and 1.5)$) to[out=-90,in=90] ($(0,1.5)+(-90:5 and 1.5)+(-.2,-.8)$) to[out=-90,in=180] ($(0,1.5)+(-90:5 and 1.5)+(0,-1.2)$) to[out=0,in=-90]  ($(0,1.5)+(-90:5 and 1.5)+(.2,-.8)$) to[out=90,in=-90] ($(0,1.5)+(-88:5 and 1.5)$);

            \draw ($(0,1.5)+(-87:5 and 1.5)$) to[out=-90,in=90] ($(0,1.5)+(-85:5 and 1.5)+(-.2,-.8)$) to[out=-90,in=180] ($(0,1.5)+(-85:5 and 1.5)+(0,-1.2)$) to[out=0,in=-90]  ($(0,1.5)+(-85:5 and 1.5)+(.2,-.8)$) to[out=90,in=-90] ($(0,1.5)+(-83:5 and 1.5)$);

            \draw ($(0,1.5)+(-77:5 and 1.5)$) to[out=-90,in=90] ($(0,1.5)+(-75:5 and 1.5)+(-.2,-.8)$) to[out=-90,in=180] ($(0,1.5)+(-75:5 and 1.5)+(0,-1.2)$) to[out=0,in=-90]  ($(0,1.5)+(-75:5 and 1.5)+(.2,-.8)$) to[out=90,in=-90] ($(0,1.5)+(-73:5 and 1.5)$);
    \end{tikzpicture}
        \caption{Projection of the fibre $\Sigma_{n+1,d+1}$ to $\D$, for $n=3$ and $d=1$.}
        \label{fig:symplAuslfibre}
    \end{figure}

Generators of morphisms between the collection of matching spheres in Definition~\ref{def:Sigma} inherit a $\Z$-grading by Proposition~\ref{prop:unperturbedendalgebra}. In fact, these are all concentrated in degree 0.

\begin{proposition}\label{prop:isomungradedpospert}
    For $I,J \in \N_{n+1,d+1}$, there is an isomorphism of graded $\kk$-algebras:
    \begin{equation}\label{eq:isomungradedpospert}
        \bigoplus_{I<J}\Hom_{\Fuk(\Sigma_{n+1,d+1})}(\Upsilon_I,\Upsilon_J)\cong \bigoplus_{I<J}\Hom_{\perf(\A_{n+1,d+1})}(\P_I,\P_J).
    \end{equation}
\end{proposition}

\begin{proof}
    The isomorphism of vector spaces
    \begin{equation*}
        \Hom_{\Fuk(\Sigma_{n+1,d+1})}(\Upsilon_I,\Upsilon_J)\cong \Hom_{\perf(\A_{n+1,d+1})}(\P_I,\P_J)
    \end{equation*}
    follows from Proposition~\ref{prop:characterisationmorphs} and Definition~\ref{def:positiveperturbations}. Compositions on the right-hand side of \eqref{eq:isomungradedpospert} give compositions on the left-hand side, as holomorphic curves contributing to $\mu_2$-products in $\F(\Phi_{n,d})$ also contribute to $\mu_2$-products in $\Fuk(\Sigma_{n+1,d+1})$. Furthermore, the process of handle attachment gives no additional holomorphic curve by a maximum principle, as all the cores of the handles are pairwise disjoint.
\end{proof}

\begin{theorem}\label{prop:uniquetopfibr}
    There exists a unique (up to Weinstein deformation equivalence) symplectic manifold $\X_{n+1,d+1}$ and a unique Lefschetz fibration $\Phi_{n+1,d+1}: \X_{n+1,d+1} \to \D$, such that $\Sigma_{n+1,d+1}$ is the generic fibre, and the framed Lagrangian spheres $\Upsilonbold$ are vanishing cycles associated to a collection of vanishing paths. Furthermore, $\X_{n+1,d+1}$ is contractible.
\end{theorem}

\begin{proof}
    We prove the claims for $d+1>2$. We note that the  $d+1=2$ case is the content of \cite[Section~1]{DiDedda23}, where we implicitly prove a stronger result: the Liouville completion of the pair $(\X_{n+1,2}, \Phi_{n+1,2})$ is isomorphic to $(\Sym^2(\C), \f_{n+1,2})$, in the sense of \cite[Section~15b]{SeidelBk}. The first part of the statement follows from the construction in \cite[Section~16e]{SeidelBk}, once one establishes the total order of the Lagrangian spheres $\Upsilon_I < \Upsilon_J$ whenever $I<J$.

    To check that $\X_{n+1,d+1}$ is contractible, it suffices to compute its homology groups in degrees $d$ and $d+1$, and check that these both vanishing (the remaining homology groups, except for $H_0$, vanish because they vanish for $\Sigma_{n+1,d+1}$, see Proposition~\ref{prop:simplex} below). The claim follows from a combination of the classical Hurewicz and Whitehead theorems \cite[Theorems~4.32 and 4.5]{Hatcher}. Note that the fundamental group of $\X_{n+1,d+1}$ vanishes for $d+1>2$ thanks to its inductive handle decomposition.
\end{proof}

Define $\Deltabold:=\{\Delta_I\}_I$ to be the collection of Lefschetz thimbles associated to the vanishing spheres $\Upsilonbold$ for the Lefschetz fibration $\Phi_{n+1,d+1}$. Equip these with appropriate brane structures, such that each morphism space is concentrated in degree 0. The following provides the desired symplectic data associated to $\A_{n+1,d+1}$.

\begin{corollary}\label{mainthmsympAus}
    There is an isomorphism of $A_{\infty}$-algebras (concentrated in degree 0)
    \begin{equation*}
        \Hom_{\F(\Phi_{n+1,d+1})}(\Delta_I,\Delta_J) \xrightarrow{\cong}  \Hom_{\perf(\A_{n+1,d+1})}(P_I,P_J)
    \end{equation*}
    and an isomorphism (also of $A_{\infty}$-algebras concentrated in degree 0)
    \begin{equation*}
        \Hom_{\F(\Phi_{n+1,d+1})}(\Deltabold) \cong \A_{n+1,d+1}
    \end{equation*}
    underlying the quasi-equivalence of triangulated $A_{\infty}$-categories
    \begin{equation*}
        \F(\Phi_{n+1,d+1}) \xrightarrow{\text{ $\simeq$ }} \perf(\A_{n+1,d+1}).
    \end{equation*}
\end{corollary}

\begin{proof}
    This follows from Proposition~\ref{prop:isomungradedpospert} and Theorem~\ref{prop:uniquetopfibr}, under \eqref{eq:Seidelisom}.
\end{proof}

\subsection{The regular fibre}
In this last section we study the regular fibre $\Sigma_{n,d}$ inductively constructed at the beginning of Section~\ref{section:3}. We assume that $\Sigma_{n,2}$ coincides with $\Sigma_{n}$ in \cite[Proposition~3.18]{DiDedda23}. Throughout this section, we denote by $\Delta_n$ the standard $n$-simplex:
\begin{equation*}
    \Delta_n:=\left\{(t_0,\dots,t_n)\in (\CR_{\geq0})^{n+1} \mid  \sum_{i=0}^{n} t_{i}=1\right\}
 \end{equation*}
 and by $\Delta_n^d$ its $d$-dimensional skeleton:
 \begin{equation*}
    \Delta_n^d:=\left\{(t_0,\dots,t_n)\in \Delta_n \mid  \sum_{h=1}^{d+1} t_{i_h}=1, \textnormal{ for some } 0\leq i_{1}<\dots<i_{d+1}\leq n \right\}.
  \end{equation*}

\begin{lemma}
    If $\X_{n,d}$ is a Liouville manifold, so is $\Sigma_{n+1,d+1}$ in Definition~\ref{def:Sigma}. In particular, $\Sigma_{n+1,d+1}$ is a Liouville hypersurface in $\X_{n+1,d+1}$, and its core is ${\Core(\X_{n,d}) \cup (\bigcup_I D^{-}_I)}$.
\end{lemma}

\begin{proof}
    The handles prescribed in Definition~\ref{def:Sigma} can be attached to $\X_{n,d}$ in a way that is compatible with its Liouville structure \cite{Weinstein}. More precisely, each handle can be equipped with a Liouville structure with respect to which the points which do not escape to infinite form the core of the handle. See \cite[Figure~1]{Weinstein}.
\end{proof}

\begin{proposition}\label{prop:simplex}
    The core of $\Sigma_{n,d}$ coincides with the $(d-1)$-dimensional skeleton of the standard $n$-simplex $\Delta_n^{d-1}$.
\end{proposition}

The proof of the above is inductive on the index $d$. We first give an explicit descriptions of the collection $\mathbf{W}$. For $1\leq i \leq n$ and $d\geq 2$, let $e_i=(0,\dots,1,\dots,0)$ denote the vertex in $\Delta_{n-1}$ with 1 in the $i^{\textnormal{th}}$ entry, and let $e_I=\cup_{h=1}^d e_{i_h}$ denote the collection of $d$ vertices in $\Delta_{n-1}$ indexed by $I\in \N_{n,d}$. We further denote by $\delta_I$ the top-dimensional face in $\Delta_{n-1}^{d-1}$ which contains $e_I$. Let $\mathbf{V}$ be the collection of matching spheres in $\Sigma_{n-1,d-1}$, indexed by $\N_{n-1,d-1}$, and assume that $\Sigma_{n-1,d-1}$ is inductively constructed under Definition~\ref{def:Sigma}. For $I\in \N_{n-1,d-1}$, let $\mathfrak{h}_I$ be the handle in $\Sigma_{n-1,d-1}$ whose core is half of the vanishing sphere $V_I$, so that $\Sigma_{n-1,d-1}=X \cup( \cup_I \mathfrak{h}_I)$ (here, the first $\cup$ denotes a handle attachment).

\begin{lemma}\label{prop:mspheresvertices}
    Assume that the core of $\Sigma_{n-1,d-1}$ is the $(d-2)$-dimensional skeleton $\Delta_{n-1}^{d-2}$ of the standard $(n-1)$-simplex, in such a way that $\X_{n,d}$ in Definition~\ref{def:Sigma} retracts (under the negative Liouville flow) to the single vertex $(1,0,\dots,0)$. Moreover, assume that the core of each $W_{1,I+1}=V_I$ coincides with the cycle in $H_{d-2}(\Delta_{n-1}^{d-2})$ passing through $(1,0,\dots,0)$ and $e_{I+1}$. Then for each ${J\in \N_{n,d}\setminus \N_{n-1,d-1}^{+1}}$, the vanishing sphere $W_J$ in Definition~\ref{def:halfmatchingspheres} retracts onto the cycle in $H_{d-2}(\Delta_{n-1}^{d-2})$ passing through $e_{J}$.
\end{lemma}

\begin{proof}
    We use the projection of $p:\Sigma_{n-1,d-1} \to \D$ depicted in Figure~\ref{fig:symplAuslfibre} (and again in Figure~\ref{fig:finalmpaths}, for $n=5$, $d=3$), under which each vanishing sphere $W_{1,I+1}$ (${I\in \N_{n-1,d-1}}$) is represented as a matching path. In this description, the assumption that $W_{1,I+1}$ retracts under the negative Liouville flow onto the cycles passing through $e_{1,I+1}$ is equivalent to the fact that $W_{1,I+1}$ retracts onto a cycle in $H_{d-2}(\Delta_{n-1}^{d-2})$ passing through $(1,0,\dots,0)$ and $\delta_{I+1}$. In particular, $W_{1,I+1}$ retracts onto a cycle in $H_{d-2}(\Delta_{n-1}^{d-2})$ containing the core of $\mathfrak{h}_I$ (and passing through $(1,0,\dots,0)$). Note that in this notation,  $\delta_{I+1}$ coincides with the core of $\mathfrak{h}_I$. We use Definition~\ref{def:halfmatchingspheres} to express the vanishing spheres $\{W_J\}_J$, for $J\in \N_{n,d}$, as matching spheres above certain matching paths. Denote by $\D'\subset p(\X_{n,d})\subset \D$ a sufficiently big disc in $p(\X_{n,d})$, containing all the images of the branch points of $p$ contained in $X_{n,d}$, depicted by a dashed line in Figure~\ref{fig:finalmpaths}.

    \begin{figure}[t]
        \centering
        \begin{tikzpicture}[scale=1.5,v/.style={draw,shape=circle, fill=black, minimum size=.5mm, inner sep=0pt, outer sep=0pt},cross/.style={cross out,draw=black, minimum size=1mm},cross/.style={cross out, draw=black, minimum size=1.5mm, inner sep=0pt, outer sep=0pt},cross/.default={1pt}]
            \draw (0,1.5) ellipse (3.5cm and 1.5cm);
            \draw[dashed] (0,1.5) ellipse (3.2cm and 1.2cm);
    
            \node at (0+3.8,1.5) {$\D$};
            \node at (0+2.9,1.5) {$\D'$};
    
            \node[cross] at (-2,1) {};
            \node[cross] at (0,1) {};
            \node[cross] at (2,1) {};
    
            \draw[color=frenchrose,line width=.9] ($(0,1.5)+(-75:5 and 1.5)+(0,.055)$) to (-2,1);
            \draw[color=orange,line width=.9pt] ($(0,1.5)+(-85:5 and 1.5)$) to (0,1);
            \draw[color=ForestGreen,line width=.9pt] ($(0,1.5)+(-90:5 and 1.5)$) to (2,1);
            \draw ($(0,1.5)+(-105:5 and 1.5)+(0,.05)$) to[out=150,in=-90] (-2.3,1.1) to[out=90,in=180] (-1.9,1.4) to[out=0,in=180] (0,.6) to[out=0,in=200] (2,1);
            \draw ($(0,1.5)+(-98:5 and 1.5)+(0,.015)$) to[out=160,in=-90] (-2.5,1) to[out=90,in=180] (-1.7,1.6) to[out=0,in=140] (0,1);
            \draw ($(0,1.5)+(-115:5 and 1.5)+(0,.165)$) to[out=145,in=-90] (-2.7,1) to[out=90,in=160] (.5,1.6) to[out=-20,in=160] (2,1);
    
            \draw[color=frenchrose,line width=.9] ($(0,1.5)+(-75:5 and 1.5)+(0,.055)$) to ($(0,1.5)+(-75:5 and 1.5)+(0,-1)$);
            \draw[color=orange,line width=.9pt] ($(0,1.5)+(-85:5 and 1.5)$) to ($(0,1.5)+(-85:5 and 1.5)+(0,-1)$);
            \draw[color=ForestGreen, line width=.9pt] ($(0,1.5)+(-90:5 and 1.5)$) to ($(0,1.5)+(-90:5 and 1.5)+(0,-1)$);
            \draw ($(0,1.5)+(-105:5 and 1.5)+(0,.05)$) to  ($(0,1.5)+(-105:5 and 1.5)+(0,-1)$);
            \draw ($(0,1.5)+(-98:5 and 1.5)+(0,.015)$) to ($(0,1.5)+(-98:5 and 1.5)+(0,-1)$) ;
            \draw ($(0,1.5)+(-115:5 and 1.5)+(0,.165)$) to ($(0,1.5)+(-115:5 and 1.5)+(0,-1)$);
            
            \node[cross] at($(0,1.5)+(-75:5 and 1.5)+(0,-1)$) {};
            \node[cross] at ($(0,1.5)+(-85:5 and 1.5)+(0,-1)$) {};
            \node[cross] at ($(0,1.5)+(-90:5 and 1.5)+(0,-1)$) {};
            \node[cross] at ($(0,1.5)+(-98:5 and 1.5)+(0,-1)$){};
            \node[cross] at ($(0,1.5)+(-105:5 and 1.5)+(0,-1)$){};
            \node[cross] at ($(0,1.5)+(-115:5 and 1.5)+(0,-1)$) {};
    
            \draw ($(0,1.5)+(-120:5 and 1.5)+(0,.245)$) to[out=-80,in=100] ($(0,1.5)+(-115:5 and 1.5)+(-.2,-.8)$) to[out=-80,in=180] ($(0,1.5)+(-115:5 and 1.5)+(0,-1.2)$) to[out=0,in=-90]  ($(0,1.5)+(-115:5 and 1.5)+(.2,-.8)$) to[out=90,in=-90] ($(0,1.5)+(-112:5 and 1.5)+(0,.12)$);
    
            \draw ($(0,1.5)+(-107:5 and 1.5)+(0,.07)$) to[out=-90,in=90] ($(0,1.5)+(-105:5 and 1.5)+(-.2,-.8)$) to[out=-90,in=180] ($(0,1.5)+(-105:5 and 1.5)+(0,-1.2)$) to[out=0,in=-90]  ($(0,1.5)+(-105:5 and 1.5)+(.2,-.8)$) to[out=90,in=-90] ($(0,1.5)+(-103:5 and 1.5)+(0,.035)$);
    
            \draw ($(0,1.5)+(-100:5 and 1.5)+(0,.02)$) to[out=-90,in=90] ($(0,1.5)+(-98:5 and 1.5)+(-.2,-.8)$) to[out=-90,in=180] ($(0,1.5)+(-98:5 and 1.5)+(0,-1.2)$) to[out=0,in=-90]  ($(0,1.5)+(-98:5 and 1.5)+(.2,-.8)$) to[out=90,in=-90] ($(0,1.5)+(-96:5 and 1.5)+(0,.005)$);
    
            \draw ($(0,1.5)+(-92:5 and 1.5)$) to[out=-90,in=90] ($(0,1.5)+(-90:5 and 1.5)+(-.2,-.8)$) to[out=-90,in=180] ($(0,1.5)+(-90:5 and 1.5)+(0,-1.2)$) to[out=0,in=-90]  ($(0,1.5)+(-90:5 and 1.5)+(.2,-.8)$) to[out=90,in=-90] ($(0,1.5)+(-88:5 and 1.5)$);
    
            \draw[] ($(0,1.5)+(-87:5 and 1.5)$) to[out=-90,in=90] ($(0,1.5)+(-85:5 and 1.5)+(-.2,-.8)$) to[out=-90,in=180] ($(0,1.5)+(-85:5 and 1.5)+(0,-1.2)$) to[out=0,in=-90]  ($(0,1.5)+(-85:5 and 1.5)+(.2,-.8)$) to[out=90,in=-90] ($(0,1.5)+(-83:5 and 1.5)+(0,.005)$);
    
            \draw[] ($(0,1.5)+(-77:5 and 1.5)+(0,.045)$) to[out=-90,in=90] ($(0,1.5)+(-75:5 and 1.5)+(-.2,-.8)$) to[out=-90,in=180] ($(0,1.5)+(-75:5 and 1.5)+(0,-1.2)$) to[out=0,in=-90]  ($(0,1.5)+(-75:5 and 1.5)+(.2,-.8)$) to[out=90,in=-90] ($(0,1.5)+(-73:5 and 1.5)+(0,.065)$);
    
            \draw[color=red, line width=.9pt] ($(0,1.5)+(-115:5 and 1.5)+(0,-1)$) to ($(0,1.5)+(-119:5 and 1.5)+(0,.22)$)to[out=150,in=-90] (-3,1) to[out=90,in=160] (.5,1.9) to[out=-20,in=140] (2,1);
    
            \draw[color=blue, line width=.9pt] ($(0,1.5)+(-115:5 and 1.5)+(0,-1)$) to ($(0,1.5)+(-118:5 and 1.5)+(0,.19)$) to[out=150,in=-90] (-2.9,1) to[out=90,in=160] (.5,1.7) to[out=-20,in=90] (1,1) to[out=-90,in=80] ($(0,1.5)+(-91:5 and 1.5)$) to ($(0,1.5)+(-90:5 and 1.5)+(0,-1)$);
    
            \draw[color=Plum, line width=.9pt] ($(0,1.5)+(-115:5 and 1.5)+(0,-1)$)to ($(0,1.5)+(-113:5 and 1.5)+(0,.13)$) to[out=20,in=150,looseness=.2] ($(0,1.5)+(-85:5 and 1.5)+(.1,.01)$) to[out=-90,in=90] ($(0,1.5)+(-85:5 and 1.5)+(0,-1)+(.1,0)$) to[out=-90,in=0] ($(0,1.5)+(-85:5 and 1.5)+(0,-1)+(0,-.1)$) to[out=180,in=-90] ($(0,1.5)+(-85:5 and 1.5)+(0,-1)+(-.1,0)$) to[out=90,in=-90] ($(0,1.5)+(-85:5 and 1.5)+(-.1,0)$) to[out=90,in=90,looseness=.7] ($(0,1.5)+(-90:5 and 1.5)+(.1,0)$) to ($(0,1.5)+(-90:5 and 1.5)+(0,-1)$);

            \draw[dotted, color=Plum, line width=1.2pt] ($(0,1.5)+(-115:5 and 1.5)+(0,-1)$)to[out=90,in=225] ($(-2,1)+(0,.25)$) to[out=45,in=90,looseness=.7] ($(0,1.5)+(-85:5 and 1.5)+(.13,.01)$) to[out=-90,in=90] ($(0,1.5)+(-85:5 and 1.5)+(0,-1)+(.13,0)$) to[out=-90,in=0] ($(0,1.5)+(-85:5 and 1.5)+(0,-1)+(0,-.15)$) to[out=180,in=-90] ($(0,1.5)+(-85:5 and 1.5)+(0,-1)+(-.13,0)$) to[out=90,in=-90] ($(0,1.5)+(-85:5 and 1.5)+(-.13,0)$) to[out=90,in=90,looseness=.7] ($(0,1.5)+(-90:5 and 1.5)+(.13,0)$) to ($(0,1.5)+(-90:5 and 1.5)+(0,-1)$);
    
    
        
    \end{tikzpicture}
            \caption{The matching path $p(W_{2,4,5})$ on the base of $\Sigma_{4,2}$.}
            \label{fig:finalmpaths}
        \end{figure}

    For simplicity, we isotope the base of $p$ so that all its critical values contained in $\D'$ have the same imaginary part (they are ``on the same horizontal line'', as in Figure~\ref{fig:finalmpaths}). Fix ${J\in\N_{n,d}\setminus \N_{n-1,d-1}^{+1}}$. The matching paths $p(V_{\widehat{j_2}})$ and $p(V_{\widehat{j_1}})$ intersect transversely at an endpoint. In this setting, the matching path for $\tau_{V_{\widehat{j_2}}}V_{\widehat{j_1}}$ can be explicitly described, and Seidel explains this in \cite[Section~16h]{SeidelBk} (see Figure~16.4 therein). In particular, it is obtained as the ``half'' Dehn twist of $p(V_{\widehat{j_1}})$ around $p(V_{\widehat{j_2}})$, and it starts at the endpoint of $p(V_{\widehat{j_1}})$ not contained in $\D'$, enters $\D'$, goes ``over'' the endpoints of $p(V_{L})$ for all $L<\widehat{j_1}$, exits  $\D'$, and ends at the endpoint of $p(V_{\widehat{j_2}})$ not in $\D'$. See Figure~\ref{fig:finalmpaths}, where we depicted $p(V_{\widehat{j_1}})$, $p(V_{\widehat{j_2}})$, and $p(\tau_{V_{\widehat{j_2}}}V_{\widehat{j_1}})$ in red, green, and blue respectively, for $n=5$, $d=3$, and $J=(2,4,5)$. The corresponding matching sphere retracts under the negative Liouville flow onto a cycle which contains the cores of $\mathfrak{h}_{\widehat{j_2}}$ and $\mathfrak{h}_{\widehat{j_1}}$ (and passes through $(1,0,\dots,0)$).

    We note that $p(\tau_{V_{\widehat{j_2}}}V_{\widehat{j_1}})$ and $p(V_{\widehat{j_3}})$ intersect transversely at an interior point. In this setting, the matching path for $\tau_{V_{\widehat{j_3}}}\tau_{V_{\widehat{j_2}}}V_{\widehat{j_1}}$ can also be explicitly described, which Seidel does in \cite[Section~18c]{SeidelBk} (see Figure~18.2 therein). In particular, $p(\tau_{V_{\widehat{j_3}}}\tau_{V_{\widehat{j_2}}}V_{\widehat{j_1}})$ starts at the endpoint of $p(\tau_{V_{\widehat{j_2}}}V_{\widehat{j_1}})$ not in $\D'$, goes ``over'' the endpoints of $p(V_{L})$ (for all $L<\widehat{j_3}$), exits $\D'$, goes ``under'' the endpoint of $p(V_{\widehat{j_3}})$ not in $\D'$, and ends at the endpoint of $p(V_{\widehat{j_2}})$ without re-entering $\D'$. This is depicted in Figure~\ref{fig:finalmpaths}, where $p(V_{\widehat{j_3}})$ and $p(\tau_{V_{\widehat{j_3}}}\tau_{V_{\widehat{j_2}}}V_{\widehat{j_1}})$ are in pink and dotted purple respectively. The matching path for $\tau_{V_{\widehat{j_3}}}\tau_{V_{\widehat{j_2}}}V_{\widehat{j_1}}$ retracts under the negative Liouville flow onto a cycle containing the cores of $\mathfrak{h}_{\widehat{j_3}}$, $\mathfrak{h}_{\widehat{j_2}}$ and $\mathfrak{h}_{\widehat{j_1}}$, and passing through $(1,0,\dots,0)$. Iterating this process, we can express $\tau_{V_{\widehat{j_k}}}\cdots\tau_{V_{\widehat{j_2}}}V_{\widehat{j_1}}$ as a matching path in $p(\Sigma_{n-1,d-1})$, in such a way that the matching sphere retracts onto a cycle containing the cores of $\mathfrak{h}_{\widehat{j_k}},\dots,\mathfrak{h}_{\widehat{j_1}}$, and passes through $(1,0,\dots,0)$. We note that we can assume the matching paths $p_k:=p(\tau_{V_{\widehat{j_k}}}\cdots\tau_{V_{\widehat{j_2}}}V_{\widehat{j_1}})$ and $p(V_{\widehat{j_{k+1}}})$ to be ``consecutive'' in a sense similar to that in Lemma~\ref{prop:vanishingDehntwists}. More precisely, the matching path $p_k$ starts at the endpoint of $p(V_{\widehat{j_1}})$ not in $\D'$, enters $\D'$, goes ``over'' the endpoints of $p(V_{L})$, for all $L \leq \widehat{j_{k+1}}$, then exits $\D'$ without going ``over'' any other critical value in $\D'$, and does not enter it again. Indeed, for $\widehat{j_{k+1}}< L < \widehat{j_{k}}$, the matching spheres $\tau_{V_{\widehat{j_k}}}\cdots\tau_{V_{\widehat{j_2}}}V_{\widehat{j_1}}$ and $V_L$ are disjoint, so the matching path of the former can cross over the endpoints of the latter without modifying the matching sphere. We depicted this in Figure~\ref{fig:finalmpaths}: the matching sphere $\tau_{V_{\widehat{j_3}}}\tau_{V_{\widehat{j_2}}}V_{\widehat{j_1}}$ (dotted purple) is disjoint from $V_{12}$ (pink), and $(1,2)<\widehat{j_3}$. The final matching path $\tau_{V_{\widehat{j_3}}}\tau_{V_{\widehat{j_2}}}V_{\widehat{j_1}}$ is the continuous purple curve in figure. Finally, for $k=d$, the matching sphere $W_J=\tau_{V_{\widehat{j_k}}}\cdots\tau_{V_{\widehat{j_2}}}V_{\widehat{j_1}}$ retracts under the negative Liouville flow onto a cycle containing the cores of $\mathfrak{h}_{\widehat{j_d}},\dots,\mathfrak{h}_{\widehat{j_1}}$. After crossing over any endpoint of disjoint matching spheres, $p(W_J)$ does not enter into $\D'$, so $W_J$ does not retract to a cycle passing through $(1,0,\dots,0)$.
\end{proof}

\begin{proof}[Proof of Proposition~\ref{prop:simplex}]
    We argue by induction on $d$. The Milnor fibre of $\Phi_{n,1}:=\f_{n,1}$, and in particular its core, consists of $n+1$ points, which is also the $0$-dimensional skeleton of the standard $n$-simplex. The case $d=2$ follows from \cite[Proposition~3.18]{DiDedda23} (see \cite[Section~1.1]{DiDedda23} for the \emph{a posteriori} description of $\Sigma_{n,2}$, $\Sigma_n$ in text). Indeed, the Milnor fibre $\Sigma_{n,2}$ of $\Phi_{n,2}=\f_{n,2}$ is homeomorphic to a Riemann surface $X_n \# Y_n$, where $X_n$ is a contractible surface, $Y_n$ is the thickening of the complete graph  associated to the set of $n-1$ vertices, glued to $X_n$ in a prescribed way. The Liouville flow is outward pointing along the boundary of $X_n \# Y_n$. The points of $X_n$ which do not escape to infinity under the Liouville flow are given by the $1$-dimensional skeleton of the $n$-simplex with one common vertex (say, for example, $(1,0,\dots,0)$). The core of $Y_n$ coincides with the cores of the handles, i.e.\ with the complete graph associated to the set of $n-1$ vertices.

    Suppose inductively that the core of $\Sigma_{n,d}$ in Setup~\ref{setup:sympdata} is $\Delta_{n}^{d-1}$. More precisely, we require that the vanishing cycles $\{V_I=W_{1,I+1} \mid I\in\N_{n,d}\}$ retract under the negative Liouville flow onto the cycles in $H_{d-1}(\Delta_{n}^{d-1})$ passing through $(1,0,\dots,0)$ and exactly $d$ other vertices. By Lemma~\ref{prop:mspheresvertices}, the collection $\mathbf{W}$ of vanishing cycles indexed by $\N_{n+1,d+1}$ retracts onto a collection of cycles in $H_{d-1}(\Delta_{n}^{d-1})$ passing through all possible subsets of size $d$ of vertices. The core of $\Sigma_{n+1,d+1}$ is the union of the point $(1,0,\dots,0)$ and the \linebreak $d$-dimensional core of the ${n+1} \choose {d+1}$ $(d+1)$-dimensional $d$-handles (attached along push-offs of $\mathbf{W}$).
  \end{proof}

\subsection{Two Lefschetz fibrations}\label{sec:appendix}
In this article, we provided two different geometric presentations of higher Auslander algebras of type A. Namely, Corollary~\ref{cor:maincor} provides a quasi-equivalence of triangulated $A_{\infty}$-categories $\F(\f_{n,d}) \simeq \perf(\A_{n,d})$, while Section~\ref{section:3} gives an explicit way to construct an abstract Lefschetz fibration ${\Phi_{n,d}:\X_{n,d} \to \D}$ such that there is a quasi-equivalence of triangulated $A_{\infty}$-categories $\F(\Phi_{n,d}) \simeq \perf(\A_{n,d})$. The following is an immediate consequence.
\begin{corollary}\label{cor:appendix}
    There is a quasi-equivalence of triangulated $A_{\infty}$-categories
\begin{equation*}
    \pushQED{\qed} \F(\f_{n,d}) \xrightarrow{\text{ $\simeq$ }} \F(\Phi_{n,d}).\qedhere \popQED
\end{equation*}
\end{corollary}
Corollary~\ref{cor:appendix} identifies the families of singularities studied in Section~\ref{section:2} and the families of Lefschetz fibrations constructed in Section~\ref{section:3} at the level of their Fukaya--Seidel categories. We conjecture a stronger identification, namely that the two families of maps are isomorphic in the sense of \cite[Section~15b]{SeidelBk}.

\begin{conjecture}\label{conj:appendix}
    There are an isomorphism between $\C^d$ and $(\X_{n,d})^{\circ}$ and an automorphism of $\D$, as exact symplectic manifolds, which are compatible with the symplectic Lefschetz fibrations $\w_{n,d}$ and the completion of $\Phi_{n,d}$.
\end{conjecture}

Proving Conjecture~\ref{conj:appendix} would amount to understanding the topology and the structure (as a symplectic hypersurface of the boundary contact manifold) of the regular fibre $\M_{n,d}$ of $\w_{n,d}$ (equipped with an appropriate collection of vanishing cycles) to a degree that is beyond the scope of the article. Furthermore, this conjectural identification would not substantially improve the current state of Theorem~\ref{thm:introthm4}. Indeed, the key feature of the Lefschetz fibrations $\Phi_{n,d}$ constructed in Section~\ref{section:3} is that they are inductively constructed, much like the higher Auslander algebras $\A_{n,d}$, not that they are defined by an explicit polynomial (which the family $\f_{n,d}$ is). For this reason, Theorem~\ref{thm:introthm4} reflects the inductive nature of Iyama's higher Auslander correspondence. Nonetheless, we prove the following as a partial result towards Conjecture~\ref{conj:appendix}.

\begin{proposition}\label{prop:identificationcores}
    If the Liouville completion of $\X_{n,d}$ in Setup~\ref{setup:sympdata} is Weinstein deformation equivalent to $\C^d$, then the completion of $\X_{n+1,d+1}$ constructed in Theorem~\ref{prop:uniquetopfibr} is Weinstein deformation equivalent to $\C^{d+1}$. Moreover, we can equip $\M_{n+1,d+1}$ with a Liouville form with respect to which $\Core(\M_{n+1,d+1})$ comprises of the $d$-dimensional skeleton of the $(n+1)$-simplex, together with a negligible component, which we call \emph{fins} (see Definition~\ref{def:fins} below). In particular, we have the following relation between $\Core(\M_{n+1,d+1})$ and the core of the Liouville completion of $\Sigma_{n+1,d+1}$:
    \begin{equation}\label{eq:coresrelation}
       ( \Sigma_{n+1,d+1})^{\circ} \hookleftarrow \Core(\Sigma_{n+1,d+1}) = \Delta_{n+1}^d \subset \Core(\M_{n+1,d+1}) \hookrightarrow \M_{n+1,d+1}.
    \end{equation}
\end{proposition}

\begin{remark}
    We remark below (Remark~\ref{rk:smoothingofcurrents}) why we expect the fins of $\Core(\M_{n+1,d+1})$ in Proposition~\ref{prop:identificationcores} to be contracted, up to an appropriate small perturbation of the Liouville form along the diagonal locus of $\Sym^d(\C)$. In particular, the ``$\subset$'' sign in \eqref{eq:coresrelation} should be replaced by ``$=$''. As it is, one should think of Proposition~\ref{prop:identificationcores} as a heuristic justification for Conjecture~\ref{conj:appendix}, more so than having significance in isolation. We also remark that Proposition~\ref{prop:identificationcores} is a rather weak evidence for Conjecture~\ref{conj:appendix}, as it only compares the cores of $F_{n+1,d+1}$ and $\Sigma_{n+1,d+1}$ as topological spaces, disregarding the symplectic structure of the relevant objects.
\end{remark}

In order to prove Proposition~\ref{prop:identificationcores}, we first provide a description of $\Core(\M_{n+1,d+1})$. Fix $\alpha=dx \wedge dy$ to be the standard area form on $\C$, associated to the Liouville form $\theta_{\C}=\frac{1}{2}(xdy-ydx)$, whose vector field is radial and outwards pointing. We equip $\Sym^d(\D)$ with an exact symplectic structure $\omega$ prescribed in \cite{Perutz} (recall that this is equivalent to the standard symplectic structure on $\C^d\cong \Sym^d(\D)$). Fix $\theta:=\theta_{\Sym^d(\C)}$ the associated Liouville form, which coincides with the smooth pushforward $\bpi_*((\theta_{\C})^{\times d})$ of the Liouville form  $\theta_{\C}$ on $\C$ away from the diagonal. We consider $\w_{n,d}$ as a symplectic Lefschetz fibration on $\Sym^d(\D)$, and we denote by $\M_{n,d}$ the regular fibre of $\w_{n,d}$. We define the following.

\begin{definition}
    For $h=1,\dots,n+1$, let $\xi^h$ denote the $(n+1)^{\textnormal{th}}$ roots of unity, with $\xi:=e^{\frac{2i\pi}{n+1}}$. For $n \geq 1$, we define $\HH_n$ to be the union of the straight segments from the origin to $\{\xi^h\}_{h}$, seen as a subset of the complex plane $\C$. For $n\geq d \geq 1$, we define $\Sym^d(\HH_n)$ to be the unordered tuple of $d$ points on $\HH_n$.
\end{definition}

\begin{notation*}
    We write $t\xi^h$ ($t\in [0,1]$) for a point on the straight segment from $0\in\C$ to $\xi^h$. A point $\mathbf{x}$ on $\Sym^d(\HH_n)$ is of the form $\{t_{i_1}\xi^{i_1},\dots,t_{i_d}\xi^{i_d}\}$, for some ${(t_{i_1},\dots,t_{i_d})\in [0,1]^d}$. Note the slight ambiguity of this notation: if some coordinate $x_h$ of $\mathbf{x}$ vanishes, we can choose \emph{any} root of unity $\xi^{i_h}$ and set $x_h= t_{i_h}\xi^{i_h}$, with $t_{i_h}=0$. This slight ambiguity is only notational, and will not affect the constructions in this section.
\end{notation*}

\begin{definition}\label{def:fins}
    We say that a point on $\Sym^d(\HH_n)$ has \emph{all coordinates on different branches} if for each pair $(h,k)$ satisfying $t_{i_h},t_{i_k}\neq 0$, $i_h\neq i_k$. We denote the set of such point by $\Sym^d(\HH_n)_{\star}$ We define the \emph{fins} of $\Sym^d(\HH_n)$ as the set of points whose coordinates are not all on different branches.
\end{definition}

We use the following description of $\Sym^d(\HH_n)_{\star}$.

\begin{lemma}\label{lem:embeddingconesimplex}
    The subset of points on $\Sym^d(\HH_n)$ which have all coordinates on different branches embeds into $(\CR_{\geq 0})^{n+1}$ under the map $\iota: \Sym^d(\HH_n)_{\star} \hookrightarrow (\CR_{\geq 0})^{n+1}$ defined by
    \begin{equation*}
        \iota(\{t_{i_1}\xi^{i_1},\dots,t_{i_d}\xi^{i_d}\})=(0,\dots,t_{i_1},\dots,t_{i_d},\dots,0),
    \end{equation*}
    where each $t_{i_h}$ is in the $i_h^{\textnormal{th}}$ coordinate. In particular, its image under $\iota$ is diffeomorphic to $\Cone(\Delta_n^{d-1})$.
\end{lemma}

\begin{proof}
    Since $t_{i_h}\in[0,1]$ for each $h$, the image of $\Sym^d(\HH_n)_{\star}$ under the continuous map $\iota$ is given exactly by the union of the $(d-1)$-dimensional faces of the $(n+1)$-dimensional cube $[0,1]^{n+1}$. Moreover, $\iota$ is bijective on its image, and the inverse map is given by
    \begin{equation*}
       (0,\dots,y_1,\dots,y_d,\dots,0) \mapsto \{y_1\xi^{j_1},\dots, y_d\xi^{j_d}\}
    \end{equation*}
    if each $y_h$ is in the $j_h^{\textnormal{th}}$ coordinate. The set of $(d-1)$-dimensional faces of the ${(n+1)}$-dimensional cube is itself diffeomorphic to $\Cone(\Delta_n^{d-1})$, under (the restriction of) a piecewise linear function, shrinking the standard cube into $\Cone(\Delta_n^{n})$.
\end{proof}

\begin{proposition}\label{prop:relativecoreMilnor}
    Up to a deformation of the Liouville form in the interior of $\w_{n,d}^{-1}(\D_d)$, the relative core of the Liouville pair $(\Sym^d(\C),\M_{n,d})$ strictly contains $\Cone(\Delta_{n}^{d-1})$.
\end{proposition}

\begin{proof}
    We first describe the relative core of the Liouville pair $(\C,\M_{n,1})$. Fix the Liouville form $\theta_{\C}=\frac{1}{2}(xdy-ydx)$ on $\C$, whose negative flow is radial and inwards pointing. We consider an isotopy of a sufficiently big disc in $\C$, containing all critical points of $\w_{n,1}$, which rotates each $c_h$ (as in Definition~\ref{def:BFShatn1}) and identifies it with $\xi^h$. The isotopy is similar to the one defined in the proof of Proposition~\ref{prop:Gabungraded}, depicted in Figure~\ref{fig:isotopybase}. The negative gradient flow under this rotation is still radial and inwards pointing. The relative core of $(\C,\M_{n,1})$ with respect to the isotoped Liouville form is given by 
    \begin{equation*}
        \Core(\C,\M_{n,1})=\{0\} \cup (\{\xi^{1},\dots,\xi^{n+1}\}\times \CR) = (\HH_n)^{\circ},
    \end{equation*}
    where $(\HH_n)^{\circ}$ denotes the Liouville completion of $\HH_n$, obtained by ``stretching'' each segment $\{t\xi^h\}_t$ in the positive direction.

    Fix a Liouville form $\theta$ and an exact symplectic form $\omega$ on $\Sym^d(\C)$, prescribed by \cite{Perutz} with respect to the choice of isotoped 1-form on $\C$ above. The negative gradient flow of the 2-current (in the sense of \cite[Section~7]{Perutz}) $\bpi_*(\alpha^{\times d})$ points inwards with respect to each variable. This allows us to identify the relative core of $(\Sym^d(\C),\M_{n,d})$ with respect to $\bpi_*(\alpha^{\times d})$ with the image under $\bpi$ of $(\Core(\C,\M_{n,1}))^{\times d}$. In particular, the relative core of the Liouville pair $(\Sym^d(\C),\M_{n,d})$ with respect to $\bpi_*(\alpha^{\times d})$ coincides with $\bpi((\Core(\C,\M_{n,1}))^{\times d})=\Sym^d(\HH_n)$.
    
    Finally, the smoothing of $\bpi_*(\alpha^{\times d})$ coincides with $\bpi_*(\alpha^{\times d})$ away from the diagonal locus. In particular, it preserves the Liouville core away from this locus. By Lemma~\ref{lem:embeddingconesimplex}, this coincides with $\Cone(\Delta_{n}^{d-1})$.
\end{proof}

\begin{remark}\label{rk:smoothingofcurrents}
    The exact form $\omega$ is prescribed by Perutz in \cite[Section~7]{Perutz} as the ``global smoothing'' of $\bpi_*(\alpha^{\times d})$, using the smoothing theory for currents. It can be constructed by perturbing the Liouville vector field $Z_{\theta}$ associated to $\bpi_*(\alpha^{\times d})$ in a neighbourhood of the diagonal. More precisely, \cite[Lemma~7.4]{Perutz} prescribes a way to construct a K\"ahler form from a collection of continuous plurisubharmonic functions $\{H_i\}_i$ defined on an open cover $\{U_i\}_i$ of $\Sym^d(\C)$. The latter define the pushforward $\bpi_*(\alpha^{\times d})$ by $\bpi_*(\alpha^{\times d})|_{U_i}=dd^{c}H_i$. The K\"ahler form is given by a collection of smooth plurisubharmonic functions $\{H_i'\}_i$ defined on a refinement $\{V_{i,j}\subset U_i\}_{i,j}$ of the open cover, where each function is given by $H_{i,j}'=H_{i,j}+\chi|_{V_{i,j}}$ ($H_{i,j}:=H_{i}|_{V_{i,j}}$) for some continuous function $\chi:\Sym^d(\C) \to \CR$ which is supported on a neighbourhood of the diagonal locus. The K\"ahler form constructed by Perutz is
    \begin{equation*}
        \omega|_{V_{i,j}}:= dd^{c}(\bpi_*(H_{i,j}'))=\bpi_*(dd^{c}(H_{i,j}))+dd^{c}\chi.
    \end{equation*}
    In particular, the Liouville 1-forms associated to $\bpi_*(\alpha^{\times d})$ and $\omega$ are given by $\theta=d^{c}H_{i,j}$ and $\theta=d^{c}H_{i,j}'$ respectively. The Liouville vector fields differ by the vector field of $d^{c}(\chi|_{V_{i,j}})$ on each chart of the open cover intersected with the neighbourhood of the diagonal. We expect that one should be able to pick an appropriate function $\chi$, for which the effect of $d^{c}(\chi|_{V_{i,j}})$ on each component is that of a small perturbation of the product Liouville flow near each $\xi^{i_h}$. The effect on the fins of $\Sym^d(\HH_n)$ should be following. Fix some $j\in\{1,\dots,n+1\}$ and some point $\{t_{i_1}\xi^{j},t_{i_2}\xi^{j},\ldots\}$ in a fin of $\Sym^d(\HH_n)$. After perturbing the Liouville vector field in a neighbourhood of the diagonal, the entire fin should be contracted (under the negative Liouville flow) to the point $\{0,0,\ldots\}$. The upshot is that, after smoothing the 2-current in a neighbourhood of the diagonal locus, we expect the relative core $\Core(\Sym^d(\C),\M_{n,d})$ to coincide with $\Sym^d(\HH_n)_{\star}$, i.e.\ with $\Cone(\Delta_{n}^{d-1})$ by Lemma~\ref{lem:embeddingconesimplex}.
\end{remark}

In particular, we have the following.

\begin{corollary}\label{cor:coreMilnor}
    Up to a deformation of the Liouville form $\theta$ in the interior of $\w_{n,d}^{-1}(\D)$, the core of the regular fibre $\M_{n,d}$ of $\w_{n,d}$ strictly contains $\Delta_{n}^{d-1}$.
\end{corollary}

\begin{proof}
    This follows from the identification between the relative skeleton of a Liouville pair $(X^{\circ},F^{\circ})$ and the cone of the skeleton of $F$. Since the relative skeleton of $(\Sym^d(\C),\M_{n,d})$ is $\Cone(\Delta_{n}^{d-1})$ by Proposition~\ref{prop:relativecoreMilnor}, the claim follows.
\end{proof}

\begin{proof}[Proof of Proposition~\ref{prop:identificationcores}]
    We refer to \cite[Chapter~4]{GompfStipsicz} for the formal definitions of $n$-dimensional Weinstein $k$-handles, as well as the related terminology we use in this proof.

    Assume that the Liouville completion $\X_{n,d}^{\circ}$ of $\X_{n,d}$ in Setup~\ref{setup:sympdata} is Weinstein deformation equivalent to $\C^d$, i.e.\ we assume that $\X_{n,d}$ is Weinstein deformation equivalent to a $2d$-dimensional ball $\D^{d}$. We recall from Definition~\ref{def:Sigma} and Theorem~\ref{prop:uniquetopfibr} that $\X_{n+1,d+1}$ is obtained from $\X_{n,d}$ in the following way:
    \begin{itemize}
        \item first attach a $2d$-dimensional Weinstein $d$-handle $\mathfrak{h}_I$ (whose core is $D_I^{-}$) is to $\X_{n,d}$ for each attaching sphere $W_I^{+}=\partial D_I^{+}$, where $I$ ranges over $\I_{n+1,d+1}$ (this gives $\Sigma_{n+1,d+1}=\D^{d} \cup \bigcup_I \mathfrak{h}_I$, where $\cup$ denotes a handle attachment);
        \item subsequently attach a $(2d+2)$-dimensional Weinstein $(d+1)$-handle $\mathfrak{z}_I$ to ${\Sigma_{n+1,d+1}\times \D}$, for each vanishing cycle $\Upsilon_I=D_I^{-} \cup_{W_I^{+}} D_I^{+}$, where $I$ ranges over $\I_{n+1,d+1}$ (this gives $\X_{n+1,d+1}=\Sigma_{n+1,d+1} \cup \bigcup_I \mathfrak{z}_I$).
    \end{itemize}
    We note that if $\mathfrak{h}_I$ is a $2d$-dimensional $d$-handle, then $\mathfrak{h}_I \times \D$ is a $(2d+2)$-dimensional \linebreak $d$-handle, whose cocore is the product of the cocore of $\mathfrak{h}_I$ and $\D$. In this case, we can think of $\mathfrak{h}_I$ as $\mathfrak{h}_I=\mathfrak{h}_I\times \{1\} \subset \mathfrak{h}_I\times \D$. 
    We argue that in $\X_{n+1,d+1}$, each $(\mathfrak{h}_I \times \D,\mathfrak{z}_I)$ form a cancelling pair of Weinstein handles, in the sense of \cite[Chapter~4]{GompfStipsicz}. A sufficient condition for this is that the belt sphere of $\mathfrak{h}_I \times \D$ (i.e.\ the boundary of the cocore of $\mathfrak{h}_I \times \D$) intersects transversely and in a single point the attaching sphere of $\mathfrak{z}_I$ \cite[Proposition~4.2.9]{GompfStipsicz}. In our case, the attaching sphere of $\mathfrak{z}_I$ is the vanishing cycle $\Upsilon_I$, while the belt sphere of $\mathfrak{h}_I \times \D$ is given by
    \begin{equation*}
        \textnormal{Belt sphere}(\mathfrak{h}_I \times \D)=  \left(\textnormal{Belt sphere}(\mathfrak{h}_I)\times \D \right)\cup \left(\textnormal{cocore}(\mathfrak{h}_I)\times \partial \D\right),
    \end{equation*}
    where $\textnormal{cocore}(\mathfrak{h}_I)$ denotes the cocore of a Weinstein handle. The $d$-handle $\mathfrak{h}_I$ is disjoint from $D^{+}_I$, so $\Upsilon_I$ can only intersect the belt sphere of $\mathfrak{h}_I \times \D$ at the points in which $D^{-}_I$ does. Since the belt sphere of $\mathfrak{h}_I$ is the boundary of its cocore, and $D^{-}_I$ is its core, these two cannot intersect. Finally, $D^{-}_I$ intersects $(\textnormal{cocore}(\mathfrak{h}_I)\times \partial \D)$ at exactly the one point $\{0\}\times \{1\}$, where $\{0\}$ denotes the transverse intersection between $\textnormal{cocore}(\mathfrak{h}_I)$ and $D^{-}_I$ (the ``centre'' of the $d$-handle $\mathfrak{h}_I$), and $\{1\}$ denotes the point $1\in\partial\D$ where we identify $\mathfrak{h}_I=\mathfrak{h}_I\times \{1\}$. In particular,
    \begin{equation*}
        \X_{n+1,d+1} = \left(\left(\D^{d} \cup \bigcup_I \mathfrak{h}_I\right)\times \D\right)\cup \bigcup_I \mathfrak{z}_I = \left(\D^d \times \D\right) \cup \left(\bigcup_I (\mathfrak{h}_I \times \D)\right)\cup\bigcup_I \mathfrak{z}_I,
    \end{equation*}
    $\X_{n+1,d+1}$ is Weinstein deformation equivalent to $\D^d \times \D\cong \D^{d+1}$ (a ball of real dimension $2d+2$), and its Liouville completion is the Weinstein deformation equivalent to $\C^{d+1}$. The second part of the claim follows from Corollary~\ref{cor:coreMilnor} and Proposition~\ref{prop:simplex}.
\end{proof}

\begin{remark}
    Assume now that the identification between $\Core(\M_{n,d})$ and $\Delta_n^{d-1}$ following from Remark~\ref{rk:smoothingofcurrents} holds. More precisely, assume that the identification between $\Core(\M_{n,d})$ and $\Core(\Sigma_{n,d})$ in Proposition~\ref{prop:identificationcores} is compatible with the contact structure of the respective boundary contact manifolds. In parallel, $\Core(\Lambda_{n}^{(d)})$ can be identified with $\Delta_n^{d-1}$ in a straightforward way. In particular, the pairs $(\X_{n,d},\Core(\M_{n,d}))$ and $(\Sym^d(\D), \Lambda_{n}^{(d)})$ would define equivalent Liouville sectors, and we would have a quasi-equivalence of pre-triangulated $A_{\infty}$-categories
    \begin{equation}\label{eq:GPSequiv}
       \W(\X_{n,d},\Core(\M_{n,d})) \xrightarrow{\text{ $\simeq$ }} \W(\Sym^d(\D), \Lambda_{n}^{(d)})
    \end{equation}
    by \cite[Corollary~3.9]{GPS2}, where the right-hand side is Sylvan and Ganatra--Pardon--Shende's partially wrapped Fukaya category. Since the left-hand side of \eqref{eq:GPSequiv} is equivalent to Ganatra--Pardon--Shende's definition of the Fukaya--Seidel category $\textnormal{FS}(\f_{n,d})$ by \cite[Proposition~8.20]{GPS2}, we would have a quasi-equivalence of pre-triangulated $A_{\infty}$-categories
    \begin{equation}\label{eq:GPSFSequiv}
        \textnormal{FS}(\f_{n,d})  \xrightarrow{\text{ $\simeq$ }} \W(\Sym^d(\D), \Lambda_{n}^{(d)}).
    \end{equation}
    The quasi-equivalence \eqref{eq:GPSFSequiv} could replace our Theorem~\ref{thm:introthm3} (again assuming the identification in Remark~\ref{rk:smoothingofcurrents}). This identification geometrically justifies the algebraic proof of Theorem~\ref{thm:introthm3}, as it follows from an equivalence of Liouville sectors. Moreover, if we assume the identification
    \begin{equation*}
        \F(\f_{n,d}) \simeq  \textnormal{FS}(\f_{n,d})
    \end{equation*}
    (generally expected, but not presently found in the literature), this would provide further geometric justification for the identification
    \begin{equation*}
        \W^{\textnormal{Aur}}(\Sym^d(\D), \Lambda_{n}^{(d)})\simeq \W(\Sym^d(\D), \Lambda_{n}^{(d)})
    \end{equation*}
    of the two named partially wrapped Fukaya categories.
\end{remark}

  \section*{Acknowledgments}
  I am thankful to my PhD supervisor Yank{\i} Lekili for suggesting this project and for his guidance throughout. This project would not have existed in its present form without the generous discussions with Ailsa Keating (who suggested the approach to the proof of Theorem~\ref{thm:introthm2}), with Tobias Dyckerhoff (who first noted the relationship between simplexes and the Milnor fibre of the singularity $\f_{n,2}$ appearing in \cite{DiDedda23}), and with Gustavo Jasso (about AR theory). I am also thankful for the useful conversations with Nick Sheridan, Matthew Habermann, Jeff Hicks, Daniil Mamaev, Noah Porcelli, and Danil Ko{\v{z}}evnikov. I am very grateful to my thesis examiners N.\ Sheridan and G.\ Jasso, and to the anonymous referees, whose comments greatly improved the content and exposition of this article, especially concerning Sections~\ref{section:quotientFC} and \ref{sec:appendix} (part of Proposition~\ref{prop:identificationcores} was suggested by one of the referees).

  \vskip .25cm
  
  This work was supported by the Faculty of Natural, Mathematical \& Engineering Sciences [NMESFS], King's College London, during its development, and by an ERC Starting Grant (award n.~850713 - HMS) and a Simons Investigator Award (n.~929034) during its conclusion.

\addcontentsline{toc}{section}{References}
{\small \bibliographystyle{alpha}
\bibliography{refs}}
\end{document}